# CONCENTRATION INEQUALITIES AND ASYMPTOTIC RESULTS FOR RATIO TYPE EMPIRICAL PROCESSES

By Evarist Giné[1] and Vladimir Koltchinskii[2]

*University of Connecticut and University of New Mexico*

Let $\mathcal{F}$ be a class of measurable functions on a measurable space $(S, \mathcal{S})$ with values in $[0, 1]$ and let

$$P_n = n^{-1} \sum_{i=1}^{n} \delta_{X_i}$$

be the empirical measure based on an i.i.d. sample $(X_1, \ldots, X_n)$ from a probability distribution $P$ on $(S, \mathcal{S})$. We study the behavior of suprema of the following type:

$$\sup_{r_n < \sigma_P f \leq \delta_n} \frac{|P_n f - P f|}{\phi(\sigma_P f)},$$

where $\sigma_P f \geq \operatorname{Var}_P^{1/2} f$ and $\phi$ is a continuous, strictly increasing function with $\phi(0) = 0$. Using Talagrand's concentration inequality for empirical processes, we establish concentration inequalities for such suprema and use them to derive several results about their asymptotic behavior, expressing the conditions in terms of expectations of localized suprema of empirical processes. We also prove new bounds for expected values of sup-norms of empirical processes in terms of the largest $\sigma_P f$ and the $L_2(P)$ norm of the envelope of the function class, which are especially suited for estimating localized suprema. With this technique, we extend to function classes most of the known results on ratio type suprema of empirical processes, including some of Alexander's results for VC classes of sets. We also consider applications of these results to several important problems in nonparametric statistics and in learning theory (including general excess risk bounds in empirical risk minimization and their versions for $L_2$-regression

Received September 2004; revised May 2005.
[1]Supported in part by NSA Grant H98230-1-0075.
[2]Supported in part by NSF Grant DMS-03-04861.
*AMS 2000 subject classifications.* Primary 60E15; secondary 60F17, 60F15, 62G08, 68T10.
*Key words and phrases.* Ratio type empirical processes, concentration inequalities, ratio limit theorems, localized sup-norms, weighted central limit theorems, VC type classes, moment bounds for empirical processes, nonparametric regression, classification.

and classification and ratio type bounds for margin distributions in classification).

**1. Introduction.** Let $\mathcal{F}$ be a class of measurable functions defined on a measurable space $(S, \mathcal{S})$ and taking values in $[0,1]$, let $X, X_1, \ldots, X_n, \ldots$, be a sequence of independent identically distributed $S$-valued random variables with probability law $P$ and let

$$P_n = n^{-1} \sum_{i=1}^n \delta_{X_i}$$

be the empirical measure based on the sequence $X_i$, that, as usual, we consider as a process on $\mathcal{F}$. Here is some notation that will be used throughout:

$$\sigma_P f \geq \mathrm{Var}_P^{1/2} f, \qquad \sigma_P := \sup_{f \in \mathcal{F}} \sigma_P f$$

(usually $\sigma_P f$ will either be the standard deviation of $f$ or $\sqrt{Pf}$). Given a continuous, strictly increasing function $\phi$ with $\phi(0) = 0$, we are interested in the behavior of suprema of the following type:

$$\sup_{r_n < \sigma_P f \leq \delta_n} \frac{|P_n f - Pf|}{\phi(\sigma_P f)}$$

for some sequences $r_n, \delta_n$. In particular, for given $r_n$ and $\delta_n$, we would like to determine a normalizing sequence $\beta_n$ such that

$$\frac{1}{\beta_n} \sup_{r_n < \sigma_P f \leq \delta_n} \frac{|P_n f - Pf|}{\phi(\sigma_P f)}$$

remains bounded or converges to a constant in probability or almost surely. We are also interested in conditions under which the sequence of stochastic processes

$$\frac{n^{1/2}|P_n f - Pf|}{\phi(\sigma_P f)} I(\sigma_P f \geq r_n), \qquad f \in \mathcal{F},$$

converges in distribution to a Gaussian process indexed by $f \in \mathcal{F}$. Such stochastic processes are often called *normalized* or *ratio type* empirical processes and the distributional convergence results are *weighted central limit theorems* for empirical processes. The study of these processes has a long history that goes back to the 1970s and 1980s when the classical case of $\mathcal{F} := \{I_{(-\infty, t]} : t \in \mathbf{R}\}$ was explored in great detail and definitive answers to most of the questions about the classical ratio type empirical processes were given; see, for example, [48]. In the late 1980s, Alexander, in a series of papers [1, 2, 3], extended this theory to ratio type empirical processes indexed by VC classes of sets $\mathcal{C}$ (i.e., for $\mathcal{F} := \{I_C : C \in \mathcal{C}\}$). He discovered that in this case the crucial role is played by the following functional characteristic of the class:

$$g_{\mathcal{C}}(\delta) := \frac{P(\bigcup_{C \in \mathcal{C}, P(C) \leq \delta} C)}{\delta} \vee 1,$$

which he called *the capacity function* of $\mathcal{C}$. This function is involved in rather sharp and subtle exponential inequalities for empirical processes indexed by VC classes proved by Alexander. The behavior of the capacity function as $\delta \to 0$ happened to be closely related to the asymptotic behavior of ratio type suprema of empirical processes. In recent years, there has been a great deal of work on the development of ratio type inequalities, primarily, in more specialized contexts of nonparametric statistics (see [31, 46]) and learning theory (see [5, 6, 7, 9, 11, 26, 27, 28, 33, 37, 38, 39], etc.). These inequalities have become one of the important ingredients in determining asymptotically sharp convergence rates in regression, classification and other nonparametric problems and they proved to be crucial in bounding the generalization error of learning algorithms based on empirical risk minimization.

In this paper, building upon our earlier work with Jon Wellner [22], we are trying to develop further a general methodology for proving exponential bounds and exploring asymptotics of ratio type empirical processes. This methodology is based on the deservedly famous Talagrand's concentration inequality [43] and on the simple idea of splitting the class $\mathcal{F}$ into slices consisting of functions for which the values of $\phi(\sigma_P f)$ are roughly the same. The empirical process on each slice is compared with its expectation using Talagrand's inequality and then all the pieces are put together using the union bound. This simple approach, called *the method of slicing* or *peeling*, proved to be rather successful in statistical applications (as in [9] or [31]) and it also allows us to obtain a number of sharp results on asymptotics of ratio type suprema (including weighted CLTs), essentially as a straightforward corollary of Talagrand's inequality. The conditions of these limit theorems are expressed in terms of expected values of localized suprema of empirical processes (suprema over the slices). To translate these conditions into a more convenient language for special function classes $\mathcal{F}$ one has to develop bounds on expected localized suprema. We prove such bounds (both upper and lower) under some conditions on random entropies of the class. Unlike most previously known bounds, the new bounds involve the $L_2$ norm of the measurable envelope of the class $\mathcal{F}$, which in applications to ratio limit theorems become the envelopes of the slices. These localized envelopes play about the same role in our theory as Alexander's capacity function plays for classes of sets (and, moreover, in the case when $\mathcal{F}$ is a class of indicators of sets the conditions on localized envelopes can be reformulated as conditions on the capacity function). We are trying to explore in this paper both the power and the limitations of the approach based on slicing and on Talagrand's inequality, and to this end we provide some examples showing in which cases the conditions we obtain are sharp. Our main goal is to provide a link (that seemed to be missing) between the results for classical empirical processes of the 1970s and 1980s extended by Alexander to VC classes of sets and more recent results on ratios developed primarily in learning theory and based on Talagrand's concentration inequality. At the moment, our method allows us to generalize a number of Alexander's theorems to classes of functions, but some other theorems and exponential inequalities seem to be beyond the reach of our approach. On the other hand, most of his specific corollaries for classical empirical processes in $\mathbf{R}^d$ can be obtained by a slight modification of our method, consisting in further decomposing each slice corresponding to a small variance into a relatively small number of subclasses with envelopes which are considerably smaller than that of the full slice. A bit surprisingly, the classical case of classes of sets of small entropy (which are needed to study the standard empirical processes) are harder to handle using Talagrand's inequality and general expectation bounds than the much more massive function classes commonly used in learning theory and nonparametric statistics. In part, this is related to the fact that the Poisson tail parts of the exponential inequalities play a more important role in this case, leading to more complicated asymptotic properties.

Finally, we provide several applications of ratio type empirical processes. First of all, we derive in a much shorter way recent results of the second named author [25] on empirical margin distributions motivated by the needs of learning theory, specifically, the analysis of large margin classifiers. Second, we give general ratio type bounds on excess risk and derive from them upper bounds on excess risk in abstract empirical risk minimization problems and in a more specialized context of regression and classification. In particular, this allows us to prove in a very easy way recent results of Tsybakov [44] on fast convergence rates in classification and, also for classification, to refine recent bounds of Massart and Nedelec [34], using a version of Alexander's capacity function.

The article is organized as follows. Section 2 contains the general exponential bounds for ratio empirical processes. Section 3 is devoted to moment bounds for empirical processes whose metric entropy with respect to the empirical $L_2$ distance is bounded by a regularly varying function independently of $P_n$; this includes, among others, VC-subgraph and VC-major classes. The reader interested in applications of the foregoing to ratios of margin distributions and to empirical risk minimization, may go directly from Section 3 to Sections 6 and 7, where we deal with these subjects. Sections 4 and 5 are devoted to rates (a.s. and in pr.), local and global moduli and limit theorems (including the central limit theorem) for ratio empirical processes.

**2. Concentration inequalities for normalized empirical processes.** In this section we derive the basic inequalities for ratio empirical processes. They are based on Talagrand's fundamental 1996 inequality, which will be formulated below. In what follows, $(S, \mathcal{S})$ is a measurable space, $P$ is a probability measure on it, $X_i$ are the coordinates $S^{\mathbf{N}} \mapsto S$, $\varepsilon_i$ are independent Rademacher variables independent of the variables $X_i$ (defined on, e.g., $([0,1], \lambda)$ and taking as $\Omega$ the product probability space $([0,1] \times S^{\mathbf{N}}, \lambda \times P^{\mathbf{N}} := \mathrm{Pr})$), $\mathcal{F}$ is a countable or suitably measurable (see, e.g., Dudley [17], Chapter 5) class of measurable functions on $S$ and $F$ is a measurable envelope of $\mathcal{F}$, that is, for all $f \in \mathcal{F}, x \in S$, $|f(x)| \leq F(x)$. For each $n$, $P_n$ is the empirical measure $n^{-1} \sum_{i=1}^{n} \delta_{X_i}$. As usual, we will also write $\|\psi(f)\|_{\mathcal{F}}$ for $\sup_{f \in \mathcal{F}} |\psi(f)|$.

TALAGRAND'S INEQUALITY. *For any measurable, uniformly bounded class of functions $\mathcal{F}$,*

$$(2.1) \quad \begin{aligned} \Pr\bigg\{ \bigg| \bigg\| \sum_{i=1}^{n} f(X_i) \bigg\|_{\mathcal{F}} - E \bigg\| \sum_{i=1}^{n} f(X_i) \bigg\|_{\mathcal{F}} \bigg| \geq t \bigg\} \\ \leq K \exp\bigg\{ -\frac{1}{K} \frac{t}{U} \log\bigg(1 + \frac{tU}{V}\bigg) \bigg\}, \end{aligned}$$

*valid for all $t > 0$, and where $K$ is a universal constant, $U$ is a uniform bound for the functions in $\mathcal{F}$ and $V$ is any number satisfying $V \geq E \sup_{f \in \mathcal{F}} \sum_{i=1}^{n} f^2(X_i)$.*

The inequality holds also for $\{X_i\}$ that are not necessarily identically distributed. The quantity $V$ is of course bounded by $n\|F\|_2^2$ if $F$ is a measurable envelope for the class $\mathcal{F}$, a trivial bound that, however, can sometimes be used. A more interesting bound that follows from randomization together with a contraction principle for Rademacher processes is the following, given by Talagrand ([42], Corollary 3.4):

$$(2.2) \quad E\bigg\| \sum_{i=1}^{n} f^2(X_i) \bigg\|_{\mathcal{F}} \leq n\sigma^2 + 8UE\bigg\| \sum_{i=1}^{n} \varepsilon_i f(X_i) \bigg\|_{\mathcal{F}},$$

where $\sigma^2 = \sup_{f \in \mathcal{F}} Ef^2(X_1)$ (see also [29]). Typically, Talagrand's inequality is used in combination with this bound for $V$.

In the sequel, throughout, we may drop the subindex $P$ in such notation as $\sigma_P f$ if no confusion arises, particularly in proofs. Given $0 < r \leq 1$ and $r < \delta \leq 1$, we set

$$\mathcal{F}(r) := \{ f \in \mathcal{F} : \sigma_P(f) \leq r \} \quad \text{and} \quad \mathcal{F}(r, s] := \mathcal{F}(s) \setminus \mathcal{F}(r);$$

for $1 < q \leq 2$ and $r < s \leq rq^l$ for some $l \in \mathbf{N}$, we let

$$\rho_j := rq^j, \qquad j = 0, \ldots, l,$$

and

$$\psi_{n,q}(u) := E\|P_n - P\|_{\mathcal{F}(\rho_{j-1},\rho_j]}, \qquad u \in (\rho_{j-1}, \rho_j], \ j = 1, \ldots, l.$$

Of course, given $\delta$ and $r$ we take $l$ to be the smallest integer larger than or equal to $\log_q(\delta/r)$. Given a continuous, strictly increasing function $\phi$ such that $\phi(0) = 0$, we also define

$$\phi_q(u) = \phi(\rho_j), \qquad u \in (\rho_{j-1}, \rho_j], j = 1, \ldots, l,$$

and

$$\beta_{n,q,\phi}(r,s] := \sup_{u \in (r,s]} \frac{\psi_{n,q}(u)}{\phi_q(u)},$$

and sometimes we may use instead the nondiscretized version, namely

$$\tilde{\beta}_{n,q,\phi}(r,s] := \sup_{u \in (r,s]} \frac{E\|P_n - P\|_{\mathcal{F}(uq^{-1},u]}}{\phi(u)}.$$

Some of the subindices of $\beta$ may be dropped in proofs. We also set

$$V_{n,q}(\rho_j) = V_n(\rho_j) := \frac{1}{n} E \left\| \sum_{i=1}^n (f(X_i) - Pf)^2 \right\|_{\mathcal{F}(\rho_{j-1},\rho_j]},$$

and note that, by (2.2) and the comment before it, if $F_j$ is a measurable envelope for $\mathcal{F}(\rho_{j-1}, \rho_j]$, then

(2.3) $$q^{-2}\rho_j^2 \leq V_n(\rho_j) \leq (PF_j^2) \wedge [\rho_j^2 + 16\psi_n(\rho_j)],$$

where for the lower bound we assume that $\sigma_P f = \mathrm{Var}_P^{1/2} f$. Finally, we let $\gamma$ be the inverse function of $\gamma^{-1}(x) := x \log(1+x)$, $0 \leq x \leq 1$. Note that

$$\gamma(x) \leq \begin{cases} \dfrac{2x}{\log(1+x)}, & \text{for } x \geq 0, \\ \dfrac{2x}{\log x}, & \text{for } x \geq 2, \\ 2\sqrt{x}, & \text{for } 0 \leq x \leq 2. \end{cases}$$

Denote

$$E_{n,q,\phi} := E \sup_{\substack{f \in \mathcal{F} \\ r < \sigma_P f \leq \delta}} \frac{|P_n f - Pf|}{\phi_q(\sigma_P f)}$$

and

$$\tau_{n,q,\phi} := \tau_{n,q,\phi}(s_1, \ldots, s_l)$$

$$:= \max_{j\,:\,s_j > 2n\overline{V}_{n,q}(\rho_j)} \frac{2s_j}{n\phi(\rho_j)\log(s_j/(n\overline{V}_{n,q}(\rho_j)))}$$

$$\vee \max_{j\,:\,s_j \leq 2n\overline{V}_{n,q}(\rho_j)} 2\sqrt{\frac{s_j \overline{V}_{n,q}(\rho_j)}{n\phi^2(\rho_j)}}$$

for any $\overline{V}_{n,q}(\rho_j) \geq V_{n,q}(\rho_j)$.

The following immediate application of Talagrand's inequality holds the key to the ratio limit theorems to be obtained below. It shows that the supremum

$$\sup_{\substack{f \in \mathcal{F} \\ r < \sigma_P f \leq \delta}} \frac{|P_n f - Pf|}{\phi_q(\sigma_P f)}$$

concentrates with high probability around both $\beta_{n,q,\phi}$ and $E_{n,q,\phi}$ with the same magnitude of the deviations (of the order $\tau_{n,q,\phi}$). In particular, it also means that $\beta_{n,q,\phi}$ and $E_{n,q,\phi}$ are within $\sim \tau_{n,q,\phi}$ of each other.

THEOREM 2.1. *With the above definitions, there exist universal constants $K, C \in (0, \infty)$ such that for any sequence $s_j$ of positive numbers*

(2.4a)
$$\Pr\left\{ \left| \sup_{\substack{f \in \mathcal{F} \\ r < \sigma_P f \leq \delta}} \frac{|P_n f - Pf|}{\phi_q(\sigma_P f)} - \beta_{n,q,\phi} \right| \geq \tau_{n,q,\phi}(s_1, \ldots, s_l) \right\}$$
$$\leq K \sum_{j=1}^{l} e^{-s_j/K}$$

*and*

(2.4b)
$$\Pr\left\{ \left| \sup_{\substack{f \in \mathcal{F} \\ r < \sigma_P f \leq \delta}} \frac{|P_n f - Pf|}{\phi_q(\sigma_P f)} - E_{n,q,\phi} \right| \geq C\tau_{n,q,\phi}(s_1, \ldots, s_l) \right\}$$
$$\leq K \sum_{j=1}^{l} e^{-s_j/K}.$$

PROOF. Set $\mathcal{F}_j := \mathcal{F}(\rho_{j-1}, \rho_j]$. Then, we can rewrite Talagrand's inequality as

(2.5)
$$\Pr\left\{ |\|P_n - P\|_{\mathcal{F}_j} - E\|P_n - P\|_{\mathcal{F}_j}| \geq \overline{V}_n(\rho_j) \gamma\left(\frac{s_j}{n\overline{V}_n(\rho_j)}\right) \right\}$$
$$\leq K e^{-s_j/K},$$

$j = 1, \ldots, l$. Hence, with probability at least $1 - \sum_{j=1}^{l} e^{-s_j/K}$,

$$\left| \sup_{\substack{f \in \mathcal{F} \\ r < \sigma_P f \leq \delta}} \frac{|P_n f - Pf|}{\phi_q(\sigma_P f)} - \beta_{n,q,\phi} \right| \leq \max_{1 \leq j \leq l} \left| \frac{\|P_n - P\|_{\mathcal{F}_j}}{\phi(\rho_j)} - \frac{\psi_{n,q}(\rho_j)}{\phi(\rho_j)} \right|$$
$$\leq \max_{1 \leq j \leq l} \frac{\overline{V}_n(\rho_j) \gamma(s_j/(n\overline{V}_n(\rho_j)))}{\phi(\rho_j)}.$$

Now, (2.4a) follows from the bounds for the function $\gamma$, namely, if $s_j > 2n\overline{V}_n(\rho_j)$ we use $\gamma(x) \leq \frac{2x}{\log x}$ to get

$$\frac{\overline{V}_n(\rho_j) \, \gamma(s_j/(n\overline{V}_n(\rho_j)))}{\phi(\rho_j)} \leq \frac{2s_j}{n\phi(\rho_j)\log(s_j/(n\overline{V}_n(\rho_j)))}$$

and otherwise we use $\gamma(x) \leq 2\sqrt{x}$ to get

$$\frac{\overline{V}_n(\rho_j) \, \gamma(s_j/(n\overline{V}_n(\rho_j)))}{\phi(\rho_j)} \leq 2\sqrt{\frac{s_j \overline{V}_n(\rho_j)}{n\phi^2(\rho_j)}}.$$

Next, we will show that

$$(2.6) \quad \Pr\left\{\left|\sup_{\substack{f\in\mathcal{F}\\ r<\sigma_Pf\leq\delta}}\frac{|P_nf - Pf|}{\phi_q(\sigma_Pf)} - \beta_{n,q,\phi}\right| \geq 2s\tau_{n,q,\phi}\right\} \leq C(\{s_j\})\exp\{-s/K\},$$

where

$$C(\{s_j\}) := K\sum_{j=1}^{l} e^{-s_j/K},$$

which is supposed to be smaller than 1 (otherwise, the inequalities of the theorem are trivial). Integration of (2.6) immediately implies that, for some $C > 0$,

$$E\left|\sup_{\substack{f\in\mathcal{F}\\ r<\sigma_Pf\leq\delta}}\frac{|P_nf - Pf|}{\phi_q(\sigma_Pf)} - \beta_{n,q,\phi}\right| \leq C\tau_{n,q,\phi}$$

and, as a consequence,

$$\beta_{n,q,\phi} \leq E_{n,q,\phi} \leq \beta_{n,q,\phi} + C\tau_{n,q,\phi}.$$

The last bound shows that in (2.4a) $\beta_{n,q,\phi}$ can be replaced by $E_{n,q,\phi}$ if we multiply $\tau_{n,q,\phi}$ by a constant, which proves (2.4b).

To prove (2.6), we again use (2.5) with $s_j$ replaced by $s_j + s$. It is enough to assume that $s, s_j \geq 2$. The right-hand side of (2.5) becomes $K\exp\{-(s + s_j)/K\}$. If $s_j > 2n\overline{V}_n(\rho_j)$ (and $s_j + s$ is even larger), we argue as in the proof of (2.4a) to get

$$\overline{V}_n(\rho_j)\gamma\left(\frac{s_j+s}{n\overline{V}_n(\rho_j)}\right) \leq \frac{2(s_j+s)}{n\phi(\rho_j)\log((s_j+s)/(n\overline{V}_n(\rho_j)))}$$

$$\leq s\frac{2s_j}{n\phi(\rho_j)\log(s_j/(n\overline{V}_n(\rho_j)))} \leq s\tau_{n,q,\phi}$$

(using $s, s_j \geq 2$). On the other hand, if $s_j \leq 2n\overline{V}_n(\rho_j)$, we use subadditivity of $\gamma$,

$$\gamma(x+y) \leq \gamma(x) + \gamma(y),$$

which follows from the inequality

$$\gamma^{-1}(x+y) \geq \gamma^{-1}(x) + \gamma^{-1}(y)$$

that is easy to prove directly. This gives

$$\overline{V}_n(\rho_j)\gamma\left(\frac{s_j+s}{n\overline{V}_n(\rho_j)}\right) \leq \overline{V}_n(\rho_j)\gamma\left(\frac{s_j}{n\overline{V}_n(\rho_j)}\right) + \overline{V}_n(\rho_j)\gamma\left(\frac{s}{n\overline{V}_n(\rho_j)}\right).$$

The first term is bounded (as in the first part of the proof of Theorem 2.1) by

$$2\sqrt{\frac{s_j\overline{V}_n(\rho_j)}{n\phi^2(\rho_j)}} \leq s\sqrt{\frac{s_j\overline{V}_n(\rho_j)}{n\phi^2(\rho_j)}} \leq \frac{1}{2}\tau_{n,q,\phi}s.$$

The second term is dominated either by

$$2\sqrt{\frac{s\overline{V}_n(\rho_j)}{n\phi^2(\rho_j)}} \leq s\sqrt{\frac{s_j\overline{V}_n(\rho_j)}{n\phi^2(\rho_j)}} \leq \frac{1}{2}\tau_{n,q,\phi}s$$

in the case when $s \leq 2n\overline{V}_n(\rho_j)$, or otherwise [if $s > 2n\overline{V}_n(\rho_j)$] by

$$\frac{2s}{n\phi(\rho_j)\log(s/(n\overline{V}_n(\rho_j)))} \leq \frac{2}{\log 2}s\frac{1}{n\phi(\rho_j)}$$

which can be further bounded by

$$\frac{2}{\log 2}s\sqrt{\frac{s_j\overline{V}_n(\rho_j)}{n\phi^2(\rho_j)}} \leq \frac{3}{2}\tau_{n,q,\phi}s$$

[since we have

$$n\overline{V}_n(\rho_j) \geq s_j/2 \geq 1 \geq 1/s_j].$$

The result now follows easily. $\square$

We may want to normalize the empirical process $P_nf - Pf$ by $\phi(\sigma_P f)$ instead of $\phi_q(\sigma_P f)$; in this case we do not obtain a concentration inequality, but two very similar deviation inequalities (one from above and one from below), particularly if $\phi$ is regular enough. The above theorem gives the following:

COROLLARY 2.2. *Assume that the continuous nondecreasing function $\phi$ satisfies that the quantity $c_{q,r,\phi} = \sup_{r \leq x \leq 1} \phi(qx)/\phi(x)$ is finite for some $1 < q \leq 2$. Then, with $s_j$ and $K$ as in Theorem 2.1 and under the same assumptions, we have both*

$$\Pr\left\{c_{q,r,\phi}^{-1} \times \sup_{\substack{f \in \mathcal{F} \\ r < \sigma_P f \leq \delta}} \frac{|P_nf - Pf|}{\phi(\sigma_P f)} \geq \beta_{n,q,\phi} + \tau_{n,q,\phi}\right\} \leq K\sum_j e^{-s_j/K}$$

and

$$\Pr\left\{\sup_{\substack{f\in\mathcal{F}\\r<\sigma_P f\leq\delta}}\frac{|P_n f-Pf|}{\phi(\sigma_P f)}\leq\beta_{n,q,\phi}-\tau_{n,q,\phi}\right\}\leq K\sum_j e^{-s_j/K}.$$

A way to use these propositions is as follows: if we let

$$b_n=\beta_{n,q,\phi}\vee\tau_{q,n,\phi},$$

where $q$, $r$, $\delta$ and $\{s_j\}$ may depend on $n$, and take $s_j=s_{j,n}$ such that $\sum_{j=1}^{l_n}e^{-s_{j,n}/K}$ tends to zero, then the sequence

$$\frac{1}{b_n}\sup_{\substack{f\in\mathcal{F}\\r<\sigma_P f\leq\delta}}\frac{|P_n f-Pf|}{\phi(\sigma_P f)},\qquad n\in\mathbf{N},$$

is stochastically bounded. The following lemma of Alexander [2] allows one to get a.s. results.

LEMMA 2.3. *With the same notation as above, let $c_n/n\downarrow$, $r_n\downarrow$, $\sqrt{n}\delta_n\uparrow$ and $u_n\downarrow$. Set*

$$A_n=\{\sqrt{n}|P_n f-Pf|\geq c_n\phi(\sigma_P f)+u_n\text{ for some }f\in\mathcal{F}\text{ with }r_n\leq\sigma_P f\leq\delta_n\}$$

*and*

$$A_n^\varepsilon=\{\sqrt{n}|P_n f-Pf|\geq(1-\varepsilon)(c_n\phi(\sigma_P f)+u_n)\text{ for some }f\in\mathcal{F}$$

$$\text{with }r_n\leq\sigma_P f\leq\sqrt{(1+\varepsilon)}\delta_n\},$$

*and assume*

$$\inf\{c_n\phi(t)/t\colon n\geq 1,\ t\in[r_n,\delta_n]\}>0.$$

*Then, if for some $\varepsilon,\theta>0$*

$$\Pr(A_n^\varepsilon)=O(1/(\log n)^{1+\theta}),$$

*we have*

$$\Pr(A_n\ i.o.)=0.$$

Sometimes, $\psi_{n,q}$ [and therefore also $\beta_{n,q,\phi}$ and $V_{n,q}(\rho_j)$] is still too large because the envelope of $\mathcal{F}_j$ is too large. Then, one may further subdivide $\mathcal{F}_j$ into $N_j$ classes $\mathcal{F}_{j,k}$ with smaller envelopes and such that $N_j$ is not too large (perhaps of the order of $\log\rho_j^{-1}$). For instance, this happens with the $d$-dimensional distribution function as we see below. One may take $\mathcal{F}_{j,k}$ to be the intersection with $\mathcal{F}_j$ of each of the components of an optimal covering of $\mathcal{F}_j$ by $L_2(P)$ balls of radius $\tau\rho_j$, $0<\tau<1$, but other subdivisions are possible; in particular, $N_j$ could be 1 for some or all $j$. We can apply the same principle as in the proof of Theorem 2.1 and get a bound that takes this into account, as follows. Let $\mathcal{F}_j = \bigcup_{k=1}^{N_j} \mathcal{F}_{j,k}$, let

$$\psi_{n,q,j,k} := E\|P_n - P\|_{\mathcal{F}_{j,k}},$$

$$\overline{\beta}_{n,q,\phi} := \max_{j,k} \frac{\psi_{n,q,j,k}}{\phi(\rho_j)},$$

$$V_{n,q,j,k} := \frac{1}{n} \left\| \sum_{i=1}^{n} (f(X_i) - Pf)^2 \right\|_{\mathcal{F}_{j,k}},$$

$$\overline{\tau}_{n,q,\phi} := \max_{j,k:s_{j,k}>2n\overline{V}_{n,q,j,k}} \frac{2s_{j,k}}{n\phi(\rho_j)\log(s_{j,k}/(n\overline{V}_{n,q,j,k}))}$$

$$\vee \max_{j,k:\, s_{j,k} \leq 2n\overline{V}_{n,q,j,k}} 2\sqrt{\frac{s_{j,k}\overline{V}_{n,q,j,k}}{n\phi^2(\rho_j)}}$$

and $\overline{V}_{n,q,j,k} \geq V_{n,q,j,k}$. Then, we have the following.

THEOREM 2.1′. *With the above definitions and letting $s_{j,k}$ be a double sequence of positive numbers, there exists a universal constant $K \in (0, \infty)$ such that*

$$(2.4') \quad \Pr\left\{ \left| \sup_{\substack{f \in \mathcal{F} \\ r < \sigma_P f \leq \delta}} \frac{|P_n f - Pf|}{\phi_q(\sigma_P f)} - \overline{\beta}_{n,q,\phi} \right| \geq \overline{\tau}_{n,q,\phi} \right\} \leq K \sum_{j=1}^{l} \sum_{k=1}^{N_j} e^{-s_{j,k}/K}.$$

The analogues of inequality (2.4b) and of the one-sided inequalities of Corollary 2.2 hold as well.

REMARK 2.4 (On the choice of $s_j$). In general, we must take

$$s_j = K \log \frac{1}{c_j}$$

with $\sum_{j=1}^{l} c_j$ small, as in this case, $\sum_{j=1}^{l} e^{-s_j/K} = \sum_{j=1}^{l} c_j$. If we take $s_j = s$ for a number of $j$'s more or less comparable to $l$, then a good choice is to take

$$s = K' \log l$$

for some $K' > K$, so that

$$\sum_{j=1}^{l} e^{-s_j/K} \leq 1/l^{K'/K-1},$$

which will tend to zero if $l \to \infty$ [so, if $\log_{q_n}(\delta_n/r_n) \to \infty$]. Another possible choice is

$$s_j = sq^{\alpha j}$$

for some $\alpha > 0$, which gives

$$\sum_{j=1}^{\infty} e^{-sq^{\alpha j}/K} = \frac{q^{\alpha}}{q^{\alpha}-1} \sum_{j=1}^{\infty} q^{-\alpha j} e^{-sq^{\alpha j}/K}(q^{\alpha j} - q^{\alpha(j-1)})$$

(2.7)
$$< \frac{1}{q^{\alpha}-1} \int_1^{\infty} e^{-sx/K}\, dx$$

$$= K\frac{1}{q^{\alpha}-1}\frac{1}{s}e^{-s/K},$$

and

$$\sum_{j=0}^{\infty} e^{-sq^{\alpha j}/K} \le \left(1 + K\frac{1}{q^{\alpha}-1}\frac{1}{s}\right)e^{-s/K}$$

(2.7′)
$$< K\frac{q^{\alpha}}{q^{\alpha}-1}\frac{1}{s}e^{-s/Kq^{\alpha}},$$

bounds that can be made small by increasing $s$. Finally, another choice is

$$s_j = s_n + K\log\log_q(q\delta/\rho_j),$$

as is easy to check. However, there does not seem to be a choice of $s_j$ that works in all situations (see, e.g., the last part of Example 2.7 below for an unexpected choice for $s_j$).

REMARK 2.5 (The role of the stratification $\mathcal{F}_{(r,\delta]} \subseteq \bigcup_{i=1}^{l} \mathcal{F}_{(\rho_{i-1},\rho_i]}$). Since on each stratum $\mathcal{F}_j := \mathcal{F}_{(\rho_{i-1},\rho_i]}$ the function of $f \mapsto \phi(\sigma_P f)$ is essentially constant [assuming that $\phi(u) \simeq \phi(qu)$], we have

$$\left\|\frac{P_n f - Pf}{\phi(\sigma_P f)}\right\|_{\mathcal{F}_j} \simeq \frac{1}{\phi(\rho_j)}\|P_n f - Pf\|_{\mathcal{F}_j},$$

which is why the terms in the bounds (2.4) do depend on the complexity of these strata [measured by $\psi_{n,q}(\rho_j)$ and $V_n(\rho_j)$, which ultimately also depends on $\psi_{n,q}(\rho_j)$], usually simpler than the complexity of $\mathcal{F}$. Instead of stratifying, we could simply apply Talagrand's inequality to the class of functions $\{\phi(r)f/\phi(\sigma_P f) : f \in \mathcal{F}, \sigma_P f \geq r\}$. But these classes are more complicated than $\mathcal{F}$ and so would be the parameters of the inequality. These parameters often depend on the $L_2$ norm of the envelope of the corresponding classes, and there may be a good advantage in using the *local envelopes* $\sup\{|f(x)| : f \in \mathcal{F}_j\}/\phi(\rho_j)$ rather than the global $\sup\{|f(x)/\phi(\sigma_P f)| : f \in$

NORMALIZED EMPIRICAL PROCESSES 13

$\mathcal{F}\}$. This advantage comes at a cost, at least with regard to distributional or in probability results: whereas the sequence $s_j$ should be large enough so that the series $\sum e^{-s_j/K}$ converges and has a small sum, we do not have to deal with this series if we apply Talagrand's inequality to the whole class. In this last case $s$ can be any number such that $e^{-s}$ is of the desired size. However, if one wants to apply Alexander's lemma, then $s$ must be of the order of $\log \log n$, which in general is comparable to $\log l$, hence to $s_j$ if $s_j$ does not depend on $j$. This cost is usually overwhelmed by the mentioned advantage, and in the worst case, the number $\tau_{n,q,\phi}$ in (2.4) is at most a factor of $\log l$, or even $\sqrt{\log l}$, larger than it would be by direct application of Talagrand, and not larger at all (except perhaps for a constant factor) if we want the probability bound to be of the order of $1/(\log l)^{1+\theta}$.

REMARK 2.6. Another approach, used, for instance, by Massart [31] or Bousquet [9], is based on stratification, but uses Talagrand's inequality only once, which is relevant to Remark 2.5, but which results in other losses when the class of functions is relatively small. We briefly describe this approach. Suppose that for all $\rho > 0$

$$E\|P_n - P\|_{\mathcal{F}(\rho)} \leq \tilde{\psi}_n(\rho),$$

where $\tilde{\psi}_n$ is a function satisfying that for some $\lambda \in (0,1)$ the function $\rho \mapsto \frac{\tilde{\psi}_n(\rho)}{\phi^\lambda(\rho)}$ is nonincreasing. Assume also that with the same $\lambda$

$$\sum_{j:\rho_j \geq r} \frac{1}{\phi(\rho_j)^{1-\lambda}} \leq c_{q,\lambda,\phi} \frac{1}{\phi(r)^{1-\lambda}}$$

for some constant $c_{q,\lambda,\phi}$. Note that these conditions immediately imply that

$$\sum_{j:\rho_j \geq r} \frac{\tilde{\psi}_n(\rho_j)}{\phi(\rho_j)} = \sum_{j:\rho_j \geq r} \frac{\tilde{\psi}_n(\rho_j)}{\phi(\rho_j)^\lambda \phi(\rho_j)^{1-\lambda}}$$

$$\leq \frac{\tilde{\psi}_n(r)}{\phi(r)^\lambda} \sum_{j:\rho_j \geq r} \frac{1}{\phi(\rho_j)^{1-\lambda}}$$

$$\leq c_{q,\lambda,\phi} \frac{\tilde{\psi}_n(r)}{\phi(r)}.$$

Consider now the class

$$\mathcal{G} := \bigcup_{j:\delta \geq \rho_j \geq r} \frac{\phi(r)}{\phi(\rho_j)} \mathcal{F}_j,$$

which is also bounded by 1. We have

$$E\|P_n - P\|_{\mathcal{G}} \leq \sum_{j:\rho_j \geq r} \frac{\phi(r)}{\phi(\rho_j)} E\|P_n - P\|_{\mathcal{F}_j}$$

$$\leq \sum_{j:\rho_j \geq r} \frac{\phi(r)}{\phi(\rho_j)} \tilde{\psi}_n(\rho_j)$$

$$\leq c_{q,\lambda,\phi} \tilde{\psi}_n(r).$$

Using (2.2), this gives

$$V_n(\mathcal{G}) := \frac{1}{n} E \left\| \sum_{i=1}^{n} (f(X_i) - Pf)^2 \right\|_{\mathcal{G}}$$

$$\leq \sup_{j:\rho_j \geq r} \frac{\phi(r_n)^2}{\phi(\rho_j)^2} \rho_j^2 + 16 E\|P_n - P\|_{\mathcal{G}}.$$

We will assume that either $\rho \mapsto \frac{\rho}{\phi(\rho)}$ is nonincreasing (case 1), or it is nondecreasing (case 2). In case 1,

$$V_n(\mathcal{G}) \leq r^2 + 16 c_{q,\lambda,\phi} \tilde{\psi}_n(r) =: \overline{V}_n(\mathcal{G}).$$

Applying Talagrand's inequality to the class $\mathcal{G}$ the same way we did in the proof of Theorem 2.1, we get that with probability at least $1 - K e^{-s/K}$

$$\left| \sup_{\substack{f \in \mathcal{F} \\ r < \sigma_P f \leq \delta}} \frac{|P_n f - Pf|}{\phi_q(\sigma_P f)} - E_{n,q,\phi} \right|$$

$$= \frac{1}{\phi(r)} |\|P_n - P\|_{\mathcal{G}} - E\|P_n - P\|_{\mathcal{G}}|$$

$$\leq I(s > 2n\overline{V}_n(\mathcal{G})) \frac{2s}{n\phi(r) \log(s/(n\overline{V}_n(\mathcal{G})))}$$

$$+ I(s \leq 2n\overline{V}_n(\mathcal{G})) 2 \sqrt{\frac{s\overline{V}_n(\mathcal{G})}{n\phi(r)^2}}$$

(2.8)

$$= I\left( \frac{s}{nr^2} > 2\left(1 + 16 c_{q,\lambda,\phi} \frac{\tilde{\psi}_n(r)}{r^2}\right)\right)$$

$$\times \frac{2s}{n\phi(r) \log(s/(nr^2(1 + 16 c_{q,\lambda,\phi} \tilde{\psi}_n(r)/r^2)))}$$

$$+ I\left( \frac{s}{nr^2} \leq 2\left(1 + 16 c_{q,\lambda,\phi} \frac{\tilde{\psi}_n(r)}{r^2}\right)\right)$$

$$\times 2 \sqrt{\frac{s}{n} \frac{r^2}{\phi(r)^2} \left(1 + 16 c_{q,\lambda,\phi} \frac{\tilde{\psi}_n(r)}{r^2}\right)}.$$

Similarly, in case 2,

$$V_n(\mathcal{G}) \le c_\phi \phi(r)^2 + 16 c_{q,\lambda,\phi} \tilde{\psi}_n(r) =: \overline{V}_n(\mathcal{G}).$$

Again, by Talagrand's inequality, with probability at least $1 - Ke^{-s/K}$

$$\left| \sup_{\substack{f \in \mathcal{F} \\ r < \sigma_P f \le \delta}} \frac{|P_n f - Pf|}{\phi_q(\sigma_P f)} - E_{n,q,\phi} \right|$$

$$\le I(s > 2n\overline{V}_n(\mathcal{G})) \frac{2s}{n\phi(r) \log(s/(n\overline{V}_n(\mathcal{G})))} + I(s \le 2n\overline{V}_n(\mathcal{G})) 2\sqrt{\frac{s\overline{V}_n(\mathcal{G})}{n\phi(r)^2}}$$

$$(2.9) = I\left(\frac{s}{n\phi(r)^2} > 2\left(c_\phi + 16 c_{q,\lambda,\phi} \frac{\tilde{\psi}_n(r)}{\phi(r)^2}\right)\right)$$

$$\times \frac{2s}{n\phi(r) \log(s/(n\phi(r)^2 (c_\phi + 16 c_{q,\lambda,\phi} \tilde{\psi}_n(r)/\phi(r)^2)))}$$

$$+ I\left(\frac{s}{n\phi(r)^2} \le 2\left(c_\phi + 16 c_{q,\lambda,\phi} \frac{\tilde{\psi}_n(r)}{\phi(r)^2}\right)\right) 2\sqrt{\frac{s}{n}\left(c_\phi + 16 c_{q,\lambda,\phi} \frac{\tilde{\psi}_n(r)}{\phi(r)^2}\right)}.$$

Bartlett and Mendelson [7] used another approach that also allowed them to apply Talagrand's inequality only once under an extra geometric assumption on the class $\mathcal{F}$ (namely, that the class is star-shaped).

EXAMPLE 2.7. As an illustration of Theorem 2.1, we will recover known results on the a.s. and the in probability behavior of

$$T_n := \sup_{1/n < t \le 1/2} \frac{|F_n(t) - t|}{\sqrt{t}},$$

where $F_n$ is the empirical distribution function corresponding to $n$ independent samples from the uniform distribution on $[0,1]$. In this case, $\mathcal{F} = \{I_{[0,t]} : 0 \le t \le 1/2\}$, $\sigma_P I_C = \sqrt{PC}$, $\phi(t) = t$, $r = r_n = 1/\sqrt{n}$, $\delta = 1/4$, $q$ is any number between 1 and 2, say 2, $l_n = (\log n/4)/2 \log 2$, $\mathcal{F}_j = \{I_{[0,t]} : t \le \rho_j^2\}$, we can take $\overline{V}_n(\rho_j) = 2\rho_j^2$ and

$$\psi_n(\rho_j) \le \frac{2}{n} E \sup_{t \le \rho_j^2} \left| \sum_{i=1}^n \varepsilon_i I_{[0,t]}(X_i) \right| \le \frac{4}{n} E \left| \sum_{i=1}^n \varepsilon_i I_{[0,\rho_j^2]}(X_i) \right| \le \frac{4\rho_j}{\sqrt{n}},$$

where the first inequality follows by symmetrization and the second by Lévy's inequality. So, the quantity $\tau_{n,q,\phi}$ of Theorem 2.1 is

$$(2.10) \qquad \max_{j \,:\, s_j > 4n\rho_j^2} \frac{2s_j}{n\rho_j \log(s_j/(2n\rho_j^2))} \vee 2 \max_{j \,:\, s_j \le 4n\rho_j^2} \sqrt{\frac{s_j}{n}}.$$

Then, if we take $s_j = K'q\log l_n \geq K''\log\log n$ with $K'' > 4K$, this bound is dominated by

$$\frac{2K''\log\log n}{\sqrt{n}\log\log\log n}$$

and

$$\sum_{j=1}^{l_n} e^{-s_j/2K} \lesssim (\log n)^{-(K''/2K-1)},$$

which, by Lemma 2.3, give

$$\limsup_{n\to\infty} \frac{\sqrt{n}\log\log\log n}{\log\log n} T_n = C < \infty \qquad \text{a.s.}$$

(as $\beta_n \leq 4/\sqrt{n}$ multiplied by the factor in front of $T_n$ tends to zero). This is sharp: Csáki [13] computed this constant, which is not zero. Suppose now we want to find the order of magnitude of $T_n$ in probability. Then we still take $s_j = K''\log\log n$ if $K''\log\log n \leq 4n\rho_j^2$ and notice that the number of the remaining $j$'s, such that $4n\rho_j^2 \leq K''\log\log n$ is of the order of $K'''\log\log\log n$; this allows us to take a smaller $s_j$ for such $j$'s, for instance, we can take $s_j = 2K\log\log\log n$ and still have

$$\sum_{j=1}^{l_n} e^{-s_j/2K} \lesssim \frac{K'''\log\log\log n}{\log\log n} + l_n^{-(K''/2K-1)} \to 0.$$

The bound (2.10) becomes of the order

$$\frac{\log\log\log n}{\sqrt{n}\log\log\log\log n} \vee \sqrt{\frac{\log\log\log n}{n}} \vee \sqrt{\frac{\log\log n}{n}} = \sqrt{\frac{\log\log n}{n}}.$$

This gives that the sequence

$$\sqrt{\frac{n}{\log\log n}} T_n, \qquad n \in \mathbf{N},$$

is stochastically bounded, a result that is also best possible since it follows from [14] that it converges in probability to a positive constant. We should remark that these results on the almost sure and the in probability size of $T_n$ can also be obtained by direct application of Talagrand's inequality to the class of functions $\{I_{[0,t]}/\sqrt{t}: t \leq 1/4\}$, although the estimation of expected values in this case is more complicated. However, they do not follow from the method developed in Remark 2.6. At any rate, this example illustrates the power of Theorems 2.1 and 2.1' when good expectation bounds are available, and also how to choose $s_j$. See Examples 4.9, 4.10 and 5.8 below for more on uniform empirical c.d.f.'s, in one or more dimensions.

Typically one wishes to normalize the empirical process by $\sqrt{\operatorname{Var}_P f}$, which corresponds to $\phi(x) = x$, or by $Pf$, which corresponds to $\phi(x) = x^2$ and $\sigma_P f = \sqrt{Pf}$ (recall $0 \leq f \leq 1$), or by a function of $\sigma_P f$ of the form $\phi(x) = xL(x)$ with $L$ slowly varying at zero. Although other situations may be considered, we will only specialize Theorem 2.1 to a small number of cases, including these. The main job consists in choosing $s_j$ so that $\sum e^{-s_j/K}$ is small. The following proposition recovers an inequality in [22].

PROPOSITION 2.8. *With the notation of Theorem 2.1 for $\sigma_P f = \sqrt{Pf}$ and $\phi(t) = t^2$, we have that, for all $s > 0$,*

$$\Pr\Biggl\{ q^{-1} \sup_{\substack{f \in \mathcal{F} \\ r^2 < Pf \leq \delta}} \left| \frac{P_n f}{Pf} - 1 \right| \geq \beta_{n,q,\phi} + 2\sqrt{\frac{s}{nr^2}(1 + 16\beta_{n,q,\phi})}$$

(2.11)
$$\vee \frac{2s}{nr^2 \log((s/(nr^2(1 + 16\beta_{n,q,\phi}))) \vee 2)} \Biggr\}$$

$$\leq K^2 \frac{1}{q^2 - 1} \frac{1}{s} e^{-s/K}$$

*and*

$$\Pr\Biggl\{ \sup_{\substack{f \in \mathcal{F} \\ r^2 < Pf \leq \delta}} \left| \frac{P_n f}{Pf} - 1 \right| \leq \beta_{n,q,\phi} - 2\sqrt{\frac{s}{nr^2}(1 + 16\beta_{n,q,\phi})}$$

(2.12)
$$\vee \frac{2s}{nr^2 \log((s/(nr^2(1 + 16\beta_{n,q,\phi}))) \vee 2)} \Biggr\}$$

$$\leq K^2 \frac{1}{q^2 - 1} \frac{1}{s} e^{-s/K}.$$

PROOF. Take $\overline{V}_n(\rho_j) = \rho_j^2(1 + 16\beta_n)$ [which is allowed by (2.3)] and $s_j = sq^{2j}$ in Corollary 2.2, and use (2.7) with $\alpha = 2$ to compute the probability bound. □

Especially important is the case $\phi(t) = t$, that is, the normalization of the empirical process at $f$ by the standard deviation of $f(X)$. The following proposition is slightly sharper (up to constants) than a similar inequality in [22], and applies in a larger range.

PROPOSITION 2.9. *Let $\phi(t) = t$. Set $c_q = \max_{1 \leq j \leq l_n} (\log j)/q^j$ and denote*

$$B_n := \frac{10\sqrt{s/q + 2c_q K}\sqrt{s + 2K \log \log_q(q\delta/r)}}{nr \log(((5s/q + 10c_q K)/(17nr^2)) \vee 10)}$$

$$\vee\, 2\sqrt{17}\sqrt{\frac{s+2K\log\log_q(q\delta/r)}{n}}.$$

(a) *If* $\beta_n = \beta_{n,q,\phi} \leq r$, *then, for any positive number* $s$,

(2.13) $$\Pr\left\{\left|\sup_{\substack{f\in\mathcal{F}\\r<\sigma_P f\leq\delta}}\frac{|P_nf-Pf|}{\phi_q(\sigma_P f)}-\beta_n\right|\geq B_n\right\}\leq 2Ke^{-s}$$

*with obvious changes in the constants if* $\beta_n \leq Cr$ *for some* $C < \infty$.

(b) *If* $\beta_n > r$, *then, for any* $s > 0$, $t > 0$,

$$\Pr\left\{\left|\sup_{\substack{f\in\mathcal{F}\\r<\sigma_P f\leq\delta}}\frac{|P_nf-Pf|}{\phi_q(\sigma_P f)}-\beta_n\right|\geq \frac{2t}{nr\log((t/(17nr\beta_n))\vee 2)}\right.$$

(2.14) $$\left.\vee\, 2\sqrt{17}\sqrt{\frac{t\beta_n}{nr}}\vee B_n\right\}$$

$$\leq K^2\frac{1}{q-1}\frac{1}{t}e^{-t/K}+2Ke^{-s},$$

*with obvious changes in the constants if* $r < C\beta_n$ *for some* $C < \infty$.

PROOF. Assume first $\beta_n \leq r$. Since $\psi_n(\rho_j)/\rho_j \leq \beta_n \leq r \leq \rho_j$ we can take $\overline{V}_n(\rho_j) = 17\rho_j^2$. Now take $s_j = s + 2K\log\log_q(\rho_j/r_n) = s + 2K\log j$ if this quantity does not exceed $34n\rho_j^2$ and five times this quantity otherwise. Then, to estimate

$$\tau_{n,q,t} = \max_{j:\, s_j > 2n\overline{V}_n(\rho_j)} \frac{2s_j}{n\rho_j\log(s_j/(n\overline{V}_n(\rho_j)))} \vee \max_{s_j\leq 2n\overline{V}_n(\rho_j)} 2\sqrt{\frac{s_j\overline{V}_n(\rho_j)}{n\rho_j^2}},$$

note that if $x \geq e^2$, then $x^{1/2}/\log x$ is nondecreasing, so that

$$\frac{10s+20K\log j}{nrq^j\log((5s+10K\log j)/(17nr^2q^{2j}))}$$

$$\leq \frac{\sqrt{10s/q+20c_qK}\sqrt{10s+20K\log\log_q(q\delta/r)}}{nr\log((5s/q+10c_qK)/(17nr^2))},$$

which gives

$$\tau_{n,q,t} \leq B_n.$$

Moreover, $K\sum_{j=1}^{l_n} e^{-s_j/K} \leq Ke^{-s}\sum_{j=1}^{\infty} 1/j^2 \leq 2Ke^{-s}$.

Assume now $r < \beta_n$. The $j$'s for which $\rho_j \geq \beta_n$ can be treated as in the previous case (where all the $\rho_j$ were larger than or equal to $\beta_n$). If $\rho_j < \beta_n$, then $V_n(\rho_j)/\rho_j \leq \rho_j + 16\beta_n \leq 17\beta_n$ and we take $\overline{V}_n(\rho_j) = 17\rho_j\beta_n$. Then, with $s_j = tq^j$, the contribution of these $j$ to $\tau_{n,q,t}$ is easily seen to be dominated by

$$\frac{2t}{nr\log((t/17nr\beta_n) \vee 2)} \vee 2\sqrt{17}\sqrt{\frac{t\beta_n}{nr}},$$

and, by (2.7), their contribution to the probability bound is dominated by $K^2\frac{q}{q-1}\frac{1}{t}e^{-t/qK}$. □

Comparing with inequality (2.8) in Remark 2.6, we see that the result in the previous proposition is better if in (2.8) we take $s$ of the order of $\log\log_q(q\delta/r)$ because $\tilde{\phi}_n(r) > \psi_n(r)$, but smaller $s$'s are possible in (2.8), and these do better than $s + 2K\log\log_q(q\delta/r)$.

By (2.3), if $\|F_j\|_2$ is comparable to $\rho_j$, then we can take $\overline{V}_n(\rho_j) = c\rho_j^2 \simeq V_n(\rho_j)$ and obtain better inequalities than (2.13) and (2.14); this is the case in Example 2.7 and, if one uses Theorem 2.1$'$, this is also the case for the c.d.f. in several dimensions (as Alexander [2] observed and we see below). This applies also to Proposition 2.8 and to the ones that follow.

REMARK 2.10. If in the previous proposition we assume

$$(2.15) \qquad r \vee \beta_n \geq \sqrt{\frac{s + 2K\log\log_q(q\delta/r)}{34n}},$$

then the Poisson term of $\tau_{n,q,t}$ can be deleted from the bounds; under this condition, if $\beta_n \leq r$, then we see that $s_j \leq 2n\overline{V}_n(\rho_j)$ for all $j$. Under condition (2.15), if $r < \beta_n$, the same is true for $\rho_j \geq \beta_n$. Thus we obtain

$$(2.16) \qquad \Pr\left\{\left|\sup_{\substack{f \in \mathcal{F} \\ r < \sigma_P f \leq \delta}} \frac{|P_n f - Pf|}{\phi_q(\sigma_P f)} - \beta_{n,q,\phi}\right| \geq 2\sqrt{17}\sqrt{\frac{s + 2K\log\log_q(q\delta/r)}{n}}\right\} \leq 2Ke^{-s},$$

if $\beta_n \leq r$, and

$$\Pr\left\{\left|\sup_{\substack{f \in \mathcal{F} \\ r < \sigma_P f \leq \delta}} \frac{|P_n f - Pf|}{\phi_q(\sigma_P f)} - \beta_{n,q,\phi}\right| \geq \frac{2t}{nr_n\log((t/(17nr\beta_n)) \vee 2)} \vee 2\sqrt{17}\sqrt{\frac{t\beta_n}{nr}}\right.$$

$$(2.17) \qquad\qquad\qquad\qquad\qquad \left.\vee 2\sqrt{17}\sqrt{\frac{s + 2K\log\log_q(q\delta/r)}{n}}\right\}$$

$$\leq K^2\frac{1}{q-1}\frac{1}{t}e^{-t/K} + 2Ke^{-s}$$

if $r < \beta_n$.

PROPOSITION 2.11. *Let $\phi(t) = t^\alpha$ for some $\alpha \in (1,2)$. Then, for any positive number $s$,*

$$\Pr\left\{\left|\sup_{\substack{f \in \mathcal{F} \\ r < \sigma_P f \leq \delta}} \frac{|P_n f - P f|}{\phi_q(\sigma_P f)} - \beta_n \right| \geq \frac{10s}{nr^\alpha \log((s/(17 n r^\alpha (r^{2-\alpha} \vee \beta_n))) \vee 10)}\right.$$

(2.18)
$$\left. \vee\, 2\sqrt{17}\sqrt{\frac{s(r^{2-\alpha} \vee \beta_n)}{nr^\alpha}}\,\right\}$$

$$\leq K^2 \frac{1}{q^\tau - 1} \frac{1}{s} e^{-s/K},$$

*where $\tau = 2(\alpha - 1)$.*

PROOF. One proceeds as in the previous proof by considering the cases $\rho_j^{2-\alpha} \geq \beta_n$ and $\rho_j^{2-\alpha} \leq \beta_n$, which correspond to $\overline{V}_n(\rho_j) = 17\rho_j^2$ and $\overline{V}_n(\rho_j) = 17\rho_j^\alpha \beta_n$, respectively; in the first case one takes $s_j = s q^{2(\alpha-1)j}$ or five times this, and in the second, $s_j = s q^{\alpha j}$.  □

The bounds for $\phi(t) = t^\alpha$, $0 < \alpha < 1$, are similar to those for $\alpha = 1$. We only state them in a case analogous to Remark 2.10.

PROPOSITION 2.12. *Let $\phi(t) = t^\alpha$, $\alpha \in (0,1)$, and assume*

(2.15′)
$$r \vee \beta_n^{1/(2-\alpha)} \geq \sqrt{\frac{s + 2K \log\log_q(q^2 \delta/r)}{n}}.$$

(a) *If $\beta_n \leq r^{2-\alpha}$, then for all $s > 0$,*

$$\Pr\left\{\left|\sup_{\substack{f \in \mathcal{F} \\ r < \sigma_P f \leq \delta}} \frac{|P_n f - P f|}{\phi_q(\sigma_P f)} - \beta_n\right|\right.$$

(2.19)
$$\left. \geq 2\sqrt{17}\sqrt{\frac{s\delta^{2(1-\alpha)} + 2K c_{q,\alpha}}{n}}\,\right\} \leq 2K e^{-s},$$

*where $c_{q,\alpha} = \sup_{0 < u \leq \delta q} u^{2(1-\alpha)} \log\log_q(q^2 \delta/u)$.*

(b) *If $\beta_n > r^{2-\alpha}$, $s, t > 0$ and*

$$B_n := 2\sqrt{17}\sqrt{\frac{s\delta^{2(1-\alpha)} + 2K c_{q,\alpha}}{n}} \vee \frac{2t}{nr^\alpha \log((t/(17 n r^\alpha \beta_n)) \vee 2)} \vee 2\sqrt{17}\sqrt{\frac{t\beta_n}{nr^\alpha}},$$

*then*

(2.20) $\Pr\left\{\left|\sup_{\substack{f \in \mathcal{F} \\ r < \sigma_P f \leq \delta}} \dfrac{|P_n f - P f|}{\phi_q(\sigma_P f)} - \beta_n\right| \geq B_n\right\} \leq K^2 \dfrac{1}{q^\alpha - 1} \dfrac{1}{t} e^{-t/K} + 2K e^{-s}.$

PROOF. Take $\overline{V}_n(\rho_j) = 17\rho_j^\alpha \beta_n$ if $\rho_j^{2-\alpha} \leq \beta_n$ and $\overline{V}_n(\rho_j) = 17\rho_j^2$ otherwise. In the first case, set $s_j = sq^{\alpha j}$ and in the second $s_j = s + 2K \log\log_q(q^2\delta/\rho_j)$ or $e^{2/\alpha}/2$ times this. $\square$

The case $\phi(t) = t^\alpha L(1/t)$ with $L$ monotone and slowly varying at infinity is also easy to handle, and we will when needed.

It should be noted that the bounds in the last three propositions are sharp only to the extent that the estimate $V_n(\rho_j) \leq \rho_j^2 + 16\psi_n(\rho_j)$ is sharp. Sometimes the class $\mathcal{F}_j$ can be further decomposed into a relatively small number of classes $\mathcal{F}_{j,k}$ for which $V_{n,j,k} \leq c\rho_j^2$, and then it is Theorem 2.1' that gives inequalities leading to sharp results.

**3. Inequalities for expected values of suprema of empirical processes under uniform, regularly varying (or slowly varying) entropy bounds.** We need good bounds for $\psi_n(\rho_j)$ in order to apply the inequalities in Section 2. In this section we do this for a large collection of classes of functions that includes the ubiquitous VC classes. In the theorems below, we assume that the functions in $\mathcal{F}$ take their values in $[-1, 1]$ and they are $P$-centered, and $F \leq 1$ denotes a measurable envelope of $\mathcal{F}$. For each $n$, we set

$$\|F\|_2 := \|F\|_{L_2(P)}, \qquad \|F\|_{2,n} := \|F\|_{L_2(P_n)}, \qquad n \in \mathbf{N},$$

and let $\sigma$ be a positive number such that

$$(3.1) \qquad \sup_{f \in \mathcal{F}} Pf^2 \leq \sigma^2 \leq \|F\|_2^2,$$

unless we specify

$$(3.2) \qquad \sigma^2 = \sup_{f \in \mathcal{F}} Pf^2.$$

We also let $H:[0,\infty) \mapsto [0,\infty)$ be *a regularly varying function of exponent* $0 \leq \alpha < 2$, *strictly increasing for $x \geq 1/2$ and such that $H(x) = 0$ for $0 \leq x < 1/2$*. Given such a function, we let the quantities $C_H$, $D_H$, $A_H$ satisfy

$$\infty > C_H \geq \sup_{x \geq 1} \frac{\int_x^\infty u^{-2}\sqrt{H(u)}\,du}{x^{-1}\sqrt{H(x)}} \vee 1,$$

$$\infty > D_H \geq \int_1^\infty u^{-2}\sqrt{H(u)}\,du,$$

$$\infty > A_H \geq \sup_{x \geq 2} \frac{\log(D_H x/(4C_H \sqrt{H(x)}))}{x^2} \vee 1.$$

Finally, if $(T,d)$ is a pseudometric space and $\varepsilon > 0$, then $N(T,d,\varepsilon)$ denotes the $\varepsilon$ covering number of $(T,d)$ (the smallest number of open balls of radius at most $\varepsilon$ needed to cover $T$) and $D(T,d,\varepsilon)$ denotes the $\varepsilon$ packing number (the largest possible number of elements in $T$ separated from each other by at least a distance of $\varepsilon$), and recall the elementary inequality

$$N(T,d,\varepsilon) \leq D(T,d,\varepsilon) \leq N(T,d,\varepsilon/2)$$

for all $\varepsilon > 0$, that we will use without further mention.

THEOREM 3.1. *If*

$$\log N(\mathcal{F}, L_2(P_n), \tau) \leq H\left(\frac{\|F\|_{2,n}}{\tau}\right) \tag{3.3}$$

*for all $\tau > 0$, $n \in \mathbf{N}$ and $\omega \in \Omega$, then there is a positive constant $C(H)$ that depends only on $A_H$, $C_H$ and $D_H$, such that*

$$\begin{aligned}
\frac{1}{C(H)} E\left\|\sum_{i=1}^n f(X_i)\right\|_{\mathcal{F}} & \\
\leq [\sqrt{n}\|F\|_2] &\wedge \left[\left(\sqrt{n}\sigma\sqrt{H\left(\frac{2\|F\|_2}{\sigma}\right)}\right)\right. \\
&\left. \vee H\left(\frac{2\|F\|_2}{\sigma} \wedge \frac{\sqrt{n}\|F\|_2}{1440 C_H}\right) \vee 1\right].
\end{aligned} \tag{3.4}$$

*The bound*

$$\frac{1}{C(H)} E\left\|\sum_{i=1}^n f(X_i)\right\|_{\mathcal{F}} \leq \left[\sqrt{n}\sigma\sqrt{H\left(\frac{2\|F\|_2}{\sigma}\right)}\right] \vee \left[H\left(\frac{2\|F\|_2}{\sigma}\right)\right], \tag{3.5}$$

*which also holds in general, is useful when $n\sigma^2 \geq c > 0$. Finally, for any $c > 0$, if*

$$n\sigma^2 \geq cH\left(\frac{2\|F\|_2}{\sigma}\right), \tag{3.6}$$

*then*

$$E\left\|\sum_{i=1}^n f(X_i)\right\|_{\mathcal{F}} \leq K(H,c)\sqrt{n}\sigma\sqrt{H\left(\frac{2\|F\|_2}{\sigma}\right)}, \tag{3.7}$$

*for a constant $K(H,c)$ that depends only on $H$ and $c$.*

PROOF. We delete the subscript $\mathcal{F}$ from norms when no confusion may arise. By standard symmetrization, $E\|\sum_{i=1}^n f(X_i)\| \leq 2E\|\sum_{i=1}^n \varepsilon_i f(X_i)\|$. Set

$$\sigma_n^2 := \|P_n f^2\|_{\mathcal{F}}.$$

The usual entropy bound for sub-Gaussian processes (for the constant, we combine the last display on page 320 of [29] with Theorem 11.17 and first display on page 322 of [29]) gives

$$E\left\|\frac{1}{\sqrt{n}}\sum_{i=1}^n \varepsilon_i f(X_i)\right\|$$
$$\leq 60 E \int_0^{2\sigma_n} \sqrt{\log N(\mathcal{F}, L_2(P_n), \tau)}\, d\tau$$

(3.8)
$$\leq 120 E \int_0^{\sigma_n} \sqrt{H\left(\frac{\|F\|_{2,n}}{\tau}\right)}\, d\tau$$

$$\leq 120 E\left[\int_0^{\sigma_n} \sqrt{H\left(\frac{2\|F\|_2}{\tau}\right)}\, d\tau\, I(\|F\|_{2,n} \leq 2\|F\|_2)\right]$$

$$+ 120 E\left[\int_0^{\sigma_n} \sqrt{H\left(\frac{\|F\|_{2,n}}{\tau}\right)}\, d\tau\, I(\|F\|_{2,n} > 2\|F\|_2)\right].$$

Now, $\int_0^{\sigma_n} \sqrt{H(\|F\|_{2,n}/\tau)}\, d\tau \leq \|F\|_{2,n} \int_0^1 \sqrt{H(1/u)}\, du \leq D_H \|F\|_{2,n}$, and therefore, Hölder's inequality followed by Bernstein's gives

(3.9)
$$E\left[\int_0^{\sigma_n} \sqrt{H\left(\frac{\|F\|_{2,n}}{\tau}\right)}\, d\tau\, I(\|F\|_{2,n} > 2\|F\|_2)\right]$$
$$\leq D_H \|F\|_2 \exp\left\{-\frac{9}{8}n\|F\|_2^2\right\} \leq \frac{D_H}{2\sqrt{n}}$$

for the second summand in (3.8). To bound the first summand, note that by concavity of $\int_0^x h(t)\, dt$ when $h \searrow$, and by the properties of $H$, if $\|\sigma_n\|_2 \leq B$,

(3.10)
$$E\left[\int_0^{\sigma_n} \sqrt{H\left(\frac{2\|F\|_2}{\tau}\right)}\, d\tau\, I(\|F\|_{2,n} \leq 2\|F\|_2)\right]$$
$$\leq E \int_0^{\sigma_n \wedge 2\|F\|_2} \sqrt{H\left(\frac{2\|F\|_2}{\tau}\right)}\, d\tau$$
$$\leq \int_0^{\|\sigma_n\|_2 \wedge 2\|F\|_2} \sqrt{H\left(\frac{2\|F\|_2}{\tau}\right)}\, d\tau$$
$$\leq \int_0^{B \wedge 2\|F\|_2} \sqrt{H\left(\frac{2\|F\|_2}{\tau}\right)}\, d\tau$$
$$\leq C_H B \sqrt{H\left(\frac{2\|F\|_2}{B \wedge 2\|F\|_2}\right)}.$$

Taking $B = \|F\|_2$ in (3.10), inequalities (3.8)–(3.10) give the bound $E\|\sum_{i=1}^n \varepsilon_i f(X_i)\| \leq 60(D_H + 2\sqrt{H(2)}C_H)\sqrt{n}\|F\|_2$, hence the first term in the minimum at the right-hand side of (3.4). This bound is accurate only if $\|F\|_2$ is comparable to $\sigma$, and useful only if $\sqrt{n}\|F\|_2$ is not too large. Otherwise, to get the remainder of bound (3.4), we use (2.2) to the effect that

$$\|\sigma_n\|_2^2 \leq \sigma^2 + \left(\frac{8}{\sqrt{n}} E \left\| \frac{1}{\sqrt{n}} \sum_{i=1}^n \varepsilon_i f(X_i) \right\| \right),$$

and take the right-hand side term of this inequality as $B$ in (3.10). Inequality (3.10) with this $B$ then gives, using (3.1),

$$E\left[ \int_0^{\sigma_n} \sqrt{H\left(\frac{\|F\|_{2,n}}{\tau}\right)} d\tau \, I(\|F\|_{2,n} > 2\|F\|_2) \right]$$

$$\leq C_H \sigma \sqrt{H\left(\frac{2\|F\|_2}{\sigma}\right)}$$

$$+ \frac{\sqrt{8}C_H}{\sqrt{n}} \sqrt{E \left\| \sum_{i=1}^n \varepsilon_i f(X_i) \right\|}$$

$$\times \left( \sqrt{H\left(\frac{2\|F\|_2}{\sigma}\right)} \wedge \sqrt{H\left(\frac{2\|F\|_2}{\sqrt{(8/n)E\|\sum_{i=1}^n \varepsilon_i f(X_i)\|} \wedge 2\|F\|_2}\right)} \right).$$

Combining this inequality with inequalities (3.8) and (3.9), and setting $E := E\|\sum_{i=1}^n \varepsilon_i f(X_i)\|$, it follows that either $\quad E \leq 120 D_H \quad$ or $\quad E \leq 360 C_H \sqrt{n} \sigma \sqrt{H\left(\frac{2\|F\|_2}{\sigma}\right)}$

or $\quad E \leq 8 \cdot 120^2 C_H^2 \left[ H\left(\frac{2\|F\|_2}{\sigma}\right) \wedge \left( H\left(\frac{\sqrt{n}\|F\|_2}{\sqrt{2E}}\right) \vee H(1) \right) \right].$

Now the result follows using elementary algebra, upon observing that if $\Psi(x) := x/H(1/\sqrt{x})$, $0 < x \leq 1$, then $\Psi^{-1}(u) \leq u(H(1/\sqrt{u}) \vee 1)$, $0 < u \leq 1/H(1)$. □

It is easy to keep track of the constants in the previous proof, but not necessarily useful.

Several remarks on the previous theorem are in order here.

(1) One may ask for similar inequalities for higher moments. In fact, Theorem 3.1 together with Proposition 3.1 in [23], yields that there exists a constant $C(H)$ that depends on $H$ only through $A_H$, $C_H$ and $D_H$, such that, under the same assumptions as in Theorem 3.1, for all $n \in \mathbf{N}$ and $p \geq 1$,

$$
(3.11) \quad \begin{aligned}E\left\|\sum_{i=1}^n f(X_i)\right\|_{\mathcal{F}}^p \\ \leq C^p(H)\bigg[\bigg\{\bigg(\sqrt{n}\,\sigma\sqrt{H\bigg(\frac{2\|F\|_2}{\sigma}\bigg)}\bigg) \vee H\bigg(\frac{2\|F\|_2}{\sigma} \wedge \frac{\sqrt{n}\|F\|_2}{1440C_H}\bigg)\bigg\}^p \\ \vee p^{p/2}(\sqrt{n}\,\sigma)^p \vee p^p\bigg].\end{aligned}
$$

(2) In (3.3) we could replace $H(\|F\|_{2,n}/\tau)$ by slightly more complicated expressions and the proof of the theorem would still yield sensible bounds; for instance, we show in Example 3.7 that for VC-major classes the right-hand side of (3.3) is of the form $H_1(\|F\|_{2,n}/\tau) + H_2(\|F\|_{2,n}/\tau)\log(A/\|F\|_{2,n})$, with $H_1$ and $H_2$ regularly varying of exponent 1, and with the whole expression monotone in $\|F\|_{2,n}$, and in this case the proof of Theorem 3.1 works with only formal changes.

(3) Note that it is the regular variation of $H$ that allows us to replace the typical entropy integrals by actual entropies in Theorem 3.1. This is significant because it turns out that a partial Sudakov inequality for Rademacher processes due to Talagrand ([29], page 114, Proposition 4.13) allows us to obtain a lower bound for expectations that in some cases is of the same order as the upper bound (3.4). Here is this inequality applied to classes of functions whose absolute values are uniformly bounded by 1:

LEMMA 3.2 (Talagrand). *There exists a universal constant $L$ such that*

$$(3.12) \quad E_\varepsilon\left\|\frac{1}{\sqrt{n}}\sum_{i=1}^n \varepsilon_i f(X_i)\right\| \geq \frac{1}{L}\delta\sqrt{\log N(\mathcal{F}, L_2(P_n), \delta)},$$

*whenever*

$$(3.13) \quad E_\varepsilon\left\|\frac{1}{\sqrt{n}}\sum_{i=1}^n \varepsilon_i f(X_i)\right\| \leq \frac{\sqrt{n}\delta^2}{L}.$$

We will apply this result with $\delta = \sigma/8$. In what follows, the function $H$ satisfies the same conditions as in Theorem 3.1. Also, we set $\|F\|_{2,n} := \|F\|_{L_2(P_n)}$.

DEFINITION 3.3. A class of functions $\mathcal{F}$ that satisfies condition (3.3), that is,

$$\log N(\mathcal{F}, L_2(P_n), \tau) \leq H\bigg(\frac{\|F\|_{2,n}}{\tau}\bigg)$$

for all $\tau > 0$, $n \in \mathbf{N}$ and $\omega \in \Omega$, is full for $H$ and $P$ if there exists $c > 0$ such that

$$\log N(\mathcal{F}, L_2(P), \sigma/2) \geq cH\left(\frac{\|F\|_2}{\sigma}\right). \tag{3.14}$$

THEOREM 3.4. *Let $\mathcal{F}$ satisfy condition (3.3). Assume*

$$n\sigma^2 \geq 2500 \vee \frac{16A_H}{9}, \qquad n\sigma^2 \geq [(672L^2) \vee 1]1920^2 C_H^2 H\left(\frac{6\|F\|_2}{\sigma}\right), \tag{3.15}$$

*where $L$ is the constant in Lemma 3.2. Then,*

$$E\left\|\sum_{i=1}^n f(X_i)\right\|_\mathcal{F} \geq \frac{\sqrt{n}\,\sigma}{32L}\sqrt{\log N(\mathcal{F}, L_2(P), \sigma/2)}. \tag{3.16}$$

*In particular, if a class $\mathcal{F}$ satisfies the entropy bound (3.3) and is full, then, for all $n$ for which conditions (3.15) hold,*

$$\frac{\sqrt{c}}{16L}\sqrt{n}\,\sigma\sqrt{H\left(\frac{\|F\|_2}{\sigma}\right)} \leq E\left\|\sum_{i=1}^n f(X_i)\right\|_\mathcal{F} \tag{3.17}$$

$$\leq 1920 C_H \sqrt{n}\,\sigma\sqrt{H\left(\frac{2\|F\|_2}{\sigma}\right)}.$$

PROOF. By Talagrand's lemma above,

$$E_\varepsilon\left\|\frac{1}{\sqrt{n}}\sum_{i=1}^n \varepsilon_i f(X_i)\right\| \geq \frac{1}{L}\frac{\sigma}{8}\sqrt{\log N(\mathcal{F}, L_2(P_n), \sigma/8)}, \tag{3.18}$$

whenever

$$E_\varepsilon\left\|\frac{1}{\sqrt{n}}\sum_{i=1}^n \varepsilon_i f(X_i)\right\| \leq \frac{\sqrt{n}\sigma^2}{64L}. \tag{3.19}$$

Now we will lowerbound the right-hand side of (3.18) and upperbound the left-hand side of (3.19) with large probability. We start with the right-hand side of (3.18). Let $D := D(\mathcal{F}, L_2(P), \sigma/2)$. By the law of large numbers applied to $D$ functions in $\mathcal{F}$ and to $F$, for all $\varepsilon > 0$ there exists $n$ and $\omega$ such that

$$D(\mathcal{F}, L_2(P), \sigma/2) \leq D(\mathcal{F}, L_2(P_n(\omega)), (1-\varepsilon)\sigma/2)$$
$$\leq N(\mathcal{F}, L_2(P_n(\omega)), (1-\varepsilon)\sigma/4)$$

and

$$\|F\|_{L_2(P_n(\omega))} \leq (1+\varepsilon)\|F\|_2,$$

so that, taking, for example, $\varepsilon = 1/5$, we obtain by (3.3) that

$$D(\mathcal{F}, L_2(P), \sigma) \le e^{H(6\|F\|_2/\sigma)}. \tag{3.20}$$

Let $f_1, \ldots, f_D$ be a maximal set of $\mathcal{F}$ satisfying $P(f_i - f_j)^2 \ge \sigma^2/4$ for all $1 \le i \ne j \le D$. By Bernstein's inequality (e.g., in the form given in [8]), since moreover $P(f_i - f_j)^4 \le 4P(f_i - f_j)^2 \le 16\sigma^2$,

$$\Pr\left\{\max_{1 \le i \ne j \le D}\left(nP(f_i - f_j)^2 - \sum_{k=1}^n (f_i - f_j)^2(X_k)\right) > \tfrac{8}{3}t + \sqrt{32tn\sigma^2}\right\} \le D^2 e^{-t}.$$

Hence, taking $t = \delta n\sigma^2$ and using $P(f_i - f_j)^2 \ge \sigma^2/4$ and (3.20),

$$\Pr\left\{\min_{1 \le i \ne j \le D} \frac{1}{n}\sum_{k=1}^n (f_i - f_j)^2(X_k) \le \sigma^2\left(\frac{1}{4} - \frac{8\delta}{3} - \sqrt{32\delta}\right)\right\} \le e^{2H(3\|F\|_2/\sigma)} e^{-\delta n\sigma^2},$$

which for $\delta = 1/(32 \cdot 8^3)$ gives

$$\Pr\left\{\min_{1 \le i \ne j \le D} P_n(f_i - f_j)^2 \le \frac{\sigma^2}{16}\right\} \le e^{H(6\|F\|_2/\sigma)} e^{-n\sigma^2/(32 \cdot 8^3)}.$$

This implies that the event $A_1$ on which

$$\begin{aligned}
N(\mathcal{F}, L_2(P_n), \sigma/8) &\ge D(\mathcal{F}, L_2(P_n), \sigma/4) \\
&\ge D = D(\mathcal{F}, L_2(P), \sigma/2) \\
&\ge N(\mathcal{F}, L_2(P), \sigma/2)
\end{aligned} \tag{3.21}$$

has probability

$$\Pr(A_1) \ge 1 - e^{H(6\|F\|_2/\sigma) - n\sigma^2/(32 \cdot 8^3)}. \tag{3.22}$$

Under the present hypotheses, (3.7) holds; actually, the proof of Theorem 3.1 gives, before desymmetrizing, that

$$E\left\|\sum_{i=1}^n \varepsilon_i f(X_i)\right\| \le 960 C_H \sqrt{n} \ \sigma\sqrt{H\left(\frac{2\|F\|_2}{\sigma}\right)},$$

in particular giving the right-hand side inequality in (3.17). Therefore, using (2.2) and (3.15), we have

$$\begin{aligned}
E\left\|\sum_{i=1}^n (f^2(X_i) - Pf^2)\right\| &\le n\sigma^2 + E\left\|\sum_{i=1}^n f^2(X_i)\right\| \\
&\le 2n\sigma^2 + 8E\left\|\sum_{i=1}^n \varepsilon_i f(X_i)\right\| \\
&\le 2n\sigma^2 + 4 \times 1920 C_H \sqrt{n} \ \sigma\sqrt{H\left(\frac{2\|F\|_2}{\sigma}\right)} \\
&\le 6n\sigma^2.
\end{aligned}$$

Hence, Bousquet's version of Talagrand's inequality ([10], Theorem 7.3; see also [32]) gives

$$\Pr\left\{\left\|\sum_{i=1}^n (f^2(X_i) - Pf^2)\right\| \geq 6n\sigma^2 + \sqrt{26tn\sigma^2} + t/3\right\} \leq e^{-t},$$

which, taking $t = 26n\sigma^2$, becomes

$$\Pr\left\{\left\|\sum_{i=1}^n (f^2(X_i) - Pf^2)\right\| \geq 41n\sigma^2\right\} \leq e^{-26n\sigma^2}.$$

[Here we could have used Talagrand's inequality (2.1) instead of Bousquet's, but the resulting bound would have been less neat.] So, the event $A_2$ where

$$(3.23) \qquad \left\|\sum_{i=1}^n f^2(X_i)\right\| < 42n\sigma^2$$

has probability

$$(3.24) \qquad \Pr(A_2) > 1 - e^{-26n\sigma^2}.$$

Also, by Bernstein's inequality, as mentioned above, the event

$$(3.25) \qquad A_3 = \{\|F\|_{2,n} \leq 2\|F\|_2\}$$

has probability

$$(3.26) \qquad \Pr(A_3) \geq 1 - \exp\{-\tfrac{9}{4}n\|F\|_2^2\}.$$

Now, on $A_2 \cap A_3$, the usual entropy bound and (3.15) give

$$E_\varepsilon \left\|\frac{1}{\sqrt{n}} \sum_{i=1}^n \varepsilon_i f(X_i)\right\| \leq 120 \int_0^{\sigma_n} \sqrt{H\left(\frac{\|F\|_{2,n}}{\tau}\right)} \, d\tau$$

$$\leq 120 \int_0^{\sqrt{42}\sigma} \sqrt{H\left(\frac{2\|F\|_2}{\tau}\right)} \, d\tau$$

$$(3.27) \qquad \leq 60\sqrt{42} \int_0^{2\sigma} \sqrt{H\left(\frac{2\|F\|_2}{\tau}\right)} \, d\tau$$

$$\leq 120\sqrt{42}\, C_H \sigma \sqrt{H\left(\frac{\|F\|_2}{\sigma}\right)}$$

$$< \frac{\sqrt{n}\, \sigma^2}{64L}.$$

It follows from (3.18)–(3.27) that

$$(3.28) \quad E\left\|\sum_{i=1}^n \varepsilon_i f(X_i)\right\|_\mathcal{F} \geq \frac{\sqrt{n}\sigma}{8L} \sqrt{\log N(\mathcal{F}, L_2(P), \sigma/2)}\, \Pr(A_1 \cap A_2 \cap A_3)$$

and that

$$\Pr(A_1 \cap A_2 \cap A_3) \geq 1 - e^{H(6\|F\|_2/\sigma) - n\sigma^2/(32 \cdot 8^3)} - e^{-26n\sigma^2} - e^{-9n\sigma^2/4}.$$

This last probability is larger than $1/2$ by the inequalities in (3.15). So, integrating in (3.28) and desymmetrizing, we obtain inequality (3.16). The left-hand side of (3.17) now follows from (3.16) and Definition 3.3, proving the theorem. □

Theorem 3.1 recovers and improves on inequalities that go back to Talagrand [42] (see also follow-up work in [20, 21, 35] and, more recently, [22], where only the first and last of these four references use the $L_2$ norm of the envelope in their inequalities). Theorem 3.4 shows that, at least for large $n$, these inequalities are sharp up to constants.

EXAMPLE 3.5. Suppose that $\mathcal{F}$ is *VC-subgraph*, that is,

$$\{\{(x,t): 0 \leq t \leq f(x)\}: f \in \mathcal{F}\} \cup \{\{(x,t): 0 \geq t \geq f(x)\}: f \in \mathcal{F}\}$$

is a VC class of sets. Or, more generally, suppose $\mathcal{F}$ is VC type, that is, there exist $A \geq e$ and $v \geq 1$ such that

$$N(\mathcal{F}, L_2(Q), \tau) \leq \left(\frac{A\|F\|_{L_2(Q)}}{\tau}\right)^v$$

for all $0 < \tau \leq 2\|F\|_{L_2(Q)}$ and all probability measures $Q$, where $F := \sup\{|f|: f \in \mathcal{F}\}$. In this case $H(u) = v \log(Au)$ is slowly varying ($\alpha = 0$) and we can take $C_H = 2$, $D_H = 2A\sqrt{v}/e$ and $A_H = A$.

Since subsets of a VC-subgraph class are also VC-subgraph, this can be applied to the class $\mathcal{F}(tq^{-1}, t] = \{f \in \mathcal{F}: tq^{-1} < \sigma_P f \leq t\}$ with its measurable envelope, say, $F_t$. Define

$$g_q(t) = \left(\frac{A\|F_t\|_2}{t}\right)^v, \qquad 0 < t \leq 1,$$

where $\|F\|_2 := \|F\|_{L_2(P)}$. Then the function $\log g_q(t)$ plays the crucial role in the expectation bound for the class $\mathcal{F}(tq^{-1}, t]$ (which is needed in the inequalities of Section 2 for ratio type suprema). In Sections 4 and 5 below, this function will be involved in conditions for limit theorems about ratio type empirical processes on VC-subgraph classes. It turns out that if $\mathcal{F}$ is a class of indicator functions (i.e., we are dealing with a VC class of sets) such that $PC \leq 1/2$ for all $I_C \in \mathcal{F}$, then

$$A^{-2}(g_q(t))^{2/v} = \frac{\Pr[\bigcup\{C: I_C \in \mathcal{F}, tq^{-1} < \sigma_P(I_C) \leq t\}]}{t^2}$$

is comparable to (in fact, posssibly smaller than) Alexander's [2] *capacity function* $g(t^2)$.

EXAMPLE 3.6. The scope of Theorem 3.1 is much larger than just VC classes. For instance, let $\mathcal{F} = \{f_n := I_{A(n)}/\log(n \vee e), n \in \mathbf{N}\}$, with $A(n) \subseteq [0,1]$ independent for Lebesgue measure and with Lebesgue measure equal to $1/2$ (introduced in [16], proof of Theorem 2.1), and let $P$ be Lebesgue measure on $[0,1]$. Then, $F := 1$ and $\sigma = 1/2$. Also, considering the $L_2(P_n)$ balls centered on the first $m$ functions, with $m$ of the order of $e^{1/\varepsilon}$, it is easy to see that $\log N(\mathcal{F}, L_2(P_n), \varepsilon)$ is of the order of a constant times $1/\varepsilon$, independently of $n$. Then, Theorem 3.1 gives that $E\|\sum_{i=1}^n f(X_i)\| \leq C\sqrt{n}$ for some fixed $c < \infty$, and this is best possible up to constants since $\mathcal{F}$ is $P$-Donsker.

Other classes whose covering numbers are not polynomial include VC-major and VC-hull (see [16] for definitions). We mention the definition of VC-major, that we use below: $\mathcal{F}$ is VC-major if the collection of sets $\{\{s \in S: f(s) \geq t\} : t \in \mathbf{R}, f \in \mathcal{F}\}$ is VC. The following bound on the entropy of such a class is, most likely, new. Note that, as in the case of VC-subgraph classes, it also involves the envelope of the class.

EXAMPLE 3.7. Suppose that $\mathcal{F}$ is a measurable VC-major class of $P$-centered functions whose absolute values are bounded by 1. Our goal will be to show that there exists $A > 0$ such that for all probability measures $Q$ and all $0 < \tau \leq 2\|F\|_{L_2(Q)}$

$$\log N(\mathcal{F}, L_2(Q), \tau) \leq \frac{A\|F\|_{L_2(Q)}}{\tau} \log\left(\frac{A\|F\|_{L_2(Q)}}{\tau}\right) \log\left(\frac{1}{\tau}\right).$$

To this end, take $t_j := (1+\tau)^{-j}$, $j \geq 0$, and let $m(\tau)$ be the smallest $j$ such that $t_j \leq \tau\|F\|_{L_2(Q)}$. Clearly,

$$m(\tau) \asymp \frac{\log(1/(\tau\|F\|_{L_2(Q)}))}{\tau}.$$

For $f \in \mathcal{F}$, define

$$f_\tau := \sum_{j=1}^{m(\tau)} t_j I(t_j < f \leq t_{j-1}).$$

If $t_j < f(x) \leq t_{j-1}$ for $j \leq m(\tau)$, then

$$0 \leq f(x) - f_\tau(x) \leq t_{j-1} - t_j = \tau t_j \leq \tau f(x) \leq \tau F(x).$$

Hence, as soon as $f(x) > \tau\|F\|_{L_2(Q)}$

$$0 \leq f(x) - f_\tau(x) \leq \tau F(x),$$

otherwise

$$0 \leq f(x) - f_\tau(x) \leq f(x) \leq \tau\|F\|_{L_2(Q)}.$$

This implies that
$$\|f - f_\tau\|_{L_2(Q)}^2 \leq \tau^2 \|F\|_{L_2(Q)}^2 + \tau^2 \|F\|_{L_2(Q)}^2 = 2\tau^2 \|F\|_{L_2(Q)}^2.$$

Denote $\mathcal{F}_\tau := \{f_\tau : f \in \mathcal{F}\}$. Since
$$\{(x,t) : f_\tau(x) \geq t\} = \bigcup_{j=1}^{m(\tau)} \{x : f(x) \geq t_{j-1}\} \times (t_j, t_{j-1}]$$

and $\mathcal{F}$ is a VC-major class, the class $\mathcal{F}_\tau$ is VC-subgraph with VC dimension bounded by $Vm(\tau)$ for some $V > 0$. Clearly, $F$ is an envelope of $\mathcal{F}_\tau$ (since $0 \leq f_\tau \leq f$ for all $f \in \mathcal{F}$). Therefore (see, e.g., [47], Theorem 2.6.7), for $\tau > 0$
$$N(\mathcal{F}_\tau; L_2(Q); \tau \|F\|_{L_2(Q)}) \leq \left(\frac{A}{\tau}\right)^{Vm(\tau)},$$

which implies
$$N(\mathcal{F}; L_2(Q); 3\tau \|F\|_{L_2(Q)}) \leq \left(\frac{A}{\tau}\right)^{Vm(\tau)}.$$

Taking into account the bound on $m(\tau)$ and changing variables $\tau \|F\|_{L_2(Q)} \mapsto \tau$, the result follows. Note that the bound can be also written as

$$\log N(\mathcal{F}, L_2(Q), \tau)$$
$$\leq \frac{A\|F\|_{L_2(Q)}}{\tau}\left[\log^2\left(\frac{A\|F\|_{L_2(Q)}}{\tau}\right)\right.$$
$$\left. + \log\left(\frac{A\|F\|_{L_2(Q)}}{\tau}\right)\log\left(\frac{1}{A\|F\|_{L_2(Q)}}\right)\right]$$
$$:= H(\|F\|_{L_2(Q)}, \tau),$$

which can be used in the proof of Theorem 3.1 (with some modifications), to give
$$E\left\|\sum_{i=1}^n f(X_i)\right\|_\mathcal{F} \lesssim (\sqrt{n}\|F\|_{L_2(P)}(1 + \sqrt{\log(A\|F\|_{L_2(P)})^{-1}}))$$
$$\wedge [(\sqrt{n}\sigma\sqrt{H(\|F\|_{L_2(P)}, \sigma)}) \vee H(\|F\|_{L_2(P)}, \sigma) \vee \sqrt{\log n}].$$

EXAMPLE 3.8. As a more specific example, consider the class $\mathcal{F}$ of nondecreasing functions from $[0,1]$ into itself. Obviously, it is a VC-major class. Let $P$ be a nonatomic probability measure on $[0,1]$ and let $G$ be its distribution function. Denote
$$\mathcal{F}_\delta := \{f \in \mathcal{F} : \sigma_P^2 f := Pf^2 \leq \delta^2\},$$

which, of course, is also a VC-major class. An easy computation shows that the envelope of $\mathcal{F}_\delta$ is

$$F_\delta(x) := \sup_{f \in \mathcal{F}, Pf^2 \leq \delta^2} f(x) = \frac{\delta}{\sqrt{P[x,1]}} \vee 1$$

(if $x$ is such that $P[x,1] > \delta^2$, then the supremum in the definition is attained at the function $f_x$ such that $f_x(y) = 0$ for $y < x$ and $f_x(y) = \frac{\delta}{\sqrt{P[x,1]}}$ for $y \geq x$; otherwise, the supremum is equal to 1). Let $x_\delta$ be such that

$$P[x_\delta, 1] = 1 - G(x_\delta) = \delta^2.$$

Then

$$\|F_\delta\|_2^2 = PF_\delta^2 = \delta^2 \int_0^{x_\delta} \frac{P(dx)}{P[x,1]} + \delta^2$$

$$= \delta^2 \int_0^{x_\delta} \frac{dG(x)}{1 - G(x)} + \delta^2 = \delta^2 \int_{\delta^2}^1 \frac{dy}{y} + \delta^2 = \delta^2 \log \frac{e}{\delta^2}.$$

Hence

$$\frac{\|F_\delta\|_2}{\delta} = \sqrt{\log \frac{e}{\delta^2}},$$

and using this together with our bound on the entropy of VC-major classes in Theorem 3.1 yields, by a simple computation,

$$E\|P_n - P\|_{\mathcal{F}_\delta} \lesssim \frac{\delta}{\sqrt{n}} (\log \delta^{-1})^{3/4} (\log \log \delta^{-1})^{1/2}$$

$$\vee \frac{1}{n} (\log \delta^{-1})^{3/2} (\log \log \delta^{-1}) \vee \frac{\sqrt{\log n}}{n}.$$

So, in spite of the fact that the entropy of the class of monotone functions is relatively large, the supremum of the empirical process over the class $\mathcal{F}_\delta$ of "small" monotone functions is of about the same size as for VC-classes of sets due to the small size of localized envelopes.

**4. Ratio limit theorems I: rates when $\phi(t) = t^\alpha$.** In this section we will derive limit theorems a.s. and in probability for general ratio empirical processes, as direct applications of the bounds in Section 2, and we will specialize these to different types of classes of functions, particularly, VC classes, for which we will use the results from Section 3.

4.1. *The case $\phi(t) = t^2$.* We begin with a law of large numbers already in [22], Theorem 6. In this case we take $\sigma_P^2 f = Pf$ (recall that the class $\mathcal{F}$ consists of functions taking values in $[0,1]$). We set

$$\beta_{n,q} := \beta_{n,q,t^2} = \sup_{1 \leq j \leq l_n} \frac{\psi_n(\rho_j)}{\rho_j^2},$$

where $\rho_j = q^j r_n$, $1 \le j \le l_n$, with $l_n = \log \log_q(\delta q / r_n)$.

THEOREM 4.1. *Let $0 < \delta \le 1$ and $r_n \searrow 0$. Let $1 < q_n \le 2$ be a nonincreasing sequence such that $\log(q_n - 1)^{-1} = o(nr_n^2)$. If $nr_n^2 \to \infty$, then the condition*

$$\beta_{n,q_n} \to 0 \tag{4.1}$$

*is necessary and sufficient for*

$$\sup_{\substack{f \in \mathcal{F} \\ r_n^2 < Pf \le \delta}} \left| \frac{P_n f}{Pf} - 1 \right| \to 0 \tag{4.2}$$

*in probability. Moreover, if $1 \ge q_n - 1 \ge (\log n)^{-\delta}$ for some $\delta > 0$, $nr_n^2 / \log \log n \nearrow \infty$ and $\beta_{n,q_n} / \sqrt{n} \searrow 0$, then condition (4.1) is necessary and sufficient for the limit in (4.2) to hold a.s.*

PROOF. The "in probability" part of the theorem follows directly from Proposition 2.8 with $s = s_n \to \infty$ such that $s_n/(nr_n^2) \to 0$ and $s_n / \log(q_n^2 - 1)^{-1} \to \infty$. The "a.s." part follows from Lemma 2.3 together with Proposition 2.8 with $s = s_n = (2+\delta) K q_n \log \log n$. □

The condition $nr_n^2 \to \infty$ is natural in this problem: it certainly is for $\mathcal{F} = \{I_{[0,t]} : 0 \le t \le 1/2\}$, since $\sum_{i=1}^n I_{[0,1/n]}(X_i) - 1 \to_d N - 1$, where $N$ is Poisson 1.

For a specialization of Theorem 4.1 to VC type classes of functions, obtained by replacing $\psi_n(\rho_j)$ in the definition of $E_{n,q}$ by its estimate from Section 3, see Theorem 10 in [22], which recovers classical results and compares to the sufficiency part of Theorem 5.1 of [2] if we restrict to VC classes of sets.

Regarding rates, Proposition 2.8 also gives immediately the following:

THEOREM 4.2. (a) *For $q_n \in (1, 2]$ nonincreasing, and with $\gamma_n := nr_n^2 / \log(q_n - 1)^{-1}$, the sequence*

$$\left( \frac{1}{\beta_{n,q_n}} \wedge \sqrt{\frac{\gamma_n}{1 + \beta_{n,q_n}}} \right. \tag{4.3}$$

$$\left. \wedge \gamma_n \log(\gamma_n(1 + \beta_{n,q_n})^{-1} \vee 2) \right) \sup_{\substack{f \in \mathcal{F} \\ r_n^2 < Pf \le \delta}} \left| \frac{P_n f}{Pf} - 1 \right|, \quad n \in \mathbf{N},$$

*is stochastically bounded.*
   (b) *With $q_n$ as in part (a), if*

$$\sup_n \beta_{n,q_n} < \infty \quad \text{and} \quad \sqrt{\gamma_n} \beta_{n,q_n} \to \infty, \tag{4.4}$$

*then*

(4.5) $$\lim_n \frac{1}{\beta_{n,q_n}} \sup_{\substack{f \in \mathcal{F} \\ r_n^2 < Pf \le \delta}} \left| \frac{P_n f}{Pf} - 1 \right| = 1 \qquad \text{in pr.}$$

Lemma 2.3 and Proposition 2.8 also give almost sure counterparts of Theorem 4.2, that we leave to the reader.

Here is an example showing that normings other than $\beta_n := \beta_{n,q_n}$ do occur in (4.3).

EXAMPLE 4.3. Our object here is to exhibit an example of a class of functions that satisfies

(4.6) $$n r_n^2 \to \infty \quad \text{and} \quad \sqrt{n} r_n \beta_n \to 0$$

for which the sequence

(4.7) $$\frac{1}{\beta_n} \sup_{f \in \mathcal{F}(r_n, \delta]} \left| \frac{P_n f}{Pf} - 1 \right|, \qquad n \in \mathbf{N},$$

is not stochastically bounded, but the sequence

(4.8) $$\sqrt{\gamma_n} \sup_{f \in \mathcal{F}(r_n, \delta]} \left| \frac{P_n f}{Pf} - 1 \right|, \qquad n \in \mathbf{N},$$

is. Let $\varepsilon_{k,j}$ be independent random variables with

$$\Pr\{\varepsilon_{k,j} = 1\} = \frac{1}{j^2} = 1 - \Pr\{\varepsilon_{k,j} = 0\}, \qquad j, k \in \mathbf{N},$$

and

$$X_k = \left( \frac{\varepsilon_{k,j}}{(\log j)^2} : k = 1, 2, \ldots \right),$$

where $\log x := \log(x \vee e)$ [and below, $\log \log x = \log \log(x \vee e^e)$]. The variables $X_k$ are i.i.d. $c_o$-valued r.v.'s.

Let

$$\mathcal{F} = \{f_j(x) = x_j : j \in \mathbf{N}\},$$

where $x_j$ is the $j$th coordinate of $x \in c_o$, so that $P f_j = (j \log j)^{-2}$ and

$$(P_n - P)(f_j) = \frac{1}{n(\log j)^2} \sum_{k=1}^n (\varepsilon_{k,j} - j^{-2})$$

and

$$\frac{P_n f_j}{P f_j} - 1 = \frac{j^2}{n} \sum_{k=1}^n (\varepsilon_{k,j} - j^{-2}).$$

Set
$$r_n = \frac{\log n}{\sqrt{n}}, \qquad \delta = \tfrac{1}{2}.$$

CLAIM 1. *There is a permissible $q_n$ such that*
$$\beta_n = \beta_{n,q_n} \lesssim \frac{\sqrt{\log \log n}}{(\log n)^2}.$$

PROOF. We can take
$$\beta_n := \sup_{u > r_n} \frac{\psi_{n,q_n}(u)}{u^2}$$
$$\simeq \sup_{u \geq r_n} E \sup_{\substack{f \in \mathcal{F} \\ u^2 \leq Pf < u^2 q_n^2}} \frac{|P_n f - Pf|}{u^2}$$

where $2 > q_n \searrow 1$ is such that
$$\log \frac{1}{q_n - 1} = o(nr_n^2) = o((\log n)^2).$$

In fact we take
$$q_n = 1 + \frac{(\log n)^2}{\sqrt{n}}.$$

In order to upperbound $\beta_n$ we note that the number of integers $j$ such that $u^2 \leq Pf < u^2 q_n^2$, that is, such that $(uq_n)^{-1} < j \log j \leq u^{-1}$, $u \geq r_n$, is dominated by
$$\frac{1}{u} - \frac{1}{uq_n} = \frac{q_n - 1}{uq_n} \leq \frac{q_n - 1}{r_n q_n} \leq \log n$$

because if $F(x) = x \log x$, $x > 1$, then $(F^{-1})'(y) < 1$. Moreover, the smallest $j$ in this range, call it $j(u)$, satisfies
$$j(u) \geq \frac{1}{uq_n \log(uq_n)^{-1}} \quad \text{and} \quad j(r_n) \geq \frac{1}{r_n q_n \log(r_n q_n)^{-1}} \geq \frac{\sqrt{n}}{2(\log n)^2}.$$

Bernstein's inequality and Lemma 2.2.10 in [47] (a convexity argument due to Pisier) then give that, for some universal constant $K$,
$$\beta_n \leq \sup_{u > r_n} E\left(\sup_{\substack{f \in \mathcal{F} \\ u^2 \leq Pf < u^2 q_n^2}} \frac{|P_n f - Pf|}{u^2}\right)$$
$$\leq \sup_{u > r_n} \frac{K}{nu^2 (\log(1/(uq_n \log(uq_n)^{-1})))^2}$$

$$\times \left[\frac{1}{3}\log(1+\log n)\right.$$
$$\left. + \sqrt{n}uq_n \log(uq_n)^{-1}\sqrt{\log(1+\log n)}\right].$$

Since this bound is the sup of a decreasing function of $u$, we have

$$\beta_n \leq \frac{5K}{nr_n^2(\log n)^2}\left[\frac{1}{3}\log(1+\log n) + (\log n)^2\sqrt{\log(1+\log n)}\right]$$
$$\leq \frac{6K\sqrt{\log\log n}}{(\log n)^2},$$

at least for all $n$ large enough. Claim 1 is proved. □

It follows from Claim 1 that:

(i) $\beta_n \to 0$, and
(ii) $\sqrt{n}r_n\beta_n \leq \frac{6K\sqrt{\log\log n}}{\log n} \to 0$,

in particular, (4.6) holds. From (i) and Theorem 4.1, we know that

$$\sup_{Pf \geq r_n^2}\left|\frac{P_nf}{Pf} - 1\right| \to 0 \qquad \text{in pr.};$$

in fact, from Theorem 4.2,

$$(\log n)^{1/2}\sup_{Pf \geq r_n^2}\left|\frac{P_nf}{Pf} - 1\right|, \qquad n \in \mathbf{N},$$

is stochastically bounded [note that $\sqrt{\gamma_n} = \sqrt{n}r_n/\sqrt{\log(q_n - 1)^{-1}}$ is of the order of $(\log n)^{1/2}$], so that (4.8) holds. Next we are going to show the following:

CLAIM 2. *For any $\lambda_n \to \infty$, the sequence*

$$\lambda_n(\log n)^{3/2}\sup_{Pf > r_n^2}\left|\frac{P_nf}{Pf} - 1\right|, \qquad n \in \mathbf{N},$$

*converges to infinity in probability, hence so does the sequence*

$$\frac{1}{\beta_n}\sup_{Pf > r_n^2}\left|\frac{P_nf}{Pf} - 1\right|, \qquad n \in \mathbf{N},$$

*by* (ii) *above, $\beta_n = o(\sqrt{\log\log n}/(\log n)^2)$.*

PROOF. Since $j \leq \sqrt{n}/(\log n)^2$ implies $Pf_j > r_n^2$, it follows that

$$\sup_{Pf > r_n^2} \left| \frac{P_n f}{Pf} - 1 \right| \geq \max_{j \leq \sqrt{n}/(\log n)^2} \frac{j^2}{n} \left| \sum_{k=1}^{n} (\varepsilon_{k,j} - j^{-2}) \right|$$
(4.9)
$$> \frac{1}{4(\log n)^4} \max_{\sqrt{n}/2(\log n)^2 < j \leq \sqrt{n}/(\log n)^2} \left| \sum_{k=1}^{n} (\varepsilon_{k,j} - j^{-2}) \right|.$$

Now we estimate this supremum. First we note that, by direct computation (or, e.g., by Hoffmann–Jørgensen's inequality), if $\xi$ is $\text{Bin}(n,p)$ and $np(1-p) > 1$, then there is a universal constant $c < \infty$ such that $E|\xi - np|^3 \leq c(np)^{3/2}$, which, by Berry–Esséen, implies

$$(4.10) \qquad |\Pr\{\xi - np \leq t\sqrt{np(1-p)}\} - \Pr\{g \leq t\}| \leq \frac{C}{\sqrt{n}}$$

for another universal constant $C$, where $g$ is standard normal. Hence, for any $A = A_n > 0$,

$$\Pr\left\{ \max_{\sqrt{n}/2(\log n)^2 < j \leq \sqrt{n}/(\log n)^2} \left| \sum_{k=1}^{n} (\varepsilon_{k,j} - j^{-2}) \right| > A \right\}$$
(4.11)
$$= 1 - \prod_j \left( 1 - \Pr\left\{ \left| \sum_{k=1}^{n} (\varepsilon_{k,j} - j^{-2}) \right| > A \right\} \right)$$
$$\geq 1 - \prod_j \left( 1 - \Pr\left\{ |g| > \frac{A}{\sqrt{n} j^{-1} \sqrt{1 - j^{-2}}} \right\} + \frac{2C}{\sqrt{n}} \right),$$

where the product is over the set of $j$'s such that $\sqrt{n}/2(\log n)^2 < j \leq \sqrt{n}/(\log n)^2$. Now,

$$\frac{A}{\sqrt{n} j^{-1} \sqrt{1 - j^{-2}}} \leq \frac{2A}{(\log n)^2}$$

and, by well-known Gaussian computations, for $2A > (\log n)^2$,

$$\Pr\left\{ |g| > \frac{A}{(\log n)^2} \right\} = \sqrt{\frac{2}{\pi}} \int_{2A/(\log n)^2} e^{-u^2/2} \, du \geq e^{-4A^2/(\log n)^4}.$$

Hence, taking

$$\frac{4A^2}{(\log n)^4} := \frac{1}{3} \log n,$$

we get

$$\Pr\left\{ |g| > \frac{A}{(\log n)^2} \right\} - \frac{C}{\sqrt{n}} \geq \frac{1}{n^{1/3}} - \frac{2C}{\sqrt{n}} = \frac{2c_n (\log n)^2}{\sqrt{n}},$$

with $c_n \to \infty$ [$c_n$ is of the order of $n^{1/6}/2(\log n)^2$]. Replacing this estimate into (4.11), gives

$$\Pr\left\{\max_{\sqrt{n}/2(\log n)^2 < j \leq \sqrt{n}/(\log n)^2}\left|\sum_{k=1}^n (\varepsilon_{k,j} - j^{-2})\right| > \frac{1}{2\sqrt{3}}(\log n)^{5/2}\right\}$$

$$\geq 1 - \left(1 - \frac{2c_n(\log n)^2}{\sqrt{n}}\right)^{\sqrt{n}/(2(\log n)^2)}$$

$$\geq 1 - e^{-c_n} \to 1.$$

Hence, by (4.9),

$$\Pr\left\{(\log n)^{3/2} \sup_{Pf > r_n^2}\left|\frac{P_n f}{Pf} - 1\right| > \frac{1}{8\sqrt{3}}\right\} \to 1,$$

proving Claim 2. □

4.2. *The case* $\phi(t) = t$. Here is a result on convergence in probability and stochastic boundedness of the "normalized" empirical process. It expands Theorem 1 in [22].

THEOREM 4.4. *Let* $\phi(t) = t$, $\delta \leq 1$, $r_n \searrow 0$ *and, for* $1 < q \leq 2$, *let* $\beta_{n,q}$ *denote* $\beta_{n,q,t}$. *Set*

$$\xi_n := \sup_{\substack{r_n < \sigma_P f \leq \delta \\ f \in \mathcal{F}}} \frac{|P_n f - P f|}{\sigma_P f}.$$

*Then the following statements hold:*

(a) *If for all* $q \in (1, \alpha)$ *for some* $\alpha > 1$, $\sqrt{\frac{\log \log 1/r_n}{n}} \vee \frac{1}{nr_n} = o(\beta_{n,q})$, *then* $\frac{\xi_n}{E\xi_n} \to 1$ *in pr.; also, there are sequences* $q_n \searrow 1$ *such that* $\frac{\xi_n}{\beta_{n,q_n}} \to 1$ *in pr.*

(b) *If for some* $q > 1$, $\sqrt{\frac{\log \log 1/r_n}{n}} \vee \frac{1}{nr_n} = O(\beta_{n,q})$, *then the sequences* $\frac{\xi_n}{E\xi_n}$ *and* $\frac{\xi_n}{\beta_{n,q}}$ *are stochastically bounded.*

(c) *If for some* $q > 1$, $\sqrt{\frac{\log \log 1/r_n}{n}} = O(\beta_{n,q})$ *and* $nr_n\beta_{n,q} \to 0$, *then the sequence* $[nr_n \log(1/(nr_n\beta_{n,q}))]\xi_n$ *is stochastically bounded.*

(d) *Let* $\frac{1}{\beta_{n,q}}\sqrt{\frac{\log \log 1/r_n}{n}} \to \infty$ *and* $\frac{1}{nr_n} = O(\beta_{n,q})$ *for some* $q > 1$; *if, moreover,* $r_n \geq \sqrt{\frac{\log \log n}{n}}$, *then* $\sqrt{\frac{n}{\log \log_{q_n} 1/r_n}}\xi_n$ *is stochastically bounded, and otherwise,*

$$\left(\sqrt{\frac{n}{\log \log_{q_n} 1/r_n}} \wedge \frac{nr_n \log((nr_n^2)^{-1} \vee 2)}{\sqrt{\log \log_{q_n} 1/r_n}}\right)\xi_n$$

*is.*

(e) Let $\frac{1}{\beta_{n,q}}\sqrt{\frac{\log\log 1/r_n}{n}} \to \infty$ and $nr_n\beta_{n,q} \to 0$ for some $q > 1$. Then, if $r_n \geq \beta_{n,q}$,

$$\left(\sqrt{\frac{n}{\log\log 1/r_n}} \wedge \frac{nr_n\log((nr_n^2)^{-1} \vee 2)}{\sqrt{\log\log 1/r_n}}\right)\xi_n$$

is stochastically bounded, and otherwise,

$$\left(nr_n\log(1/(nr_n\beta_{n,q})) \wedge \sqrt{\frac{n}{\log\log 1/r_n}} \wedge \frac{nr_n\log((nr_n^2)^{-1} \vee 2)}{\sqrt{\log\log 1/r_n}}\right)\xi_n$$

is.

PROOF. (a) In this case condition (2.15) is satisfied and we can apply inequalities (2.16) and (2.17) with $s = s_n \to \infty$ so that $s = O(\log\log(1/r_n))$ and $t = t_n \to \infty$ so that $t_n = o(nr_n\beta_{n,q})$; then the lower bounds $\tau_{n,q,t}$ for

$$\left|\sup_{\substack{r_n < \sigma_P f \leq \delta \\ f \in \mathcal{F}}} \frac{|P_n f - Pf|}{\phi_q(\sigma_P f)} - \beta_{n,q}\right|$$

in these inequalities are $o(\beta_{n,q})$. Now the result follows because $t \leq \phi_q(t) \leq qt$ and $\beta_{n,q} \leq E_{n,q} \leq \beta_{n,q} + C\tau_{n,q}$ [see the proof of Theorem 2.1 for this last inequality, which holds when the probability in (2.4a) is less than 1].

(b) Follows from similar considerations.

(c) In this case (2.15) is still satisfied (at least up to a multiplicative constant whose only effect is in the multiplicative constants in the probability inequalites). Then, since $n\beta_{n,q}^2 \to \infty$, necessarily $\beta_n > r_n$ from some $n$ on, and inequality (2.17) applies. Under the hypotheses of (c), we have (with $\ll$ signifying "little o")

$$\sqrt{\frac{\log\log_{q_n}(1/r_n)}{n}} \lesssim \beta_n \ll \sqrt{\frac{\beta_n}{nr_n}} \ll \frac{1}{nr_n\log(nr_n\beta_n)^{-1}}$$

so that $t$ times this last term is the dominant one in $\tau_{n,q,t}$, inequality (2.17).

(d) In the first case, $\beta_{n,q} \ll \sqrt{n^{-1}\log\log n} \leq r_n$ and (2.16) applies. Otherwise we must use (2.13) and (2.14); for (2.14), note that if $\beta_{n,q} > r_n$, then

$$\frac{1}{nr_n\log(10 \vee (nr_n\beta_{n,q})^{-1})} \lesssim \sqrt{\frac{\beta_{n,q}}{nr_n}} \lesssim \beta_{n,q} \ll \sqrt{\frac{\log\log n}{n}}.$$

(e) Follows using (2.13) and (2.14), from similar easy considerations. □

A similar result for the a.s. size of $\xi_n$ can be obtained as well. One applies the same principles but makes sure that Lemma 2.3 is satisfied. For instance, direct application of Remark 2.10 and Lemma 2.3 gives that if

$$\frac{\beta_{n,q}}{\sqrt{n}} \searrow, n^{3/2} r_n \nearrow \quad \text{and} \quad \sqrt{\frac{\log\log 1/r_n \vee \log\log n}{n}} \vee \frac{\log\log n}{nr_n} = o(\beta_{n,q})$$

for $q \in (1, \alpha)$, then

$$\limsup_{n \to \infty} \frac{\xi_n}{E\xi_n} \leq 1 \quad \text{a.s.}$$

To show that this lim sup is actually equal to 1, apply Remark 2.10 for $n_k = e^k$ and Borel–Cantelli, as in the second part of the proof of Theorem 2 in [22].

EXAMPLE 4.5. We now modify Example 4.3 to show that the condition

$$\sqrt{\frac{\log\log_{q_n} 1/r_n}{n}} \vee \frac{1}{nr_n} = o(\beta_{n,q})$$

from part (a) of Theorem 4.4 has some degree of sharpness (the problem with absolute sharpness is that we are not using the exact value of $\beta_{n,q}$ to violate the condition, but only an upper estimate), so that our example will satisfy

$$nr_n \beta_{n,q} \leq K < \infty,$$

but it might well be that actually $nr_n \beta_{n,q} \to 0$, which might be too strong a violation of the condition.

We consider

$$\xi_n = \sup_{r_n < \sigma_P f \leq \delta} \frac{|P_n f - Pf|}{\sigma_P f}$$

from Theorem 4.4 with, for example, $\delta = 1/8$, and (see also Section 2)

$$\tilde{\beta}_n = \sup_{u \in (r_n, \delta]} \frac{1}{u} E\left(\sup_{u/q_n < \sigma_P f \leq u} |P_n f - Pf|\right),$$

with $q_n \searrow 1$ and $r_n \searrow 0$. Take $\varepsilon_{k,j}$ and $\mathcal{F}$ as in Example 4.3, and $X_k = (\frac{\varepsilon_{k,j}}{\log j} : j = 2, \ldots) \in c_0$. Since $\text{Var}_P(f_j) = \frac{1}{j^2(\log j)^2} - \frac{1}{j^4(\log j)^2}$ is of the order of $1/(j \log j)^2$ for large $j$ and $\tilde{\beta}_n^{-1} \gg \sqrt{n}$, taking $\sigma_P^2(f) = 1/(j \log j)^2$ will be equivalent to taking $\sigma_P^2(f) = \text{Var}_P(f)$. Now define

$$r_n = \frac{1}{\sqrt{n \log\log n}}, \qquad q_n = 1 + r_n \log n.$$

If, for $u \in (r_n, s]$ we set

$$J_u = \{j : u/q_n < \sigma_P(f_j) \leq u\}$$
$$= \{j : u/q_n < 1/(j \log j) \leq u\},$$

we have

$$j(u) := \min\{j : j \in J_u\} \geq \frac{1}{u \log u^{-1}}$$

and

$$\mathrm{Card}(J_u) \leq \frac{q_n}{u} - \frac{1}{u} \leq \frac{q_n - 1}{r_n} \leq \log n.$$

Using Bernstein and Lemma 2.2.10 in [47], as in Example 4.3, we obtain

$$\frac{1}{u} E\left(\sup_{u/q_n < \sigma_P f \leq u} |P_n f - Pf|\right)$$
$$\leq \frac{K}{nu \log(1/(u \log u^{-1}))} \left[\frac{1}{3} \log \log n + \sqrt{n} u \log u^{-1} \sqrt{\log \log n}\right].$$

Taking the sup over $u \in (r_n, \delta)$ we obtain that

$$\tilde{\beta}_n \leq K \sqrt{\frac{\log \log n}{n}}$$

for some other $K$ and for all $n$ large enough. So, as mentioned, we have

$$nr_n \tilde{\beta}_n \leq K.$$

Now we show that $\xi_n/\tilde{\beta}_n \to \infty$ in probability (actually faster than any rate $A_n$ such that $\frac{A_n}{\sqrt{\log n/\log \log n}} \to 0$). First we observe that

$$2K \frac{\xi_n}{\tilde{\beta}_n} \geq \frac{1}{\log n} \sup_{(1/2)\sqrt{n \log \log n}/\log n \leq j \leq \sqrt{n \log \log n}/\log n} \left|\sum_{k=1}^n (\varepsilon_{k,j} - j^{-2})\right|,$$

and then we easily check, proceeding as in the previous example, that

$$\Pr\left\{2K \frac{\xi_n}{\tilde{\beta}_n} > \frac{1}{4} \sqrt{\frac{\log n}{\log \log n}}\right\} \to 1$$

as $n \to \infty$.

Next, we will consider the case of VC classes of functions, for which we will obtain a result that, although it falls short of recovering the full strength of Theorem 3.1 in [2] when restricted to classes of sets, still gives best possible results up to constants in the classical situation of the uniform empirical c.d.f. in several dimensions, indicators of intervals for the uniform, and halfspaces for the normal (Corollaries 3.5, 3.7 and 3.9 there).

We refer to Example 3.5 for the definition of VC-*subgraph* classes of functions and recall that, by a result of [40], reproduced in [15], if $\mathcal{F}$ is a bounded

VC-subgraph, there exist $A \geq e$ and $v \geq 1$ such that, for every subclass $\mathcal{G} \subseteq \mathcal{F}$, if $G$ is a measurable envelope for $\mathcal{G}$, then

$$N(\mathcal{G}, L_2(Q), \tau) \leq \left(\frac{A\|G\|_{L_2(Q)}}{\tau}\right)^v$$

for all probability measures $Q$ and $0 < \tau \leq 2\|G\|_{L_2(Q)}$. Hence, by Theorem 3.1 there exists a constant $1 \leq K_1 < \infty$ such that if $\mathcal{F}$ is such a class and, moreover, it is suitably measurable and consists of functions taking values in $[0,1]$, then for all $\mathcal{G} \subseteq \mathcal{F}$,

$$E\left\|\sum_{i=1}^n (f(X_i) - Pf)\right\|_{\mathcal{G}}$$

(4.12)
$$\leq K_1 \left[\sqrt{n}\|G\|_2 \wedge \left(\sqrt{n}\sigma_{\mathcal{G}}\sqrt{\log \frac{A\|G\|_2}{\sigma_{\mathcal{G}}}} \right.\right.$$
$$\left.\left. \vee \log\left(\frac{A\|G\|_2}{\sigma_{\mathcal{G}}} \wedge \sqrt{n}\|G\|_2\right) \vee 1\right)\right].$$

In particular this applies to the classes $\mathcal{F}(tq^{-1}, t] = \{f \in \mathcal{F}: tq^{-1} < \sigma_P f \leq t\}$. Letting $F_t$ denote a measurable envelope of $\mathcal{F}(tq^{-1}, t]$, we define (as in Example 3.5)

(4.13)
$$g_q(t) = \left(\frac{A\|F_t\|_2}{t}\right)^v, \qquad 0 < t \leq 1,$$

where $\|F\|_2 := \|F\|_{L_2(P)}$.

Assume that $\sigma_P(f)$ [which is always $\geq \mathrm{Var}_P^{1/2}(f)$] satisfies the following condition:

$$\forall f \in \mathcal{F} \qquad \|f\|_2 \leq C\sigma_P(f)$$

with some constant $C > 0$.

Recall that, given $f_-, f_+ \in L_2(P)$, the set

$$[f_-, f_+] := \{f \in L_2(P): f_- \leq f \leq f_+\}$$

is called an $L_2(P)$-bracket of size (or of order) $\delta > 0$ iff $\|f_+ - f_-\|_2 \leq \delta$. It will be said that $\mathcal{F}$ satisfies the *local bracketing condition* iff there exists a constant $K > 0$ such that for all $f \in \mathcal{F}$ and $0 < \delta \leq \sigma_P(f)/K$ there exists an $L_2(P)$-bracket $[f_-, f_+]$ of size $K\delta$ such that

$$\{g \in \mathcal{F}: \|g - f\|_2 \leq \delta\} \subset [f_-, f_+]$$

[in other words, $L_2(P)$ balls of radius $\delta$ are to be covered by $L_2(P)$-brackets of size $K\delta$].

Given $0 < r < \delta < 1$ and $0 < q \leq 2$, with $\rho_j = rq^j$, $j = 0, 1, \ldots, l = \log_q(\delta q/r)$, we also define

$$(4.14) \qquad w = w(r) = \max_{0 \leq j \leq l}(\log\log_q(\delta q/\rho_j)) \vee (\log g_q(\rho_j)).$$

In the following theorem, we go as far as we can toward extending Alexander's [2] Theorem 3.1 to classes of functions.

THEOREM 4.6. *Let $1 < q \leq 2$ and $r_n \to 0$. Let $\mathcal{F}$ be a VC-subgraph class satisfying the local bracketing condition. Then the following hold with $w_n := w(r_n)$.*

(a) *If $\liminf_n nr_n^2/w_n > 0$ (infinity not excluded), then the sequence*

$$\sqrt{\frac{n}{w_n}} \sup_{\substack{f \in \mathcal{F} \\ r_n < \sigma_P f \leq \delta}} \frac{|P_n f - Pf|}{\sigma_P f}$$

*is stochastically bounded.*

(b) *If $\lim_n nr_n^2/w_n = 0$ and $\log \frac{w_n}{nr_n^2} = O(e^{\tau' w_n})$ for some $\tau' > 0$, then*

$$\frac{nr_n \log(w_n/nr_n^2)}{w_n} \sup_{\substack{f \in \mathcal{F} \\ r_n < \sigma_P f \leq \delta}} \frac{|P_n f - Pf|}{\sigma_P f}$$

*is stochastically bounded.*

(c) *These statements with stochastic boundedness replaced by $\limsup$ finite a.s. (in fact a constant) also hold with $w_n$ changed to $\overline{w}_n = w_n \vee (\log \log n)$ in assumptions and conclusions, under the extra hypothesis that $\overline{w}_n/n^2 \downarrow$ and $\frac{\overline{w}_n}{n^{3/2} r_n \log(\overline{w}_n/(nr_n^2) \vee 2)} \downarrow$.*

To prove the theorem, we start by adapting Theorem 2.1' to this situation, using the bound from Section 3. Let us set up the simplifying notation

$$\mathcal{F}_j := \mathcal{F}(\rho_{j-1}, \rho_j]$$

and denote as $F_j$ a measurable envelope of $\mathcal{F}_j$.

LEMMA 4.7. *Let $\mathcal{F}$ be a (measurable) VC-subgraph class of functions taking values in $[0, 1]$, and let $A$, $v$, $K_1$, $0 < r < \delta$, $0 < q \leq 2$, $\rho_j$, $\mathcal{F}_j$, $l$, $g_q$ and $w$ be as above. Assume further that for each $j$ for which $n\rho_j^2 < w$, $\mathcal{F}_j = \bigcup_{k=1}^{N_j} \mathcal{F}_{j,k}$ and $\mathcal{F}_{j,k}$ has an envelope $F_{j,k}$ satisfying $\|F_{j,k}\|_2 \leq K_2 \rho_j$ for some $K_2 \geq 1$ [i.e., $\mathcal{F}_j$ decomposes into $N_j$ $L_2(P)$-brackets of size of the order of $\rho_j$]. Then, for $s \leq 2K_2^2 w$,*

$$\Pr\left\{ \sup_{\substack{f \in \mathcal{F} \\ r < \sigma_P f \leq \delta}} \frac{|P_n f - Pf|}{\sigma_P f} \geq \frac{qe^2 s I(nr^2 < w)}{nr(1 \vee \log(e^2 s/(K_2^2 nr^2)))} \right.$$

$$+ 2q(17K_1 + K_1K_2 + K_2^2)\sqrt{\frac{2w}{n}}\bigg\}$$

(4.15)
$$\leq Kw \exp\left(\frac{-34K_1 w}{K}\right)$$
$$+ K\left(\max_{j:\,n\rho_j^2 < w} N_j\right)\left(1 + \frac{1}{2}\log_q \frac{w}{nr^2}\right) I(nr^2 < w) \exp\left(\frac{-s}{K}\right).$$

PROOF. We will apply Theorem 2.1′. Set $J = \{1, \ldots, l\}$, $J_1 = \{j \in J : n\rho_j^2 < w\}$ and $J_2 = J \setminus J_1$. For $j \in J_1$, we define

$$\psi_{n,j,k} = \frac{1}{n} E \bigg\| \sum_{i=1}^n (f(X_i) - Pf) \bigg\|_{\mathcal{F}_{j,k}}$$

and

$$V_{n,j,k} = \frac{1}{n} E \bigg\| \sum_{i=1}^n (f(X_i) - Pf)^2 \bigg\|_{\mathcal{F}_{j,k}}$$

which we upperbound by (4.12) and (2.3) as

$$\psi_{n,j,k} \leq \frac{1}{\sqrt{n}} K_1 \|F_{j,k}\|_2 \leq K_1 K_2 \frac{\rho_j}{\sqrt{n}}$$

and

$$V_{n,j,k} \leq \|F_{j,k}\|_2^2 \leq K_2^2 \rho_j^2 := \overline{V}_{n,j,k}.$$

For $j \in J_2$, by (4.12),

$$\frac{\psi_n(\rho_j)}{\rho_j} \leq K_1\left(\sqrt{\frac{w}{n}} \vee \frac{w}{n\rho_j}\right) \leq K_1 \sqrt{\frac{w}{n}}$$

(note that $w \geq 1$) and, by (2.3),

$$V_n(\rho_j) \leq \rho_j^2 + 16\psi_n(\rho_j) \leq 17K_1 \rho_j^2 := \overline{V}_n(\rho_j).$$

Then,

$$\overline{\beta}_{n,q,t} = \max_{\substack{(j,k):\,j \in J_1 \\ k \leq N_j}} \frac{\psi_{n,j,k}}{\rho_j} \vee \max_{j \in J_2} \frac{\psi_n(\rho_j)}{\rho_j} \leq K_1 K_2 \sqrt{\frac{w}{n}}.$$

Now, for $j \in J_2$ we take

$$s_j = 34K_1 w \leq 2n\overline{V}_n(\rho_j)$$

so that the contribution of $J_2$ to $\overline{\tau}_{n,q,t}$ in (2.4′) is just $34K_1\sqrt{2w/n}$ and the contribution to the probability bound, $K/e^{(34K_1/K-1)w}$. For $j \in J_1$, we take $s_{j,k} = s \leq 2K_2^2 w$ if $w \leq nr^2$ and $s_{j,k} = e^2 s$ otherwise. Then the contribution of $\{(j,k): j \in J_1\}$ to $\overline{\tau}_{n,q,t}$ is

$$\frac{e^2 s I(w > nr^2)}{nr \log(e^2 s/(K_2^2 nr^2))} + K_2 \sqrt{\frac{2s}{n}},$$

on account of the fact that $\sqrt{x}/\log x$ is increasing for $x > e^2$, whereas the contribution to the probability bound is dominated by

$$K\left(\max_{j \in J_1} N_j\right)(\operatorname{Card} J_1) e^{-s/(2K)} \leq 2K\left(\max_{j \in J_1} N_j\right)\left(1 + \log_q \frac{w}{nr^2}\right) e^{-s/K}.$$

Collecting bounds, the lemma follows. $\square$

In the previous proof, we could take $s = e^2 K_2 w$ which dominates $e^2 n \overline{V}_{n,j,k}$ for $j \in J_1$, and obtain

$$\Pr\left\{\sup_{\substack{f \in \mathcal{F} \\ r < \sigma_P f \leq \delta}} \frac{|P_n f - Pf|}{\sigma_P f} \geq \frac{2qe^2 K_2^2 w I(nr^2 < w)}{nr \log(e^2 w/(nr^2))} + qK_1(34 + K_2)\sqrt{\frac{2w}{n}}\right\}$$

$(4.15')\quad \leq Kw \exp(-34K_1 w)$

$$+ 2K\left(\max_{j:\, n\rho_j^2 < w} N_j\right)\left(1 + \log_q \frac{w}{nr^2}\right) I(nr^2 < w) \exp(-e^2 K_2 w/K),$$

but in situations when $N_j$ is small (e.g., a constant, as in the case of the uniform empirical c.d.f. in $\mathbf{R}$) we should take $s$ of a smaller order.

Note that replacing $w$ by $cw$, $0 < c < \infty$, in the hypothesis of the previous lemma yields the same conclusion up to constants.

Theorem 4.6 follows at once from $(4.15')$:

PROOF OF THEOREM 4.6. First, for any slice $\mathcal{F}_j$, we construct a partition $\{\mathcal{F}_{j,k}: 1 \leq k \leq N_j\}$, as needed in Lemma 4.7. To this end, consider a minimal covering of $\mathcal{F}_j$ with $L_2(P)$-balls of radius $\rho_j/(Kq)$ and define $\mathcal{F}_{j,k}$ as the intersection of $\mathcal{F}_j$ with the $k$th ball in the partition (if it is empty, discard it). By the bound on the covering numbers of subclasses of a VC-subgraph class, the number $N_j$ will be upperbounded by $(Kq)^v g_q(\rho_j)$, which is in turn upperbounded by $ce^{\tau w(\rho_j)}$ for some $c, \tau > 0$. By the local bracketing condition, for any $k = 1, \ldots, N_j$ there exists an $L_2(P)$-bracket $[f_{j,k,-}, f_{j,k,+}]$ of size $K\rho_j$ covering the class $\mathcal{F}_{j,k}$. If we set $F_{j,k} := f_{j,k,+}$, $F_{j,k}$ becomes an envelope of $\mathcal{F}_{j,k}$ and for arbitrary $f \in \mathcal{F}_{j,k}$

$$\|F_{j,k}\|_2 \leq \|f_{j,k,+} - f_{j,k,-}\|_2 + \|f_{j,k,-}\|_2$$
$$\leq K\rho_j + \|f\|_2 \leq K\rho_j + C\sigma_P(f) \leq (K + C)\rho_j := K_2 \rho_j.$$

Now we are in a position to apply directly inequality $(4.15')$ [together with Lemma 2.3 for part (c)]: increase if necessary the constant $K_2$ so that $K_2/K - \tau - \tau'$ is positive for (a) and (b), and is larger than 1 for (c); for (c) we should also increase $K_1$ so that $34K_1/K > 1$. The a.s. limit is a constant by Borel–Cantelli. $\square$

REMARK 4.8. 1. Suppose that for all $t$ the "slice" $\mathcal{F}(tq^{-1}, t]$ is full for $P$ and $H(u) = v \log Au$ with $\sigma = t$ (recall Definition 3.3 and Theorem 3.4). Then, in the case (a) of Theorem 4.6 and under additional assumption

$$\frac{w_n}{\log \log_q(1/r_n)} \to \infty,$$

we have for some $C > 0$

$$\Pr\left\{C^{-1} \leq \sqrt{\frac{n}{w_n}} \sup_{\substack{f \in \mathcal{F} \\ r_n < \sigma_P f \leq \delta}} \frac{|P_n f - Pf|}{\sigma_P f} \leq C\right\} \to 1 \quad \text{as } n \to \infty.$$

For example, it is clear that the (VC) class $\mathcal{C}$ of all closed (or open) intervals in $[0, 1]$ is full for the uniform distribution and so are any of the slices $\{C \in \mathcal{C} : tq^{-1} < \sqrt{PC} \leq t\}$. Then, $g_q(t) \simeq 1/t^4$. Take $r_n = \sqrt{(\log n)/n}$, which yields $w_n = \overline{w}_n \simeq \log n$, so that Theorem 4.6(b) and (c) give

$$\limsup_{n \to \infty} \sqrt{\frac{n}{\log n}} \sup_{C \in \mathcal{C} : \log n/n \leq P(C) \leq 1/2} \frac{|P_n(C) - P(C)|}{\sqrt{P(C)}} = L < \infty \quad \text{a.s.}$$

Then, the class $\mathcal{C}$ being full, the above limit implies that $L > 0$, a result first obtained by Shorack and Wellner [41] ($L < \infty$), Yukich [49] ($L > 0$) and Alexander ([2], Corollary 3.9, equation (3.12)), where he also obtains it in several dimensions. See also [30].

2. Note that the conclusions of Theorem 4.6 are also true if we only assume (instead of the local bracketing condition) that $\mathcal{F}$ is as in Lemma 4.7, except that now, the bracketing condition of this lemma holds for all $\rho_j$ with $N_j \leq ce^{\tau w(\rho_j)}$ for some $c, \tau > 0$. In principle, the condition $\|F_{j,k}\|_2 \lesssim \rho_j$ can be replaced by weaker assumptions on the local envelopes $F_{j,k}$, which would give rise to different rates.

Alexander [2] does not have an equivalent of the local bracketing assumption in his Theorem 3.1 for VC classes of sets. At this moment, we do not know whether this assumption is needed because of our method (based on combining Talagrand's concentration inequalities and expectation bounds of Section 3), or if it is unavoidable in some form for function classes. However, this assumption holds in all classical examples of classes of sets to which Alexander's Theorem 3.1 applies. Suppose, for instance, that $S = [0, 1]^d$ for some $d \geq 1$ and $P$ has a density that is uniformly bounded and bounded away from 0 on $S$. For $C \subset S$ closed and $\delta > 0$, let $C^\delta$ be the set of all points in $S$ that are within a distance $< \delta$ from $C$ and let $C^{-\delta}$ be the set of all points $x$ such that the closed ball of radius $\delta$ around $x$ is included in $C$. Denote $h$ the Hausdorff distance between closed subsets of $S$, that is,

$$h(C_1, C_2) := \inf\{\delta > 0 : C_1 \subset C_2^\delta, C_2 \subset C_1^\delta\}.$$

Let $\mathcal{C}$ be a VC class of closed *convex* subsets of $S$ such that, for some $K > 0$ and for all $C_0 \in \mathcal{C}$ with $P(C_0) > 0$,

$$K^{-1}h(C, C_0) \leq P(C \triangle C_0) \leq Kh(C, C_0), \qquad C \in \mathcal{C},$$

as soon as $P(C \triangle C_0) < P(C_0)/K$. The upper bound of this inequality always holds for convex sets (see [17], pages 269–270), but the lower bound is satisfied only for special classes of sets (balls, rectangles, etc). Denote $\sigma_P(I_C) := \sqrt{P(C)}$. Then the class $\mathcal{F} := \{I_C : C \in \mathcal{C}\}$ satisfies the local bracketing assumption. The proof easily follows from several simple properties of convex sets described on pages 269–270 of [17]. Indeed, if $C_0 \in \mathcal{C}$ and $0 < \delta < \sqrt{P(C_0)/K}$, then $P(C \triangle C_0) < \delta^2 < P(C_0)/K$ implies that $h(C, C_0) \leq KP(C \triangle C_0) < K\delta^2$. It follows that, for $\sigma = K\delta^2$, $C_0^{-\sigma} \subset C \subset C_0^\sigma$. Hence,

$$\{I_C : \|I_C - I_{C_0}\|_2 = \sqrt{P(C \triangle C_0)} \leq \delta\} \subseteq [I_{C_0^{-\sigma}}, I_{C_0^\sigma}].$$

Since also with some constant $K'$

$$P(C_0^\sigma \setminus C_0^{-\sigma}) \leq K'\sigma,$$

the above inclusion provides a bracket of the size needed in the local bracketing condition. Quite similarly, one can check the condition for VC-subgraph classes of concave (i.e., with a convex subgraph) functions on $[0, 1]$ as well as for some other examples of function classes.

As an illustration, we apply Theorem 4.6 to the uniform empirical c.d.f. in $\mathbf{R}^d$ ([2], Corollary 3.5).

EXAMPLE 4.9 (*The finite-dimensional uniform empirical c.d.f.*). Let $P$ be Lebesgue measure on $[0, 1]^d$, $d \geq 1$, denote by $x^i$ the coordinates of points $\mathbf{x} \in \mathbf{R}^d$, let $\mathcal{F} = \{I_{[0,\mathbf{x}]} : 0 \leq x^i \leq 1, \prod_{i=1}^d x^i \leq 1/2\}$ and take $\sigma_P(I_{[0,\mathbf{x}]}) := (\prod_{i=1}^d x^i)^{1/2}$. Then $\mathcal{F}$ is VC of index $v = d + 1$ ([17], Corollary 4.5.11) so that (4.12) holds with this $v$, and some $A$. It is also easy to see that $\|F_j\|_2^2 = P\{x^1 \cdots x^d \leq \rho_j^2\} \simeq 2^{d-1}\rho_j^2(\log \rho_j^{-1})^{d-1}/(d-1)!$, so that $g(\rho_j) \simeq (\log \rho_j^{-1})^{(d^2-1)/2}$. The local bracketing condition holds by the argument given before the example for convex sets in general. So, we can apply Theorem 4.6 with $w_n \simeq \log \log r_n^{-1}$ and $\overline{w}_n \simeq \log \log n$, assuming, for (c), that $\log \log r_n^{-1}$ is not larger than a constant times $\log \log n$ for all large $n$. The conclusion is: suppose $r_n \to 0$ and $\log \log r_n^{-1}/\log \log n \to c \geq 0$; then,

$$\liminf_n \frac{nr_n^2}{\log \log n} > 0 \implies \limsup_n \sqrt{\frac{n}{\log \log n}} \times \sup_{r_n^2 < \prod_{i=1}^d x_i \leq 1/4} \frac{|F_n(\mathbf{x}) - F(\mathbf{x})|}{\sqrt{\prod_{i=1}^d x_i}} < \infty \qquad \text{a.s.}$$

and, assuming $\log\log n/(n^{3/2}r_n \log(\log\log n/(nr_n^2)))\downarrow$,

$$\lim_n \frac{nr_n^2}{\log\log n} = 0 \implies \limsup_n \frac{nr_n \log(\log\log n/(nr_n^2))}{\log\log n}$$
$$\times \sup_{r_n^2 < \prod_{i=1}^d x_i \leq 1/4} \frac{|F_n(\mathbf{x}) - F(\mathbf{x})|}{\sqrt{\prod_{i=1}^d x_i}} < \infty \quad \text{a.s.}$$

In particular, this last limit allows us to recover the tightness part of a limit theorem of [19], as follows. For dimension $d \geq 2$ and $\varepsilon > 0$, take $r_n = \sqrt{\frac{\varepsilon}{n(\log n)^{d-1}}}$. Then, the last inequality gives

$$\limsup_n \frac{\sqrt{n}}{(\log n)^{(d-1)/2}} \sup_{\varepsilon(n(\log n)^{d-1})^{-1} < \prod_{i=1}^d x_i \leq 1/4} \frac{|F_n(\mathbf{x}) - F(\mathbf{x})|}{\sqrt{\prod_{i=1}^d x_i}} < \infty \quad \text{a.s.}$$

A simple computation shows that if $\xi_i$ are $d$ independent random variables uniform on $[0,1]$, then

$$\Pr\left\{\prod_{i=1}^d \xi_i \leq \frac{\varepsilon}{n(\log n)^{d-1}}\right\} \simeq \frac{\varepsilon}{(d-1)!},$$

which, by another simple computation, allows us to conclude from the previous limit that the sequence

$$\frac{\sqrt{n}}{(\log n)^{(d-1)/2}} \sup_{\prod_{i=1}^d x_i \leq 1/4} \frac{|F_n(\mathbf{x}) - F(\mathbf{x})|}{\sqrt{\prod_{i=1}^d x_i}}$$

is stochastically bounded.

EXAMPLE 4.10 (Example 2.7, revisited). Theorem 4.6 essentially does not distinguish between stochastic and a.s. boundedness for the empirical c.d.f. However, Lemma 4.7 does when $d = 1$. For $d = 1$ we can take $s$ of the order of $\log(w_n/nr_n^2)$ since $N_j$ is constant (as we did in Example 2.7). If $r_n \gtrsim \log\log\log n/\sqrt{n\log\log n}$, then $w_n \simeq \log\log n$, and $s/nr_n$ is dominated by $\sqrt{\frac{\log\log n}{n}}$, so that Lemma 4.7 shows that the sequence

$$\sqrt{\frac{n}{\log\log n}} \sup_{r_n^2 < x \leq 1/2} \frac{|F_n(x) - x|}{\sqrt{x}}$$

is stochastically bounded. Then, since $\Pr\{\min_{i\leq n} X_i \leq \varepsilon/n\} \leq \varepsilon$, we get that, for $d = 1$, the sequence

$$\sqrt{\frac{n}{\log\log n}} \sup_{0 < x \leq 1/2} \frac{|F_n(x) - x|}{\sqrt{x}}$$

is stochastically bounded. By the limiting result of Eicker [18] (see also [14, 24]) this rate is exact.

4.3. *The cases $\phi(t) = t^\alpha$, $\alpha \neq 1, 2$.* For $\alpha \in (1, 2)$ direct application of Proposition 2.11 gives the following analogue of Theorems 4.1 and 4.2:

THEOREM 4.11. *Let $\alpha \in (1, 2)$, $0 < \delta \leq 1$ and $r_n \downarrow 0$, and set*

$$\xi_n := \sup_{\substack{r_n < \sigma_P f \leq \delta \\ f \in \mathcal{F}}} \frac{|P_n f - Pf|}{(\sigma_P f)^\alpha}.$$

(a) *If $nr_n^\alpha \to \infty$, then the condition $\beta_{n,q} := \beta_{n,q,t^\alpha} \to 0$ for some $q > 1$ is necessary and sufficient for $\xi_n \to 0$ in probability (a.s., if we also have $nr_n^\alpha / \log \log n \to \infty$).*

(b) *If $\sup_n \beta_{n,q} < \infty$ and $nr_n^\alpha \beta_{n,q} \to \infty$ for all $q \in (1, 1 + \delta)$ for some $\delta > 0$, then*

$$\frac{\xi_n}{E\xi_n} \to 1$$

*in probability (and the convergence is a.s. if $nr_n^\alpha \beta_{n,q} / \log \log n \to \infty$).*

(c) *For any $q > 1$, the sequence*

$$\left(\frac{1}{\beta_{n,q}} \wedge nr_n^\alpha \log((nr_n^\alpha (r_n^{2-\alpha} \wedge \beta_{n,q})^{-1}) \vee 2) \wedge \sqrt{\frac{nr_n^\alpha}{r_n^{2-\alpha} \vee \beta_{n,q}}}\right) \xi_n$$

*is stochastically bounded.*

For VC classes of functions, adapting the proofs of Lemma 4.7 and Theorem 4.6 to the case of $\alpha \in (1, 2]$ only gives the obvious: for instance, that if in Theorem 4.6 we replace $\sigma_P f$ by $(\sigma_P f)^\alpha$ in the displays, then multiplying by $r_n^{\alpha-1}$ the corresponding expressions produces sequences that are stochastically bounded [or a.s. bounded in part (c)]. This observation applies, for instance, to give the tightness part of the remaining cases in [19], namely, that, just as in Example 4.9, for $1/2 \leq \nu \leq 1$ (the case $\nu = 1/2$ is covered by that example), the sequence

$$\frac{n^{1-\nu}}{(\log n)^{\nu(d-1)}} \sup_{\prod_{i=1}^d x_i \leq 1/4} \frac{|F_n(\mathbf{x}) - F(\mathbf{x})|}{\sqrt{\prod_{i=1}^d x_i}}$$

is stochastically bounded, where we assume $d \geq 2$. Extensions of Example 4.10 to powers of $x$ different from $1/2$ are equally easy to get in the case $d = 1$ (they are omitted).

For $\alpha \in (0, 1)$, we make the rates explicit only under condition $(2.15')$ and the result is a direct consequence of Proposition 2.12.

THEOREM 4.12. *Let $\alpha \in (0,1)$, $0 < \delta \leq 1$ and $r_n \downarrow 0$, set*

$$\xi_n := \sup_{\substack{r_n < \sigma_P f \leq \delta \\ f \in \mathcal{F}}} \frac{|P_n f - P f|}{(\sigma_P f)^\alpha},$$

*and assume*

$$r \vee \beta_{n,q}^{1/(2-\alpha)} \geq \sqrt{\frac{3K \log \log_q(q^2 \delta/r)}{n}}$$

*for all $q \in (1, \tau)$, for some $\tau > 1$. Then:*

(a) *if $nr_n^\alpha \to \infty$, then the condition $\beta_{n,q} \to 0$ for some $q > 1$ is necessary and sufficient for $\xi_n \to 0$ in pr.;*

(b) *if $\sqrt{n} \beta_{n,q} \to 0$ for $q \in (1, \tau)$, then $\xi_n / E \xi_n \to 1$ in pr., and there are sequences $q_n \searrow 1$ for which $\xi_n / \beta_{n,q_n} \to 1$ in pr.;*

(c) *if $\beta_{n,q} \leq C(n^{-1/2} \wedge r_n^{2-\alpha})$ for some $C < \infty$ and $q > 1$, then the sequence $\sqrt{n} \xi_n$ is stochastically bounded, and*

(d) *if, for some $0 < C < \infty$ and $q > 1$, $r_n^{2-\alpha} \leq C \beta_{n,q} \leq n^{-1/2}$, then the sequence*

$$(\sqrt{n} \wedge nr_n^\alpha \log((nr_n^\alpha \beta_{n,q}) \vee 2) \wedge \sqrt{nr_n^\alpha \beta_{n,q}^{-1}}) \xi_n$$

*is stochastically bounded.*

For VC classes of functions, one obtains analogues of Lemma 4.7 and Theorem 4.6 for $\alpha \in (0,1)$ as follows: under the hypotheses of Lemma 4.7, the bound in (4.15) holds for the probability

$$\Pr\left\{ \sup_{\substack{f \in \mathcal{F} \\ r < \sigma_P f \leq \delta}} \frac{|P_n f - P f|}{(\sigma_P f)^\alpha} \geq C(K_1, K_2, q, \alpha) \right.$$

$$\left. \times \left[ \frac{s I(nr^2 < w)}{nr^\alpha (1 \vee \log(s/(nr^2)))} + \delta^{1-\alpha} \sqrt{\frac{w}{n}} \right] \right\}.$$

And under the hypotheses of Theorem 4.6, except that we replace $1/2 \geq \delta > r_n$ by $1/2 \geq \delta_n > r_n$, we have:

(a) *if $\liminf_n n(r_n^\alpha \delta_n^{1-\alpha})^2 / w_n > 0$ (infinity not excluded), then the sequence*

$$\frac{1}{\delta_n^{1-\alpha}} \sqrt{\frac{n}{w_n}} \sup_{\substack{f \in \mathcal{F} \\ r_n < \sigma_P f \leq \delta}} \frac{|P_n f - P f|}{(\sigma_P f)^\alpha}$$

*is stochastically bounded;*

(b) if $\lim_n n(r_n^\alpha \delta_n^{1-\alpha})^2/w_n = 0$ and $\log \frac{w_n}{nr_n^2} = O(e^{\tau' w_n})$ for some $\tau' > 0$, then

$$\frac{nr_n \log(w_n/nr_n^\alpha)}{w_n} \sup_{\substack{f \in \mathcal{F} \\ r_n < \sigma_P f \leq \delta}} \frac{|P_n f - Pf|}{(\sigma_P f)^\alpha}$$

is stochastically bounded;

(c) the corresponding statements for asymptotic a.s. boundedness under monotonicity conditions analogous to those in Theorem 4.6(c).

**5. Ratio limit theorems II: asymptotic continuity moduli and weighted central limit theorems.** These two types of limit theorems usually involve functions of the form $\phi(t) = tL(1/t)$ where $L$ is nondecreasing and slowly varying at infinity.

5.1. *Local and global moduli.* Local asymptotic moduli in probability for general classes of functions were already treated in [22], Theorems 4, 5 and 9. Here we will only derive an a.s. general result which is the companion to Theorem 4 in the just mentioned reference. As usual, $\mathcal{F}$ is a measurable class of functions taking values on $[0,1]$.

Following Alexander [2], a local asymptotic modulus of the empirical process over $\mathcal{F}$ at 0 is an increasing function $\omega$ for which there exist $r_n < \delta_n < 1$ both nonincreasing, with $\sqrt{n}\delta_n$ nondecreasing such that

$$(5.1) \qquad \limsup_n \sup_{\substack{r_n < \sigma_P f \leq \delta_n \\ f \in \mathcal{F}}} \frac{|\nu_n(f)|}{\omega(\sigma_P f)} < \infty \qquad \text{a.s.,}$$

where

$$\nu_n(f) := \sqrt{n}(P_n f - Pf)$$

is the empirical process indexed by $\mathcal{F}$.

Although the results below are under our general assumption that the functions in the class $\mathcal{F}$ take values in $[0,1]$, in the case when $\sigma_P(f) := \sqrt{\operatorname{Var}_P(f)}$ and hence $\sigma_P(f+c) = \sigma_P(f)$, $\sigma_P(cf) = |c|\sigma_P(f)$ for any constant $c$, a simple rescaling allows one to deal also with arbitrary uniformly bounded classes of functions. This is of importance in the case of global moduli. A global asymptotic modulus of the empirical process over $\mathcal{F}$ is any local modulus for $\mathcal{F}' = \{f - g : f, g \in \mathcal{F}\}$ at 0, with $\delta_n \gg \sqrt{n^{-1}\log n}$.

THEOREM 5.1. *Let $q > 1$, $r_n < \delta_n < 1$ both nonincreasing, with $\sqrt{n}\delta_n$ nondecreasing, and let $\omega$ be a bounded nondecreasing function on $[0,1]$ such that*

$$\omega(t) \geq \sqrt{n}\psi_{n,q}(t), \qquad t \in [r_n, \delta_n],$$

*for all $n$, and satisfying that $\omega(u)/u \downarrow$,*

$$\sup_n \frac{\delta_n \sqrt{\log\log n \vee \log\log_q(\delta_n q^2/r_n)}}{\omega(\delta_n)} < \infty,$$

$$\sup_n \frac{\log\log n \vee \log\log_q(\delta_n q^2/r_n)}{\sqrt{n}\omega(r_n)} < \infty$$

*and that these two sequences decrease when divided by $n$. Then, the limit (5.1) holds.*

PROOF. We apply Theorem 2.1 and Lemma 2.3. Let $K$, which we can assume to be larger than 1, be as in (2.4a). We take [see (2.3)]

$$\overline{V}_{n,q}(\rho_j) = L(\rho_j^2 + 16\omega(\rho_j)/\sqrt{n}) \geq L(\rho_j^2 + 16\psi_{n,q}(\rho_j)),$$

for $j = 1, \ldots, \ell_n$, where $\ell_n$ is the smallest integer $j$ such that $\rho_j = r_n q^j \geq \delta_n$, and where $L$ is the largest of $K$ and the second supremum above. Then, if we take $s_j = 2K(\log\log n + \log\log_q(\delta_n q^2/\rho_j))$, we have $s_j \leq 2n\overline{V}_{n,q}(\rho_j)$, and inequality (2.4a) directly gives

$$\Pr\left\{q^{-1} \sup_{\substack{f \in \mathcal{F} \\ r_n < \sigma_P f \leq \delta_n}} \frac{\nu_n(f)}{\omega(\sigma_P f)} \geq 1 + 2\max_j \sqrt{\frac{s_j \overline{V}_{n,q}(\rho_j)}{\omega^2(\rho_j)}}\right\} \leq \frac{K}{(\log n)^2}.$$

Now the theorem follows from Lemma 2.3 and the hypotheses on $\omega$. □

Let $F_{q,u} = \sup\{|f| : f \in \mathcal{F}, u/q < \sigma_P f \leq u\}$, $1 < q \leq 2$, $0 < u < 1$, be the local envelopes for $\mathcal{F}$, and define $g_q(r)$ as any nonincreasing function satisfying $q\|F_{q,r}\|_{L_2(P)}/r \leq g_q(r) \leq q/r$. By proceeding as in Theorem 9 in [22], Theorem 5.1 gives that, for any bounded VC class of functions, the function $\omega_1(t) = t\sqrt{\log\log(1/t) + \log g_q(t)}$ is a local asymptotic modulus at 0 and that the function $\omega_0(t) = t\sqrt{\log(1/t)}$ is a global modulus, thus generalizing Theorem 4.1 of [2] to classes of functions (and demonstrating the same difference between local and global continuity moduli as in the classical cases of Brownian motion, Brownian bridge and univariate empirical process). For the global modulus, one takes $g_q(u) = q/u$.

5.2. *Central limit theorems.* We consider here weighted CLTs for empirical processes in the spirit of Alexander [3]. Let $\psi$ be a strictly increasing continuous function such that

(5.2) $$\psi(0) = 0 \quad \text{and} \quad \lim_{t \to 0} \frac{\psi(t)}{t} = \infty.$$

We call such a function a *weight*. We will find conditions on $\psi$ and a decreasing sequence $r_n$ so that, for a $P$-Donsker class of uniformly bounded functions $\mathcal{F}$, we have

$$\frac{\nu_n(f)}{\psi(\sigma_P f)} I(\sigma_p f > r_n) \xrightarrow{\mathcal{L}} \frac{G_P(f)}{\psi(\sigma_P f)}$$

in $\ell_\infty(\mathcal{F} \setminus \mathcal{F}_0)$ and the limiting process $G_P/\psi \circ \sigma_P$ is sample continuous on $\mathcal{F} \setminus \mathcal{F}_0$ for the pseudo distance $d_P(f,g) = \sigma_P(f-g)$, where $\mathcal{F}_0 := \{f \in \mathcal{F} : \sigma_P f = 0\}$. (For definitions of $P$-Donsker or CLT$(P)$ classes, pre-Gaussian classes, and others associated to uniform central limit theorems, see, e.g., [17] or [47].)

We need to comment on condition (5.2). For classes of sets, this condition is necessary for $G_P/\psi \circ \sigma_P$ to be a.s. in $\ell_\infty$ (Lemma 5.1 in [4]) but this is not so for classes of functions: just consider $\mathcal{F} = \{\alpha f : 0 < \alpha \leq 1\}$ for some bounded function $f$. Then, $G_P(\alpha f)/\sigma_P(\alpha f)$ does not depend on $\alpha$ and the sample paths are just constants. However, if the class $\mathcal{F}$ is sufficiently rich, then (5.2) is also necessary; for instance, assume that $\mathcal{F}$ is *convex* and symmetric (i.e., $f_i \in \mathcal{F}$ and $\sum_{\text{finite}} |\lambda_i| \leq 1$ implies $\sum \lambda_i f_i \in \mathcal{F}$), and that the subspace of $L_2(\Omega)$ generated by the process $G_P(f)$, $f \in \mathcal{F}$, is infinite dimensional (if it were finite dimensional, we would be in the case of the finite-dimensional central limit theorem). Then, by Gram–Schmidt orthogonalization, there exists an infinite sequence of functions $f_i$ in $\mathcal{F}$ such that $\sigma_P(f_i) \neq 0$ and $EG_P(f_i)G_P(f_j) = 0$ if $i \neq j$. But then $G_P(f_i)/\sigma_P(f_i)$ are i.i.d. $N(0,1)$ and their sup is infinite with probability 1. We are thus justified in assuming condition (1) for our weights.

Another useful remark is the following:

LEMMA 5.2. *Assuming (5.2) and $\mathcal{F}$ $P$-pre-Gaussian, if $G_P/\psi \circ \sigma_P$ is $d_P$ sample continuous on $\mathcal{F} \setminus \mathcal{F}_0$ (meaning that it has a version with bounded and $d_P$-uniformly continuous sample paths), then*

$$\lim_{f \in \mathcal{F} \setminus \mathcal{F}_0, \sigma_P f \to 0} \frac{G_p(f)}{\psi(\sigma_P f)} = 0 \qquad a.s.$$

PROOF. If $\sigma_P(f)$, $f \in \mathcal{F} \setminus \mathcal{F}_0$, is bounded away from zero, then there is nothing to prove. Otherwise, let $f_n \in \mathcal{F}$ be such that $\sigma_P f_n \to 0$. Then, the sequence $G_P(f_n)/\psi(\sigma_P f_n)$ is a.s. Cauchy by hypothesis, and since $E(\frac{G_P(f_n)}{\psi(\sigma_P f_n)})^2 \to 0$ by (5.2), it also converges to zero in probability. Hence, this sequence converges to zero a.s. $\square$

The following proposition, which is analogous to Theorem 4.2 in [4], will allow use of the inequality in Theorem 2.1. From now on we will assume without loss of generality that $\mathcal{F}_0$ *is empty* and that the functions in $\mathcal{F}$ take values in $[0,1]$.

LEMMA 5.3. *Let $\mathcal{F}$ be a measurable class of functions, let $\psi$ be a weight function as defined above and let $r_n \to 0$, $r_n > 0$. Then,*

$$\frac{\nu_n(f)}{\psi(\sigma_P f)} I(\sigma_p f > r_n) \overset{\mathcal{L}}{\to} \frac{G_P(f)}{\psi(\sigma_P f)}$$

*in $\ell_\infty(\mathcal{F})$ and the limiting process $G_P/\psi \circ \sigma_P$ is $d_P$ sample continuous on $\mathcal{F}$ if and only if both*

$$\mathcal{F}_{\geq r} := \{f \in \mathcal{F} : \sigma_P(f) \geq r\} \text{ is } P\text{-Donsker}$$

*and*

(5.3) $$\lim_{\delta \to 0} \limsup_n \Pr\left\{\sup_{\substack{f \in \mathcal{F} \\ r_n < \sigma_P f \leq \delta}} \frac{|\nu_n(f)|}{\psi(\sigma_P f)} > \varepsilon\right\} = 0.$$

PROOF. If the weighted processes converge in law and the limit is $d_P$ sample continuous, then, by the continuous mapping theorem, $\mathcal{F}_{\geq r}$ is $P$-Donsker. Also, by the portmanteau lemma,

$$\limsup_n \Pr\left\{\sup_{\substack{f \in \mathcal{F} \\ r_n < \sigma_P f \leq \delta}} \frac{|\nu_n(f)|}{\psi(\sigma_P f)} \geq \varepsilon\right\} \leq \Pr\left\{\sup_{f \in \mathcal{F}, \sigma_P f \leq \delta} \frac{|G_P(f)|}{\psi(\sigma_P f)} \geq \varepsilon\right\},$$

which, by Lemma 5.2, tends to zero as $\delta \to 0$. The direct part follows as in [3]. □

THEOREM 5.4. *Let $r_n \to 0$, $0 < r_n < 1/2$, and let $\psi$ be a weight function such that $\sup_{0 < x \leq 1/2} \psi(2x)/\psi(x) = C < \infty$. Assume*

(5.4) $$\lim_{\delta \to 0} \limsup_n \sup_{r \in (r_n, \delta]} \frac{r\sqrt{\log \log_{q_n} 1/r}}{\psi(r)} = 0$$

*and*

(5.5) $$\frac{\log \log_{q_n} 1/r_n}{\psi(r_n)\sqrt{n}} \to 0 \qquad \text{as } n \to \infty,$$

*where $2 \geq q_n \searrow 1$ or $q_n \equiv c$. Then, the conditions $\mathcal{F}_{\geq r} \in CLT(P)$ for all $r > 0$ and*

(5.6) $$\lim_{\delta \to 0} \limsup_n \sup_{r \in (r_n, \delta]} \frac{\sqrt{n}\psi_{n,q_n}(r)}{\psi(r)} = 0$$

*are necessary and sufficient for the process $\frac{G_P(f)}{\psi(\sigma_P f)}$, $f \in \mathcal{F}$, to be $d_P$ sample continuous and for*

(5.7) $$\frac{n^{1/2}(P_n - P)(f)}{\psi(\sigma_P(f))} I(\sigma_P(f) > r_n) \overset{\mathcal{L}}{\to} \frac{G_P(f)}{\psi(\sigma_P(f))} \qquad \text{in } \ell_\infty(\mathcal{F}).$$

PROOF. By Lemma 5.3, the proof is basically the same as that of Theorem 5.1. Define

$$\varepsilon(n,\delta) = \frac{\log\log_{q_n} 1/r_n}{\psi(r_n)\sqrt{n}} \vee \sup_{r \in (r_n, \delta q_n)} \frac{\sqrt{n}\psi_{n,q_n}(r)}{\psi(r)}$$

which, by (5.5) and (5.6), satisfies $\lim_{\delta \to 0} \limsup_n \varepsilon(n,\delta) = 0$, and then,

$$\overline{V}_{n,q_n}(\rho_j) = K[\rho_j^2 + 16\varepsilon(n,\delta)\psi(\rho_j)/\sqrt{n}]$$
$$\geq K[\rho_j^2 + 16\psi_{n,q_n}(\rho_j)],$$

where $K \geq 1$ comes from (2.4a), and where, as usual, $\rho_j = r_n q_n^j$ so that it depends on $n$ even if we do not show it. Note that $\overline{V}_{n,q_n}(\rho_j)$ is admissible in Theorem 2.1 by (2.3). Now set $s_j = 2K[t(n,\delta) + \log\log_{q_n}(\delta q^2/\rho_j)]$, with $t(n,\delta) = \min[\log\log_{q_n} r_n^{-1}, \inf_{r_n \leq r \leq \delta q_n} \psi(r)/r]$, which satisfies

(5.8) $$\lim_{\delta \to 0} \liminf_n t(n,\delta) = \infty,$$

because $r_n \to 0$ and $\psi$ is a weight function. Note also that, by the hypotheses, $s_{n,j} \leq 2n\overline{V}_{n,q_n}(\rho_j)$ for all $1 \leq j \leq \ell_n(\delta)$, where $\ell_n(\delta)$ is the smallest integer $j$ such that $r_n q^j \geq \delta$. Therefore, Theorem 2.1 or Corollary 2.2 gives

$$\Pr\left\{C^{-1} \sup_{\substack{f \in \mathcal{F} \\ r_n < \sigma_P f \leq \delta}} \frac{|\nu_n(f)|}{\psi(\sigma_P f)} \geq \sup_j \frac{\sqrt{n}\psi_{n,q_n}(\rho_j)}{\psi(\rho_j)}\right.$$
$$\left. + 2\max_j \sqrt{\frac{s_j \overline{V}_{n,q_n}(\rho_j)}{\psi^2(\rho_j)}}\right\} \leq Ke^{-2t(n,\delta)}.$$

Now, condition (5.6) implies that $\lim_{\delta \to 0} \limsup_n \sup_j \frac{\sqrt{n}\psi_{n,q_n}(\rho_j)}{\psi(\rho_j)} = 0$ and moreover, since

$$\frac{1}{2K^2} \frac{s_j \overline{V}_{n,q_n}(\rho_j)}{\psi^2(\rho_j)} \leq \frac{\rho_j}{\psi(\rho_j)} + \frac{32\varepsilon(n,\delta)\log\log_{q_n} r_n^{-1}}{\sqrt{n}\psi(r_n)}$$
$$+ \frac{\rho_j^2 \log\log_{q_n}(\delta q^2/\rho_j)}{\psi^2(\rho_j)},$$

equations (5.4) and (5.5) imply

(5.9) $$\lim_{\delta \to 0} \limsup_n \max_j \sqrt{\frac{s_j \overline{V}_{n,q_n}(\rho_j)}{\psi^2(\rho_j)}} = 0.$$

Therefore, (5.3) holds, and by Lemma 5.3, so does (5.7). Conversely, if (5.6) does not hold, while we still have (5.8) and (5.9), the term $\limsup_n \sup_j \frac{\sqrt{n}\psi_{n,q_n}(\rho_j)}{\psi(\rho_j)}$ stays bounded away from zero for a sequence $\delta_k \to 0$, and this implies, by the second inequality in Corollary 2.2, that (5.3) does not hold, and therefore, by Lemma 5.3, neither does (5.7). □

In order to apply the above theorem, one needs to have reasonable estimates of $\psi_{n,q}(r)$, and it is here where the results in Section 3 may become useful.

In the case of the classical VC-subgraph classes (the uniform empirical distribution function, indicators of intervals for the uniform law on the unit cube, or half-spaces for the normal), the above theorem does not give best possible results, just as in the case of $\phi(x) = x$ in Section 4.2 [see (4.6)–(4.8)]. As in that section, we will prove a theorem that handles these cases, but we will only apply it to the multidimensional empirical c.d.f.

For a VC-subgraph class $\mathcal{F}$, or more generally for a VC type class [i.e., one such that for any $\mathcal{G} \subseteq \mathcal{F}$ and any probability measure $Q$, $N(\mathcal{G}, L_2(Q), \tau) \leq (A\|G\|_{L_2(Q)}/\tau)^v$ for some $A \geq e$, $v \geq 1$ and all $0 < \tau \leq 2\|G\|_{L_2(Q)}$], and given $0 < r < \delta < 1$ and $1 < q \leq 2$, let $g_q(t)$ and $w$ be as defined in (4.13) and (4.14).

THEOREM 5.5. *Let $\mathcal{F}$ be a VC-subgraph class satisfying the local bracketing condition and let $r_n \to 0$, $q \in (1, 2]$. Let*

$$(5.10) \qquad \phi(t) = tL(1/t), \qquad 0 \leq t \leq 1,$$

*with $L(u) \nearrow \infty$ as $u \nearrow \infty$ and $u^\tau L(1/u)$ nondecreasing for some $0 < \tau < 1$. Assume $\mathcal{F}_{\geq r} \in CLT(P)$ for all $r > 0$. Then, the conditions*

$$(5.12) \quad \lim_{u \to 0} \frac{\sqrt{\log g_q(u)}}{L(1/u)} = 0 \quad \text{and} \quad \lim_{n \to \infty} \frac{w_n}{\sqrt{n} r_n L(1/r_n) \log(w_n/nr_n^2)} = 0,$$

*where $w_n = w(r_n)$, imply that the Gaussian process $G_P(f)/\phi(\sigma_P f)$, $f \in \mathcal{F}$, is sample continuous and*

$$\frac{\nu_n(f)}{\phi(\sigma_P f)} I(\sigma_P f > r_n) \overset{\mathcal{L}}{\to} \frac{G_P(f)}{\phi(\sigma_P f)} \qquad \text{in } \ell_\infty(\mathcal{F}).$$

For the proof, we begin with the analogue of Lemma 4.7.

LEMMA 5.6. *Let $\mathcal{F}$ be a VC type class of functions satisfying the same hypotheses as in Lemma 4.7 (the bracketing properties). Let $L, \phi, \tau$ be as in Theorem 5.5 and set $\gamma = 2/(1 - \tau)$. Let $0 < r < \delta < 1$ and $q \in (1, 2]$. Then, there is $C = C(K_1, K_2, q)$ such that, for all $n \in \mathbf{N}$,*

$$\Pr\left\{\sup_{\substack{r < \sigma_P f \leq \delta \\ f \in \mathcal{F}}} \frac{\nu_n(f)}{\phi(\sigma_P f)} \geq C\left[\max_{\sqrt{w/n} \leq u < \delta q} \frac{u + \sqrt{\log g_q(u)}}{L(1/u)}\right.\right.$$

$$+ \frac{e^{\gamma} w I(nr^2 < w)}{\sqrt{n} r L(1/r) \log(e^{\gamma} w/(nr^2))} \bigg]\bigg\}$$

(5.11)
$$\leq K w e^{-34 K_1 w/K}$$
$$+ K \bigg( \max_{j\,:\, n\rho_j^2 < w} N_j \bigg) \bigg( 1 + \frac{1}{2} \log_q \frac{w}{nr^2} \bigg) I(nr^2 < w) e^{-e^{\gamma} K_2^2 w/K},$$

*with notation as in Lemma* 4.7.

PROOF. We will apply Theorem 2.1'. As in Lemma 4.7, let $J_1 = \{j : n\rho_j^2 < w\}$ and $J_2 = J \setminus J_1$, where $J = \{1, \ldots, \ell\}$. Then, on $J_1$,

$$\psi_{n,j,k} \leq K_1 K_2 \frac{\rho_j}{\sqrt{n}}, \qquad \overline{V}_{n,j,k} = K_2^2 \rho_j^2,$$

and on $J_2$, by (4.12),

$$\psi_n(\rho_j) \leq K_1 \frac{\rho_j \sqrt{\log g_q(\rho_j)}}{\sqrt{n}}, \qquad \overline{V}_n(\rho_j) = 17 K_1 \rho_j^2.$$

So,

$$\sqrt{n}\,\overline{\beta}_{n,q,\phi} = \sqrt{n} \bigg[ \max_{(j,k):\,j\in J_1} \frac{\psi_{n,j,k}}{\rho_j L(1/\rho_j)} \vee \max_{j \in J_2} \frac{\psi_n(\rho_j)}{\rho_j L(1/\rho_j)} \bigg]$$

$$\leq \frac{K_1 K_2}{L(\sqrt{n/w})} \vee \max_{\sqrt{w/n} \leq \rho_j < \delta q} \frac{\sqrt{\log g_q(\rho_j)}}{L(1/\rho_j)}$$

$$\leq K_1 K_2 \max_{\sqrt{w/n} \leq u < \delta q} \frac{\sqrt{\log g_q(u)}}{L(1/u)}$$

since $\log g_q \geq 1$. For $j \in J_2$, we take $s_j = 2n \overline{V}_n(\rho_j) = 34 K_1 n \rho_j^2 \geq 34 K_1 w$, so that the contribution of $J_2$ to $\overline{\tau}_{n,j,\phi}$ is

$$2 \max_{j \in J_2} \sqrt{\frac{s_j \overline{V}_n(\rho_j)}{\rho_j^2 L^2(1/\rho_j)}} \leq 34 \sqrt{2} K_1 \max_{\sqrt{w/n} \leq u < \delta q} \frac{u}{L(1/u)},$$

whereas the contribution to the probability bound is

$$K(\mathrm{Card}\, J_2) e^{-34 K_1 w/K} \leq K \log(\delta q/r) e^{-34 K_1 w/K} \leq K w e^{-34 K_1 w/K}.$$

To keep things simple, on $J_1$ we take

$$s_j = e^{\gamma} s := e^{\gamma} K_2^2 w > e^{\gamma} K_2^2 n \rho_j^2 = e^{\gamma} \overline{V}_{n,j,k}.$$

Then, the contribution of $J_1$ to $\overline{\tau}_{n,j,\phi}$ is

$$\max_{j \in J_1} \frac{e^{\gamma} K_2^2 w}{\sqrt{n} r q^j L(1/r q^j) \log(e^{\gamma} w/(nr^2 q^{2j}))} \leq \frac{e^{\gamma} K_2^2 w I(nr^2 < w)}{\sqrt{n} r L(1/r) \log(e^{\gamma} w/(nr^2))},$$

where the inequality holds because $x^\tau L(1/x)$ is increasing and $u^{(1-\tau)/2}/\log u$ is increasing for $u > e^{2/(1-\tau)} = e^\gamma$. Finally, the contribution of $J_1$ to the probability bound is

$$K\left(\max_{j \in J_1} N_j\right)(\operatorname{Card} J_1)e^{-e^\gamma K_2^2 w/K} \leq K\left(\max_{j \in J_1} N_j\right)\left(1 + \frac{1}{2}\log\frac{w}{nr^2}\right)e^{-e^\gamma K_2^2 w/K}.$$

Now the lemma follows from collecting the above bounds and plugging them into inequality (2.4′), Theorem 2.1′. □

REMARK 5.7. The bound (5.11) can be slightly refined by taking $s_j = e^\gamma s \leq e^\gamma K_2^2 w$ for $j \in J_1$ such that $s \geq K_2^2 n\rho_j^2$ and $s_j = 34K_1\rho_j^2 n$ for the remaining $j$'s in $J_1$. Then, the contribution of $J_1$ to $\bar{\tau}$ becomes

$$\frac{e^\gamma s I(K_2^2 nr^2 < s)}{\sqrt{n}rL(1/r)(1 \vee \log(e^\gamma s/(K_2^2 nr^2)))} + 2K_2\frac{\sqrt{s}}{L(\sqrt{n/w})},$$

and the probability contribution is

$$K\left(\max_{j \in J_1} N_j\right)\left(1 + \frac{1}{2}\log\frac{w}{nr^2}\right)e^{-e^\gamma s/K}.$$

The resulting inequality is analogous to that of Lemma 4.7, whereas (5.11) is more similar to (4.15′).

Theorem 5.5 is an immediate consequence of Lemmas 5.3 and 5.6.
The following example shows how this theorem recovers the (sufficiency part of the) results in Example 2.9 of [3].

EXAMPLE 5.8 (The finite-dimensional uniform empirical c.d.f.). In the case of $\mathcal{F} = \{I_{[0,\mathbf{x}]}: 0 \leq x^i \leq 1, \prod_{i=1}^d x^i \leq 1/2\}$, $P$ being the uniform measure on $[0,1]^d$, as in Example 4.9, $g_q(\rho_j) \simeq (\log \rho_j^{-1})^{(d^2-1)/2}$ and the class satisfies the local bracketing condition. As long as $\log \log r_n^{-1}$ is of the same order as $\log \log n$, we have $w_n \simeq \log \log n$. Then, the first condition (5.12) requires $L(u) \gg (\log \log u)^{1/2}$ as $u \to \infty$. To illustrate, take $L(u) = (\log \log u)^\alpha$ for $\alpha > 1/2$. Then, the second condition (5.12) readily implies the CLT for any $r_n$ satisfying $r_n \gg \frac{(\log \log n)^{1-\alpha}}{\sqrt{n}\log \log \log n}$, which is best possible for $d > 1$ [3]. However, as in the case of Theorem 4.6, this is not sharp for $d = 1$: in this case, since in Lemmas 4.7 and 5.6 $N_j = constant$, we can take a smaller $s$ and still have the probability bound that results from Remark 5.7 tend to zero, for instance, $s = \log(w_n/nr_n^2)$. It is easy to see that if $L(u)/\log \log u \to \infty$ as $u \to \infty$, then Remark 5.7 implies the CLT for $r_n$ of a strictly smaller order than $1/\sqrt{n}$ ($r_n = \log \log \log n/\sqrt{n \log \log n}$), which by the same argument as in Example 4.10, implies that we can take $r_n = 0$, namely, one obtains the well-known Čibisov–O'Reilly CLT for the weighted uniform empirical c.d.f. in $\mathbf{R}$ [3, 12, 36].

So, Theorem 5.5, perhaps complemented by a modification along the lines of Remark 5.7, does give results comparable to those in [3] for the classical classes of sets and, moreover, it applies as well to classes of functions.

**6. Applications I: ratios of margin distributions.** The goal of this section is to suggest a much easier approach to the proofs of some of the results of Koltchinskii [25] on bounding margin distributions. The motivation and the terminology come from learning theory: functions $f$ below represent what is known as "classification margins." "Large margin algorithms" tend to output functions (classifiers) $f$ whose empirical distribution is shifted in the positive direction. The question is whether the true distribution is also shifted in the same direction. Since we are interested in the values of the margin for which these distribution functions are small, it is natural to study their ratios. See [25] for a detailed discussion.

Let
$$F_f(\delta) := P\{f \leq \delta\}, \qquad F_{n,f}(\delta) := P_n\{f \leq \delta\}.$$

Suppose that $\mathcal{F}$ is a class of functions such that
$$\forall \varepsilon > 0 \qquad \log N(\mathcal{F}; L_2(P_n); \varepsilon) \leq \left(\frac{D}{\varepsilon}\right)^\alpha$$

with some constants $D > 0$ and $\alpha \in (0, 2)$.

For two distribution functions $F, G$ and interval $(a, b)$, define
$$M_{a,b}(F; G) := \log \inf\{c > 1 : \forall t \in (a,b) : F(t) \leq cG(ct) \text{ and } G(t) \leq cF(ct)\}.$$

If $F, G$ are distribution functions on the positive real line [i.e., $F(0) = G(0) = 0$], then $M(F; G) := M_{0,+\infty}(F, G)$ is a metric (a multiplicative version of Lévy distance). We want to study the closeness of $F_{n,f}$ to $F_f$ in distances of this type uniformly in $f \in \mathcal{F}$. Unfortunately, the metric $M$ itself cannot be used even in the case of a single function $f$ (the range of $t$'s in the definition has to be restricted). However, define for $\lambda > 0$
$$\delta_n(f; \lambda) := \inf\{\delta \geq n^{-1} : \delta^{2\alpha/(2+\alpha)} F_f(\delta) \geq \lambda n^{-2/(2+\alpha)}\}.$$

THEOREM 6.1. *If $\lambda_n \to \infty$ as $n \to \infty$ and*
$$\sup_{f \in \mathcal{F}} P\{f \geq t\} \to 0 \qquad as\ t \to \infty,$$

*then*
$$\sup_{f \in \mathcal{F}} M_{\delta_n(f;\lambda_n), +\infty}(F_{n,f}, F_f) \to 0 \qquad as\ n \to \infty\ a.s.$$

PROOF. The proof is based on a couple of inequalities that follow from Proposition 2.8 of Section 2. Namely, it will be shown that for all $q > 1$ with some constant $c > 0$ depending only on $D$ and $q$ and with an absolute constant $K > 0$ we have

$$\Pr\{\exists f \in \mathcal{F} : \delta_n(f; D^{2\alpha/(2+\alpha)}\sigma^{-2}) \leq \delta \text{ and } F_f(\delta) \geq (1-c\sigma)^{-1}F_{n,f}((1+\sigma)\delta)\}$$
$$\leq K \frac{q^2}{q^2-1} \frac{1}{t} e^{-t/Kq^2}$$

and

$$\Pr\{\exists f \in \mathcal{F} : \delta_n(f; D^{2\alpha/(2+\alpha)}\sigma^{-2}) \leq \delta \text{ and } F_{n,f}(\delta) \geq (1+c\sigma)F_f((1+\sigma)\delta)\}$$
$$\leq K \frac{q^2}{q^2-1} \frac{1}{t} e^{-t/Kq^2}$$

for all $t > 0$, $\sigma \in (0, 1]$ and

$$\delta \leq \frac{Dn^{1/2}}{t^{(2+\alpha)/2\alpha}} =: A_n(t).$$

To this end, given $\delta > 0$, define a function $\varphi$ that is equal to 1 on $(-\infty, \delta]$, equal to 0 on $[(1+\sigma)\delta, +\infty)$ and is linear in between. Clearly, $\varphi$ is Lipschitz with constant $L = \frac{1}{\sigma\delta}$. Denote

$$r_n := \frac{1}{\sigma^{2/(2+\alpha)}} \frac{A^{\alpha/(2+\alpha)}}{n^{1/(2+\alpha)}},$$

where $A = DL$.

Then $F_f(\delta) \geq r_n^2$ iff $\delta \geq \delta_n(f; D^{2\alpha/(2+\alpha)}\sigma^{-2})$ and, hence, for $\delta \geq \delta_n(f; D^{2\alpha/(2+\alpha)}\sigma^{-2})$ we also have

$$P(\varphi \circ f)^2 \geq F_f(\delta) \geq r_n^2.$$

Define

$$\Delta_n := \sup_{P(\varphi \circ f)^2 \geq r_n^2} \left| \frac{P_n(\varphi \circ f)}{P(\varphi \circ f)} - 1 \right|.$$

Then for $\delta \geq \delta_n(f; D^{2\alpha/(2+\alpha)}\sigma^{-2})$

$$F_f(\delta) \leq P(\varphi \circ f) \leq (1-\Delta_n)^{-1} P_n(\varphi \circ f) \leq (1-\Delta_n)^{-1} F_{n,f}((1+\sigma)\delta)$$

and

$$F_{n,f}(\delta) \leq P_n(\varphi \circ f) \leq (1+\Delta_n) P(\varphi \circ f) \leq (1+\Delta_n) F_f((1+\sigma)\delta).$$

To prove the inequalities it remains to obtain a bound for $\Pr\{\Delta_n \geq c\sigma\}$, which is done using Proposition 2.8. First note that since $\varphi$ is Lipschitz with constant $L$, we have for the class $\varphi \circ \mathcal{F} := \{\varphi \circ f : f \in \mathcal{F}\}$

$$\forall \varepsilon > 0 \quad \log N(\varphi \circ \mathcal{F}; L_2(P_n); \varepsilon) \leq \log N(\mathcal{F}; L_2(P_n); \varepsilon/L)$$
$$\leq \left(\frac{DL}{\varepsilon}\right)^\alpha = \left(\frac{A}{\varepsilon}\right)^\alpha.$$

By Theorem 3.1,

$$E \sup_{P(\varphi \circ f)^2 \leq r^2} |(P_n - P)(\varphi \circ f)| \leq C\left[\frac{A^{\alpha/2}}{\sqrt{n}} r^{1-\alpha/2} \vee \frac{A^\alpha}{nr^\alpha}\right].$$

Under the assumption

$$r \geq r_n \geq \frac{A^{\alpha/(2+\alpha)}}{n^{1/(2+\alpha)}},$$

the first term dominates, so we have

$$E_n := \sup_{r \geq r_n} \frac{1}{r^2} E \sup_{P(\varphi \circ f)^2 \leq r^2} |(P_n - P)(\varphi \circ f)| \leq C \frac{A^{\alpha/2}}{\sqrt{n}} r_n^{-1-\alpha/2} = C\sigma,$$

by the definition of $r_n$. Using Proposition 2.8, we get

$$\Pr\bigg\{\Delta_n \geq qC\sigma + 2q\sqrt{\frac{t}{nr_n^2}(1 + 16C\sigma)}$$
$$\vee \frac{2qt}{nr_n^2 \log(t/(nr_n^2(1+16C\sigma)) \vee 2)}\bigg\} \leq K\frac{q^2}{q^2-1}\frac{1}{t}e^{-t/Kq^2}.$$

Now if

$$\delta \leq \frac{Dn^{1/2}}{t^{(2+\alpha)/2\alpha}}$$

then, by a simple computation, $\frac{t}{nr_n^2} \leq \sigma^2$, so we easily get with some constant $c > 0$ and for $\sigma \in (0,1]$

$$\Pr\{\Delta_n \geq c\sigma\} \leq K\frac{q^2}{q^2-1}\frac{1}{t}e^{-t/Kq^2}.$$

The inequalities now follow. We will use them for $\delta_j = q^{-j} \in [n^{-1}, A_n(t)]$ to prove that on an event $E$ with

$$\Pr(E) \geq K \log_q(nA_n(t))\frac{q^2}{q^2-1}\frac{1}{t}e^{-t/Kq^2}$$

we have

$$\forall j \, \forall \delta \in (\delta_{j+1}, \delta_j] \quad F_f(\delta) \leq F_f(\delta_j) \leq (1-c\sigma)^{-1}F_{n,f}((1+\sigma)\delta_j)$$
$$\leq (1-c\sigma)^{-1}F_{n,f}((1+\sigma)q\delta)$$

and

$$\forall j \ \forall \delta \in (\delta_{j+1}, \delta_j] \qquad F_{n,f}(\delta) \leq F_{n,f}(\delta_j) \leq (1+c\sigma)F_f((1+\sigma)\delta_j)$$
$$\leq (1+c\sigma)F_f((1+\sigma)q\delta),$$

which implies that on this event

$$\sup_{f \in \mathcal{F}} M_{\delta_n(f;D^{2\alpha/(2+\alpha)}\sigma^{-2}),A_n(t)}(F_{n,f};F_f) \leq (\sigma q + q - 1) \vee \frac{c\sigma}{1-c\sigma}.$$

Choosing $t = t_n = 2Kq^2 \log n$ and using the Borel–Cantelli lemma we get

$$\limsup_{n \to \infty} \sup_{f \in \mathcal{F}} M_{\delta_n(f;D^{2\alpha/(2+\alpha)}\sigma^{-2}),A_n(t_n)}(F_{n,f};F_f) \leq (\sigma q + q - 1) \vee \frac{c\sigma}{1-c\sigma} \qquad \text{a.s.},$$

and since $\sigma > 0$ and $q > 1$ are arbitrary and, under the condition $\lambda_n \to \infty$, for large enough $n$,

$$\delta_n(f;\lambda_n) \geq \delta_n(f;D^{2\alpha/(2+\alpha)}\sigma^{-2}),$$

we get

$$\sup_{f \in \mathcal{F}} M_{\delta_n(f;\lambda_n),A_n(t_n)}(F_{n,f};F_f) \to 0 \qquad \text{a.s.}$$

It now follows from the definitions that to prove

$$\sup_{f \in \mathcal{F}} M_{\delta_n(f;\lambda_n),+\infty}(F_{n,f};F_f) \to 0 \qquad \text{a.s.}$$

it suffices to check that

$$\sup_{f \in \mathcal{F}} M_{B_n,+\infty}(F_{n,f};F_f) \to 0 \qquad \text{a.s.}$$

for any sequence $B_n$ such that

$$\frac{A_n(t_n)}{B_n} \to \infty.$$

We have by conditions

$$\tau_n := \sup_{f \in \mathcal{F}} P\{f \geq B_n\} \to 0$$

and it also follows from (3.1) in [25] that a.s.

$$\eta_n := \sup_{f \in \mathcal{F}} P_n\{f \geq B_n\} \to 0.$$

For all $f \in \mathcal{F}$ and all $\delta \geq B_n$, we have

$$F_f(\delta) \geq 1 - \tau_n \quad \text{and} \quad F_{n,f}(\delta) \geq 1 - \eta_n.$$

Let $c > 1$. Then a.s. for all large enough $n$ (such that $\tau_n \leq 1 - c^{-1}$ and $\eta_n \leq 1 - c^{-1}$) for all $f \in \mathcal{F}$ and all $\delta \geq B_n$

$$F_f(\delta) \leq 1 \leq cF_{n,f}(\delta) \quad \text{and} \quad F_{n,f}(\delta) \leq 1 \leq cF_f(\delta),$$

implying

$$\sup_{f \in \mathcal{F}} M_{B_n, +\infty}(F_{n,f}; F_f) \leq \log c \leq c - 1,$$

and the result follows. $\square$

**7. Applications II: excess risk bounds in empirical risk minimization.** In this section, we discuss the problem of minimizing $Pf$ over the class $\mathcal{F}$ that is interpreted in learning theory as a risk minimization problem (e.g., in the regression or classification setting). Since the distribution $P$ is typically unknown, it has to be replaced by *empirical risk minimization*

$$P_n f \longrightarrow \min, \qquad f \in \mathcal{F}.$$

For simplicity, assume that $\hat{f}_n$ is a precise solution of the above problem, that is, it is an empirical risk minimizer (the results can be easily modified if it is only an approximate solution). Given $f \in \mathcal{F}$, let

$$\mathcal{E}_P(f) := Pf - \inf_{g \in \mathcal{F}} Pg.$$

This quantity is often called the *excess risk* of $f$. It is of interest to obtain bounds on the excess risk $\mathcal{E}_P(\hat{f}_n)$ of the empirical risk minimizer $\hat{f}_n$. It is also of interest to have some control of the ratios $\frac{\mathcal{E}_{P_n}(f)}{\mathcal{E}_P(f)}$ uniformly in $\mathcal{F}$.

The bounds given below are modifications of recent results of Koltchinskii [26]. Let

$$\mathcal{F}(\delta) := \{f \in \mathcal{F} : \mathcal{E}_P(f) \leq \delta\}$$

be the $\delta$-minimal set of $P$. For

$$\rho_P^2(f, g) \geq P(f - g)^2 - (P(f - g))^2,$$

define the diameter of the set $\mathcal{F}(\delta)$

$$D(\delta) := D_P(\delta) := \sup_{f,g \in \mathcal{F}(\delta)} \rho_P(f, g).$$

Also define

$$\psi_n(\delta) := E \sup_{f,g \in \mathcal{F}(\delta)} |(P_n - P)(f - g)|.$$

Let

$$\beta_n(r) := \sup_{\rho \geq r} \frac{\psi_n(\rho)}{\rho}$$

and

$$\Delta(r) := \sup_{\rho \geq r} \frac{D^2(\rho)}{\rho}.$$

Finally, for $s > 0$, denote

$$\gamma_n(r,s) := \beta_n(r) + 2\sqrt{\frac{s}{nr}(\Delta(r) + 16\beta_n(r))}$$

$$\vee \frac{2s}{nr \log((s/(nr(\Delta(r) + 16\beta_n(r)))) \vee 2)}.$$

THEOREM 7.1. *There exists a constant $K > 0$ such that for $q > 1$, $s > 0$ and $r > 0$ satisfying the condition*

$$q\gamma_n(r;s) < 1,$$

*the following inequality holds:*

$$(7.1) \quad \Pr\left\{ \sup_{\substack{f \in \mathcal{F} \\ \mathcal{E}_P(f) \geq r}} \left| \frac{\mathcal{E}_{P_n}(f)}{\mathcal{E}_P(f)} - 1 \right| \geq q\gamma_n(r,s) \right\} \leq K \frac{q}{q-1} \frac{1}{s} e^{-s/Kq}.$$

*Moreover, let $\tilde{f}_n \in \mathcal{F}$ be a data-dependent function such that*

$$\mathcal{E}_{P_n}(\tilde{f}_n) \leq (1 - q\gamma_n(r;s))r.$$

*Then*

$$(7.2) \quad \Pr\{\mathcal{E}_P(\tilde{f}_n) \geq r\} \leq K \frac{q}{q-1} \frac{1}{s} e^{-s/Kq}.$$

*In particular, (7.2) holds for $\tilde{f}_n = \hat{f}_n$.*

PROOF. As before, denote

$$\rho_j := rq^j, \quad j = 1, \ldots, l,$$

with $l$ being the smallest natural number such that $\rho_l \geq 1$. Let

$$\mathcal{F}_j := \{f - g : f, g \in \mathcal{F}(\rho_j)\}.$$

The key ingredient of the proof is the following inequality:

$$(7.3) \quad \Pr\left\{ \max_{1 \leq j \leq l} \frac{\|P_n - P\|_{\mathcal{F}_j}}{\rho_j} \geq \gamma_n(r,s) \right\} \leq K \frac{q}{q-1} \frac{1}{s} e^{-s/Kq}.$$

Its proof is a straightforward modification of the proof of (2.4a) of Theorem 2.1 with further bounding as in (2.11) of Proposition 2.8 (but taking

$s_j = sq^j$), so we skip the details of the derivation. The only difference is that the bound on

$$V_n(\rho_j) := \frac{1}{n} E \left\| \sum_{i=1}^n (f(X_i) - Pf)^2 \right\|_{\mathcal{F}_j}$$

now involves the diameter of the set $\mathcal{F}(\rho_j)$:

$$V_n(\rho_j) \leq D^2(\rho_j) + 16\psi_n(\rho_j).$$

Now on the event

$$E := \left\{ \max_{1 \leq j \leq l} \frac{\|P_n - P\|_{\mathcal{F}_j}}{\rho_j} \leq \gamma_n(r,s) \right\}$$

we have for all $1 \leq j \leq l$ the following implication:

$$f \in \mathcal{F}(\rho_j) \setminus \mathcal{F}(\rho_{j-1}) \implies \forall \sigma \in (0, \rho_j) \forall g \in \mathcal{F}(\sigma)$$

$$\begin{aligned}
\mathcal{E}_P(f) &\leq P(f-g) + \sigma \\
&\leq P_n(f-g) + \sigma + \|P_n - P\|_{\mathcal{F}_j} \\
&\leq \mathcal{E}_{P_n}(f) + \sigma + \rho_j \gamma_n(r,s) \\
&\leq \mathcal{E}_{P_n}(f) + \sigma + q\mathcal{E}_P(f)\gamma_n(r,s).
\end{aligned}$$

Since $\sigma > 0$ is arbitrary, this implies that on the event $E$ for all $f \in \mathcal{F}$ with $\mathcal{E}_P(f) \geq r$,

$$\frac{\mathcal{E}_{P_n}(f)}{\mathcal{E}_P(f)} \geq 1 - q\gamma_n(r,s).$$

Since $\mathcal{E}_{P_n}(\hat{f}_n) = 0$, under the condition $1 - q\gamma_n(r,s) > 0$, we must have on the event $E$ $\mathcal{E}_P(\hat{f}_n) < r$. Therefore, we have on the event $E$ the following implication:

$$\begin{aligned}
f \in \mathcal{F}(\rho_j) \setminus \mathcal{F}(\rho_{j-1}) \implies \mathcal{E}_{P_n}(f) \\
= P_n f - P_n \hat{f}_n \leq Pf - P\hat{f}_n + \|P_n - P\|_{\mathcal{F}_j} \\
\leq \mathcal{E}_P(f) + \rho_j \gamma_n(r,s) \leq \mathcal{E}_P(f)(1 + q\gamma_n(r,s)),
\end{aligned}$$

which means that on the event $E$ for all $f \in \mathcal{F}$ with $\mathcal{E}_P(f) \geq r$

$$\frac{\mathcal{E}_{P_n}(f)}{\mathcal{E}_P(f)} \leq 1 + q\gamma_n(r,s).$$

Since by (7.3)

$$\Pr(E^c) \leq K \frac{q}{q-1} \frac{1}{s} e^{-s/Kq},$$

inequality (7.1) now follows. Inequality (7.2) is an obvious consequence of (7.1) since the assumptions

$$\mathcal{E}_{P_n}(\tilde{f}_n) \leq (1 - q\gamma_n(r;s))r$$

and $\mathcal{E}_P(\tilde{f}_n) \geq r$ lead to

$$\frac{\mathcal{E}_{P_n}(\tilde{f}_n)}{\mathcal{E}_P(\tilde{f}_n)} \leq 1 - q\gamma_n(r,s). \qquad \square$$

If we define

$$\omega_n(\delta) := E \sup_{f,g \in \mathcal{F}, \rho_P(f,g) \leq \delta} |(P_n - P)(f - g)|,$$

then

$$\psi_n(\delta) \leq \omega_n(D(\delta)).$$

As a result, under the assumptions

$$\omega_n(\delta) \leq C\delta^{1-\rho} n^{-1/2}$$

and

$$D(\delta) \leq C\delta^{1/(2\kappa)}$$

with some $C > 0, \rho \in (0,1), \kappa \geq 1$, Theorem 7.1 gives a convergence rate of $\mathcal{E}_P(\tilde{f}_n)$ to 0 of the order

$$O(n^{-\kappa/(2\kappa+\rho-1)}),$$

a very typical rate in regression and classification problems.

7.1. *Regression.* For simplicity and in order to directly use the above bounds, we consider only regression models with bounded response. Let $(X,Y)$ be a random couple taking values in $S \times [0,1]$. The regression function

$$g_0(x) := E(Y|X=x), \qquad x \in S,$$

takes its values in $[0,1]$ and minimizes the functional $g \mapsto E(Y - g(X))^2$. The problem of estimating $g_0$ becomes a risk minimization problem $Pf \longrightarrow \min$ if one defines $P$ as the distribution of $(X,Y)$ and relates to each $g$ on $S$ the function $f$ on $S \times [0,1]$ as follows:

$$f(x,y) := f_g(x,y) := (y - g(x))^2, \qquad (x,y) \in S \times [0,1].$$

Given a class $\mathcal{G}$ of measurable functions from $S$ into $[0,1]$ and a sample $(X_1,Y_1),\ldots,(X_n,Y_n)$ of i.i.d. copies of $(X,Y)$, one can define a least-squares estimate of $g_0$ as a solution $\hat{g}_n$ of the following minimization problem:

$$n^{-1} \sum_{j=1}^{n} (Y_j - g(X_j))^2 \longrightarrow \min, \qquad g \in \mathcal{G},$$

which is equivalent to minimizing $P_n f$ over the class $\mathcal{F} := \{f_g : g \in \mathcal{G}\}$, $P_n$ being the empirical measure based on the sample $(X_1,Y_1),\ldots,(X_n,Y_n)$. This will allow us to use the bounds of Theorem 7.1. First suppose that $g_0 \in \mathcal{G}$. Then, by a simple and direct computation,

$$\mathcal{E}_P(f_g) = E(Y - g(X))^2 - E(Y - g_0(X))^2 = \|g - g_0\|_{L_2(\Pi)}^2,$$

where $\Pi$ is the distribution of $X$. Therefore,

$$\mathcal{F}(\delta) = \{f \in \mathcal{F} : \mathcal{E}_P(f) \leq \delta\} = \{f_g : \|g - g_0\|_{L_2(\Pi)}^2 \leq \delta\}.$$

Also, if $g_1, g_2 \in \mathcal{G}$, then

$$\begin{aligned}P(f_{g_1} - f_{g_2})^2 &= E((Y - g_1(X))^2 - (Y - g_2(X))^2)^2 \\ &= E(g_1(X) - g_2(X))^2(2Y - g_1(X) - g_2(X))^2 \\ &\leq 4\|g_1 - g_2\|_{L_2(\Pi)}^2 =: \rho_P^2(f_{g_1}, f_{g_2}),\end{aligned}$$

since $Y, g_1(X), g_2(X) \in [0,1]$. It immediately follows that the $\rho_P$-diameter of $\mathcal{F}(\delta)$ satisfies the following bound: $D(\delta) \leq 4\sqrt{\delta}$ and as a result we have $\Delta(r) \leq 16$.

Next, the usual symmetrization inequality gives

$$\psi_n(\delta) = E \sup_{g_1, g_2 \in \mathcal{G}, \|g_1 - g_0\|_{L_2(\Pi)}^2 \leq \delta, \|g_2 - g_0\|_{L_2(\Pi)}^2 \leq \delta} |(P_n - P)(f_{g_1} - f_{g_2})|$$

$$\leq 4E \sup_{g \in \mathcal{G}, \Pi(g - g_0)^2 \leq \delta} \left| n^{-1} \sum_{i=1}^n \varepsilon_i ((Y_i - g(X_i))^2 - (Y_i - g_0(X_i))^2) \right|,$$

and, using a Rademacher comparison inequality (e.g., [29], Theorem 4.12), this can be bounded further by

$$8E \sup_{g \in \mathcal{G}, \Pi(g - g_0)^2 \leq \delta} \left| n^{-1} \sum_{i=1}^n \varepsilon_i (g(X_i) - g_0(X_i)) \right| =: \tilde{\psi}_n(\delta).$$

The inequality of Theorem 4.12 in [29] is used as follows: for fixed $X_i, Y_i$, define $A_i := (Y_i - g_0(X_i))^2$, $\phi_i(u) := (A_i - u)^2 - A_i^2$ and, using the fact that $\phi_i$ are Lipschitz functions on $[0,1]$, upper bound

$$E_\varepsilon \sup_{g \in \mathcal{G}, \Pi(g - g_0)^2 \leq \delta} \left| n^{-1} \sum_{i=1}^n \varepsilon_i \phi_i(g(X_i) - g_0(X_i)) \right|.$$

Define now

$$\tilde{\beta}_n(r) := \sup_{\rho \geq r} \frac{\tilde{\psi}_n(\rho)}{\rho}$$

and

$$\tilde{\gamma}_n(r,s) := \tilde{\beta}_n(r) + 8\sqrt{\frac{s}{nr}(1 + \tilde{\beta}_n(r))} \vee \frac{2s}{nr \log((s/(16nr(1 + \tilde{\beta}_n(r)))) \vee 2)}.$$

Theorem 7.1 immediately implies that as soon as $g_0 \in \mathcal{G}$ and $q\tilde{\gamma}_n(r,s) < 1$ we have

$$\Pr\{\|\hat{g}_n - g_0\|^2_{L_2(\Pi)} \geq r\} \leq K \frac{q}{q-1} \frac{1}{s} e^{-s/Kq}. \tag{7.4}$$

Clearly, a similar bound holds for approximate least-squares estimates (as in Theorem 7.1). It is also possible and easy to handle the case $g_0 \notin \mathcal{G}$ and to bound in this case $\|\hat{g}_n - g_0\|^2_{L_2(\Pi)}$ by

$$K\left(\inf_{g \in \mathcal{G}} \|g - g_0\|^2_{L_2(\Pi)} + r\right)$$

with high probability, but we do not give this type of bound here (see, e.g., [26]). We conclude this brief discussion of regression problems with a couple of specific examples where the expectation bounds are used to derive the value of $r$ in (7.4).

EXAMPLE 7.2. Let $\mathcal{G}$ be a VC-subgraph class of measurable functions from $S$ into $[0,1]$. Let $F_\delta : S \mapsto [0,1]$ be a measurable envelope of the class $\{g - g_0 : \Pi(g - g_0)^2 \leq \delta\}$. Denote

$$\tau(\delta) := \frac{\|F_\delta\|_{L_2(\Pi)}}{\sqrt{\delta}}.$$

Applying Theorem 3.1 to VC-subgraph classes gives

$$\tilde{\psi}_n(\delta) \leq \frac{K}{\sqrt{n}} \sqrt{\delta \log \tau(\delta)},$$

assuming $\frac{\log \tau(\delta)}{\delta} \geq n$. Therefore, we have (under a natural assumption that the function $\delta \mapsto \frac{\log \tau(\delta)}{\delta}$ is nonincreasing)

$$\tilde{\beta}_n(r) \leq \frac{K}{\sqrt{n}} \sqrt{\frac{\log \tau(r)}{r}}$$

for $r$ larger than or equal to the solution $r_n$ of the equation

$$\frac{\log \tau(r)}{r} = n.$$

Then for $r = C(r_n + \frac{s}{n})$ with large enough $C$ and for $q = 2$, we have $q\tilde{\gamma}_n(r,s) < 1$ and the following bound holds for the least-squares estimate $\hat{g}_n$:

$$\Pr\left\{\|\hat{g}_n - g_0\|^2_{L_2(\Pi)} \geq C\left(r_n + \frac{s}{n}\right)\right\} \leq K \frac{1}{s} e^{-s/K}.$$

Since $\tau(\delta) \leq \frac{1}{\sqrt{\delta}}$, this always gives the convergence rate at least as good as $O(\frac{\log n}{n})$ for least-squares estimators picked from VC-subgraph classes.

However, if $\tau(\delta)$ is smaller, one can get an improvement on the logarithmic factor. In particular, if $\mathcal{G}$ is a subset of a finite-dimensional space of functions on $S$ of dimension $d$, then one can find an orthonormal system of functions $e_1, \ldots, e_d$ in $L_2(\Pi)$ such that $\mathcal{G} \subset \text{l.s.}(e_1, \ldots, e_d)$. Then we have

$$\sup_{\|g-g_0\|^2_{L_2(\Pi)} \leq \delta} |g - g_0|(x) = \sup_{\sum_{j=1}^d (\alpha_j - \alpha_j^0)^2 \leq \delta} \left| \sum_{j=1}^d (\alpha_j - \alpha_j^0) e_j(x) \right|$$

$$\leq \sup_{\sum_{j=1}^d (\alpha_j - \alpha_j^0)^2 \leq \delta} \left( \sum_{j=1}^d (\alpha_j - \alpha_j^0)^2 \right)^{1/2} \left( \sum_{j=1}^d e_j^2(x) \right)^{1/2}$$

$$\leq \sqrt{\delta} \left( \sum_{j=1}^d e_j^2(x) \right)^{1/2}.$$

If we set

$$F_\delta(x) := \sqrt{\delta} \left( \sum_{j=1}^d e_j^2(x) \right)^{1/2} \wedge 1,$$

this implies $\|F_\delta\| \leq \sqrt{d\delta}$ and as a result $\tau(\delta) \leq \sqrt{d}$, which gives the correct convergence rate $O(n^{-1})$.

EXAMPLE 7.3. Let $\mathcal{G}$ denote the set of all monotone step functions from $[0,1]$ into itself with a finite number of jumps. For a fixed $g_0 \in \mathcal{G}$, say with $m$ jumps, the class $\{g - g_0 : g \in \mathcal{G}\}$ is VC-major ($g_0$ defines a partition of $[0,1]$ into $m$ intervals; on each of these intervals $g - g_0$ is monotone and hence $\{\{g - g_0 \geq t\} : g \in \mathcal{G}, t \in \mathbf{R}\}$ is a VC class with VC dimension depending on $m$). Arguing as in Example 3.8, we can show that

$$\tilde{\psi}_n(\delta) \leq \frac{K}{\sqrt{n}} \sqrt{\delta} \left( \log \frac{1}{\delta} \right)^{3/4} \left( \log \log \frac{1}{\delta} \right)^{1/2}$$

$$\vee \frac{K}{n} \left( \log \frac{1}{\delta} \right)^{3/2} \log \log \frac{1}{\delta} \vee \frac{\sqrt{\log n}}{n},$$

which implies

$$\tilde{\beta}_n(r) \leq \frac{K}{\sqrt{nr}} \left( \log \frac{1}{r} \right)^{3/4} \left( \log \log \frac{1}{r} \right)^{1/2}$$

$$\vee \frac{K}{nr} \left( \log \frac{1}{r} \right)^{3/2} \log \log \frac{1}{r} \vee \frac{\sqrt{\log n}}{nr}.$$

Let us take $q = 2$. Then it is easy to conclude that $q\tilde{\gamma}_n(r,s) < 1$ as soon as

$$r \geq C \frac{s + (\log n)^{3/2} \log \log n}{n}$$

with sufficiently large constant $C$ (which will depend on the number of jumps of $g_0$!). Hence, if we take an estimate $\tilde{g}_n$ such that

$$n^{-1}\sum_{j=1}^{n}(Y_j - \tilde{g}_n(X_j))^2 \leq \inf_{g\in\mathcal{G}} n^{-1}\sum_{j=1}^{n}(Y_j - g(X_j))^2 + \frac{(\log n)^{3/2}\log\log n}{2n},$$

then Theorem 7.1 implies that

$$\Pr\left\{\|\tilde{g}_n - g_0\|_{L_2(\Pi)}^2 \geq C(g_0)\frac{s + (\log n)^{3/2}\log\log n}{n}\right\} \leq K\frac{1}{s}e^{-s/K},$$

with some constants $C(g_0)$ and $K$. In particular, the bound implies that

$$E\|\tilde{g}_n - g_0\|_{L_2(\Pi)}^2 = O\left(\frac{(\log n)^{3/2}\log\log n}{n}\right).$$

Since the constant $C(g_0)$ tends to infinity as the number of jumps of the function $g_0$ tends to infinity, the above bound cannot be made uniform in $g_0 \in \mathcal{G}$ (and, in fact, the convergence rate of $\sup_{g_0\in\mathcal{G}} E\|\tilde{g}_n - g_0\|_{L_2(\Pi)}^2$ to 0 is much slower for *any* estimator $\tilde{g}_n$). Results of this type (in a slightly different context and with an improvement on the logarithmic factors) can be found, for instance, in [45] and references therein.

7.2. *Classification.* In classification problems, one deals with random couples $(X,Y)$ in $S \times \{0,1\}$, $X$ being an observable instance and $Y$ an unobservable binary label assigned to this instance. Functions $g$ from $S$ into $\{0,1\}$ are called classifiers. The generalization error of a classifier $g$ is defined as

$$\Pr\{Y \neq g(X)\} = P\{(x,y) : y \neq g(x)\},$$

where $P$ is the joint distribution of $(X,Y)$. It is well known that the minimal possible generalization error (the Bayes risk) is attained at a classifier

$$g_0(x) := I(\eta(x) \geq 1/2),$$

where $\eta(x) := E(Y|X=x)$ is the regression function. Since the distribution $P$ of $(X,Y)$ and hence also the regression function $\eta$ are unknown, a reasonable approach to classification is to minimize the training error

$$n^{-1}\sum_{j=1}^{n} I(Y_j \neq g(X_j)) = P_n\{(x,y) : y \neq g(x)\},$$

based on i.i.d. training examples sampled from $P$, over a suitable class $\mathcal{G}$ of classifiers. For simplicity, we assume that $g_0 \in \mathcal{G}$. Denote $\hat{g}_n$ a minimizer of the training error over the class $\mathcal{G}$. Thus, the classification problem becomes a version of empirical risk minimization and one can use Theorem 7.1 to study the size of the excess risk

$$\mathcal{E}(\hat{g}_n) := P\{(x,y) : y \neq \hat{g}_n(x)\} - P\{(x,y) : y \neq g_0(x)\}.$$

As before, $\Pi$ denotes the distribution of $X$. If, with some $\kappa \geq 1$ and $c > 0$, for all $g \in \mathcal{G}$

$$P\{(x,y) : y \neq g(x)\} - P\{(x,y) : y \neq g_0(x)\} \geq c\Pi^\kappa \{x : g(x) \neq g_0(x)\},$$

then the diameter $D(\delta) \leq C\delta^{1/(2\kappa)}$. This holds, for instance, if for all $t > 0$

$$\Pi\{x : 0 < |\eta(x) - 1/2| \leq t\} \leq Ct^\alpha$$

and in this case $\kappa = \frac{1+\alpha}{\alpha}$ [44]. Under the standard condition that the $\varepsilon$-entropy of the class $\mathcal{G}$ grows as $O(\varepsilon^{-2\rho})$ [with several possible kinds of entropy involved and with $\rho \in (0,1)$] Theorem 7.1 yields a bound on the excess risk of the order $O(n^{-\kappa/(2\kappa+\rho-1)})$ as in [44]. The main difference with the $L_2$-regression problem where $\kappa = 1$ is that in classification $\kappa$ can take any value greater than or equal to 1 leading to the whole spectrum of convergence rates. If there exists $h > 0$ such that

$$\forall x \in S \qquad |\eta(x) - 1/2| \geq h,$$

then it is easy to see that

$$P\{(x,y) : y \neq g(x)\} - P\{(x,y) : y \neq g_0(x)\} \geq ch\Pi\{x : g(x) \neq g_0(x)\},$$

so we do have $\kappa = 1$. This case of well-separated classes was looked at in the recent paper of Massart and Nedelec [34]. Define $f_g(x,y) := I(y \neq g(x))$ and $\mathcal{F} := \{f_g : g \in \mathcal{G}\}$. We are using the distance

$$\rho_P(f_{g_1}, f_{g_2}) := \Pi^{1/2}(g_1 - g_2)^2.$$

Then we have the following bound for the diameter $D(\delta)$:

$$D(\delta) \leq C\left(\frac{\delta}{h}\right)^{1/2}$$

implying

$$\Delta(r) \leq \frac{C}{h}.$$

Suppose that $\mathcal{C} := \{\{g = 1\} : g \in \mathcal{G}\}$ is a VC class and $C_0 := \{g_0 = 1\}$. Define a local version of Alexander's capacity function:

$$\tau(\delta) := \frac{\Pi(\bigcup_{C \in \mathcal{C}, \Pi(C \triangle C_0) \leq \delta}(C \triangle C_0))}{\delta}.$$

Then

$$\psi_n(\delta) \leq K\sqrt{\frac{V}{nh}\delta \log \tau\left(\frac{\delta}{h}\right)}$$

and as a result

$$\beta_n(r) \leq K\sqrt{\frac{V}{nhr} \log \tau\left(\frac{r}{h}\right)}.$$

To satisfy the condition $q\gamma_n(r,s) < 1$ (say, with $q=2$) it is enough to take

$$r = C\left[h\varphi\left(\frac{nh^2}{V}\right) + \frac{s}{nh}\right],$$

where $\varphi$ denotes the inverse of the function

$$r \mapsto \frac{\log \tau(r)}{r}$$

and $C$ is a sufficiently large constant. Now it is easy to check that

$$\varphi\left(\frac{nh^2}{V}\right) \leq \frac{V}{nh^2}\log \tau\left(\frac{V}{nh^2}\right),$$

yielding the following bound:

$$\Pr\left\{\mathcal{E}(\hat{g}_n) \geq C\left[\frac{V}{nh}\log \tau\left(\frac{V}{nh^2}\right) + \frac{s}{nh}\right]\right\} \leq K\frac{1}{s}e^{-s/K}.$$

If we replace $\tau(r)$ by the trivial upper bound $\frac{1}{r}$, this gives one of the results of Massart and Nedelec [34]: the excess risk is bounded by

$$K\frac{V}{nh}\log \frac{nh^2}{V}.$$

In the case of smaller $\tau$, it is a slight improvement of their bound. It is easy to see that for some classes of sets and probability measures $P$ $\tau$ can be even bounded, leading to the bound on excess risk of the order $O(\frac{1}{nh})$. For instance, as in Section 4, suppose that $S = [0,1]^d$ for some $d \geq 1$ and $P$ has a density that is uniformly bounded and bounded away from $0$ on $S$. As before, $h$ is the Hausdorff distance between subsets of $S$. Let $\mathcal{C}$ be a VC class of convex subsets of $S$ and $C_0 \in \mathcal{C}$, $P(C_0) > 0$. Suppose that for some $K > 0$

$$K^{-1}h(C, C_0) \leq P(C \triangle C_0) \leq Kh(C, C_0), \qquad C \in \mathcal{C}.$$

Recall that the upper bound always holds for convex sets ([17], pages 269–270), but the lower bound holds only for special classes of sets (balls, rectangles, etc.). Then the function $\tau$ is uniformly bounded. The proof easily follows from the same type of argument as in Section 4 (before Example 4.9).

**Acknowledgment.** We thank an anonymous referee for several suggestions and references.

## REFERENCES

[1] ALEXANDER, K. S. (1985). Rates of growth for weighted empirical processes. In *Proc. of the Berkeley Conference in Honor of Jerzy Neyman and Jack Kiefer* **II** (L. Le Cam and R. Olshen, eds.) 475–493. Wadsworth, Belmont, CA. MR0822047

[2] ALEXANDER, K. S. (1987). Rates of growth and sample moduli for weighted empirical processes indexed by sets. *Probab. Theory Related Fields* **75** 379–423. MR0890285

[3] ALEXANDER, K. S. (1987). The central limit theorem for weighted empirical processes indexed by sets. *J. Multivariate Anal.* **22** 313–339. MR0899666

[4] ALEXANDER, K. S. (1987). The central limit theorem for empirical processes on Vapnik–Červonenkis classes. *Ann. Probab.* **15** 178–203. MR0877597

[5] BARTLETT, P., BOUSQUET, O. and MENDELSON, S. (2002). Localized Rademacher complexities. *Computational Learning Theory. Lecture Notes in Comput. Sci.* **2375** 44–58. Springer, Berlin. MR2040404

[6] BARTLETT, P. and LUGOSI, G. (1999). An inequality for uniform deviations of sample averages from their means. *Statist. Probab. Lett.* **44** 55–62. MR1706315

[7] BARTLETT, P. and MENDELSON, S. (2003). Empirical risk minimization. Preprint.

[8] BIRGÉ, L. and MASSART, P. (1998). Rates of convergence for minimum contrast estimators. *Bernoulli* **4** 329–375. MR1653272

[9] BOUSQUET, O. (2002). Concentration inequalities and empirical processes theory applied to the analysis of learning algorithms. Ph.D. thesis, Ecole Polytechnique, Paris.

[10] BOUSQUET, O. (2003). Concentration inequalities for sub-additive functions using the entropy method. In *Stochastic Inequalities and Applications* 213–247. Birkhäuser, Basel. MR2073435

[11] BOUSQUET, O., KOLTCHINSKII, V. and PANCHENKO, D. (2002). Some local measures of complexity of convex hulls and generalization bounds. *Computational Learning Theory. Lecture Notes in Comput. Sci.* **2375** 59–73. Springer, Berlin. MR2040405

[12] ČIBISOV, D. M. (1964). Some limit theorems on the limiting behavior of empirical distribution functions. *Selected Translations Math. Statist. Probab.* **6** 147–156.

[13] CSÁKI, E. (1977). The law of hte iterated logarithm for normalized empirical distribution functions. *Z. Wahrsch. Verw. Gebiete* **38** 147–167. MR0431350

[14] CSÖRGŐ, M., CSÖRGŐ, S., HORVÁTH, L. and MASON, D. (1986). Weighted empirical and quantile processes. *Ann. Probab.* **14** 31–85. MR0815960

[15] DE LA PEÑA, V. and GINÉ, E. (1999). *Decoupling*: *From Dependence to Independence*. Springer, New York. MR1666908

[16] DUDLEY, R. M. (1987). Universal Donsker classes and metric entropy. *Ann. Probab.* **20** 1968–1982. MR1188050

[17] DUDLEY, R. M. (1999). *Uniform Central Limit Theorems*. Cambridge Univ. Press. MR1720712

[18] EICKER, F. (1979). The asymptotic distribution of suprema of the standardized empirical process. *Ann. Statist.* **7** 116–138. MR0515688

[19] EINMAHL, J. H. J. (1996). Extension to higher dimensions of the Jaeschke–Eicker result on the standardized empirical process. *Comm. Statist. Theory Methods* **25** 813–822. MR1380620

[20] Einmahl, U. and Mason, D. (2000). An empirical process approach to the uniform consistency of kernel type function estimators. *J. Theor. Probab.* **13** 1–37. MR1744994

[21] Giné, E. and Guillou, A. (2001). On consistency of kernel density estimators for randomly censored data: Rates holding uniformly over adaptive intervals. *Ann. Inst. H. Poincaré Probab. Statist.* **37** 503–522. MR1876841

[22] Giné, E., Koltchinskii, V. and Wellner, J. (2003). Ratio limit theorems for empirical processes. In *Stochastic Inequalities and Applications* 249–278. Birkhäuser, Basel. MR2073436

[23] Giné, E., Latała, R. and Zinn, J. (2000). Exponential and moment inequalities for U-statistics. In *High Dimensional Probability* **II** (E. Giné, D. M. Mason and J. Wellner, eds.) 13–38. Birkhäuser, Boston. MR1857312

[24] Jaeschke, D. (1979). The asymptotic distribution of the supremum of the standardized empirical distribution function on subintervals. *Ann. Statist.* **7** 108–115. MR0515687

[25] Koltchinskii, V. (2003). Bounds on margin distributions in learning problems. *Ann. Inst. H. Poincaré Probab. Statist.* **39** 943–978. MR2010392

[26] Koltchinskii, V. (2006). Local Rademacher complexities and oracle inequalities in risk minimization. *Ann. Statist.* To appear.

[27] Koltchinskii, V. and Panchenko, D. (2000). Rademacher processes and bounding the risk of function learning. In *High Dimensional Probability* **II** (E. Giné, D. Mason and J. Wellner, eds.) 443–459. Birkhäuser, Boston. MR1857339

[28] Koltchinskii, V. and Panchenko, D. (2002). Empirical margin distributions and bounding the generalization error of combined classifiers. *Ann. Statist.* **30** 1–50. MR1892654

[29] Ledoux, M. and Talagrand, M. (1991). *Probability in Banach Spaces.* Springer, Berlin. MR1102015

[30] Mason, D., Shorack, G. R. and Wellner, J. (1983). Strong limit theorems for oscillation moduli of the uniform empirical process. *Z. Wahrsch. Verw. Gebiete* **65** 83–97. MR0717935

[31] Massart, P. (2000). Some applications of concentration inequalities in statistics. *Ann. Fac. Sci. Tolouse Math. (6)* **9** 245–303. MR1813803

[32] Massart, P. (2000). About the constants in Talagrand's concentration inequalities for empirical processes. *Ann. Probab.* **28** 863–884. MR1782276

[33] Massart, P. (2005). Concentration inequalities with applications to model selection and statistical learning. *Lecture on Probability Theory and Statistics. Ecole d'Eté de Probabilités de Saint-Flour XXXIV-2003. Lecture Notes in Math.* To appear. Available at www.math.u-psud.fr/~massart/.

[34] Massart, P. and Nedelec, E. (2006). Risk bounds for statistical learning. *Ann. Statist.* **34**(5).

[35] Mendelson, S. (2001). Learning relatively small classes. *Computational Learning Theory. Lecture Notes in Comput. Sci.* **2111** 273–288. Springer, Berlin. MR2042041

[36] O'Reilly, N. E. (1974). On the weak convergence of empirical processes in sup-norm metrics. *Ann. Probab.* **2** 642–651. MR0383486

[37] Panchenko, D. (2002). Concentration inequalities in product spaces and applications to statistical learning theory. Ph.D. dissertation, Univ. New Mexico, Albuquerque.

[38] Panchenko, D. (2002). Some extensions of an inequality of Vapnik and Chervonenkis. *Electron. Comm. Probab.* **7** 55–65. MR1887174

[39] PANCHENKO, D. (2003). Symmetrization approach to concentration inequalities for empirical processes. *Ann. Probab.* **31** 2068–2081. MR2016612
[40] POLLARD, D. (1984). *Convergence of Stochastic Processes*. Springer, New York. MR0762984
[41] SHORACK, G. R. and WELLNER, J. A. (1982). Limit theorems and inequalities for the uniform empirical process indexed by intervals. *Ann. Probab.* **10** 639–652. MR0659534
[42] TALAGRAND, M. (1994). Sharper bounds for Gaussian and empirical processes. *Ann. Probab.* **22** 28–76. MR1258865
[43] TALAGRAND, M. (1996). New concentration inequalities in product spaces. *Invent. Math.* **126** 505–563. MR1419006
[44] TSYBAKOV, A. (2004). Optimal aggregation of classifiers in statistical learning. *Ann. Statist.* **32** 135–166. MR2051002
[45] VAN DE GEER, S. A. (1993). Hellinger-consistency of certain nonparametric maximum likelihood estimators. *Ann. Statist.* **21** 14–44. MR1212164
[46] VAN DE GEER, S. A. (2000). *Applications of Empirical Process Theory*. Cambridge Univ. Press. MR1739079
[47] VAN DER VAART, A. W. and WELLNER, J. A. (1996). *Weak Convergence and Empirical Processes*. Springer, New York. MR1385671
[48] WELLNER, J. A. (1978). Limit theorems for the ratio of the empirical distribution function to the true distribution function. *Z. Wahrsch. Verw. Gebiete* **45** 108–123. MR0651392
[49] YUKICH, J. E. (1987). Some limit theorems for empirical processes indexed by functions. *Probab. Theory Related Fields* **74** 71–90. MR0863719

DEPARTMENT OF MATHEMATICS
UNIVERSITY OF CONNECTICUT
STORRS, CONNECTICUT 06239
USA
E-MAIL: gine@math.uconn.edu

DEPARTMENT OF MATHEMATICS AND STATISTICS
UNIVERSITY OF NEW MEXICO
ALBUQUERQUE, NEW MEXICO 87131
USA
E-MAIL: vlad@math.unm.edu